\documentclass[11pt]{article}
\usepackage[a4paper,margin=1in]{geometry}
\usepackage{amsmath,amssymb,amsthm,mathtools}
\usepackage{enumitem}
\usepackage[hidelinks,
 pdftitle={Profile decomposition, strict binding and global compactness for a critical Grushin embedding},
 pdfauthor={Shengbing Deng and Dehong Yang},
 pdfkeywords={Baouendi--Grushin operator; profile decomposition; critical embedding; concentration compactness; sharp constants}]{hyperref}
\usepackage{cite}

\numberwithin{equation}{section}

\newtheorem{theorem}{Theorem}[section]
\newtheorem{proposition}[theorem]{Proposition}
\newtheorem{lemma}[theorem]{Lemma}
\newtheorem{corollary}[theorem]{Corollary}

\newtheorem{remark}[theorem]{Remark}

\title{Profile decomposition, strict binding and global compactness
for a critical Grushin embedding}
\author{Shengbing Deng\thanks{Corresponding author. Email: \texttt{shbdeng@swu.edu.cn}.}\ \ and\ 
Dehong Yang}
\date{}

\begin{document}
\maketitle

\begin{abstract}
We study attainment of the best constant in
\begin{equation*}
 \int_{\mathbb R^{m+n}}
 \bigl(|\nabla_xu|^2+|x|^2|\nabla_yu|^2+a|u|^2\bigr)
 \geq S_a\left(\int_{\mathbb R^{m+n}}
 (1+|x|^2)^{n/(Q-2)}|u|^q\right)^{2/q},
\end{equation*}
where \(m,n\geq1\), \(Q=m+2n\), \(q=2Q/(Q-2)\), and \(a>0\).
Using annular interpolation and a transverse-average recurrence with
\(L^q\) contraction \(2^{-m/q}\), we decompose bounded sequences into
finite, Grushin and Euclidean profiles.  Energy and critical mass split,
with a power bound in the number of profiles on the remaining mass.
We prove compactness of normalized minimizing sequences modulo
\(y\)-translations under
\(S_a<\min\{S_G,a^{n/Q}S_E(1)\}\), where \(S_G\) and \(S_E(1)\)
are the tangent constants.  The strict-binding set is empty or an open
interval; we prove nonemptiness for \(n=1\) and even \(m\geq6\).
At finite positive endpoints of a nonempty binding interval, we
characterize attainment by first-order binding gaps.
We obtain a finite critical-profile decomposition of Palais--Smale
sequences with a remainder vanishing in energy.
\end{abstract}

\noindent\textbf{Keywords.} Baouendi--Grushin operator; profile
decomposition; critical embedding; concentration compactness; sharp
constants.

\medskip
\noindent\textbf{2020 Mathematics Subject Classification.} 35H20, 35J70,
35B33, 46E35.
\section{Introduction}

Let \(m,n\geq1\), \(N=m+n\), and \(z=(x,y)\in\mathbb R^m\times
\mathbb R^n\).  Set
\begin{equation*}
 Q=m+2n,\qquad q=\frac{2Q}{Q-2},\qquad
 \nabla_Gu=(\nabla_xu,|x|\nabla_yu).
\end{equation*}
For \(a>0\), define
\begin{align*}
 \mathcal E_a(u)
 &=\int_{\mathbb R^N}
 \bigl(|\nabla_xu|^2+|x|^2|\nabla_yu|^2+a|u|^2\bigr)\,dz,\\
 \mathcal M(u)
 &=\int_{\mathbb R^N}(1+|x|^2)^{n/(Q-2)}|u|^q\,dz.
\end{align*}
Let \(H_a\) be the completion of \(C_c^\infty(\mathbb R^N)\) for
\(\|u\|_{H_a}^2=\mathcal E_a(u)\), and put
\begin{equation*}
 S_a=\inf_{u\in H_a\setminus\{0\}}
 \frac{\mathcal E_a(u)}{\mathcal M(u)^{2/q}}.
\end{equation*}
All functions are real-valued.  The spaces \(H_a\), \(a>0\), have
the same underlying space and equivalent norms.
We study attainment of \(S_a\) and compactness of normalized minimizing
sequences modulo translations in \(y\).

The concentration--compactness principle of Lions
\cite{Lions1984a,Lions1984b,Lions1985a,Lions1985b} separates
concentration from vanishing and splitting.
The decompositions of Solimini \cite{Solimini1995} and G\'erard
\cite{Gerard1998} identify the individual profiles.
Refined Sobolev inequalities and cocompactness provide the local
estimates needed for extraction
\cite{TintarevFieseler2007,Tintarev2013,PalatucciPisante2014}.
The homogeneous and inhomogeneous Sobolev cases are treated in
\cite{Okumura2022a,Okumura2022b}; noncompact manifolds are considered
in \cite{SandeepTintarev2019}.
For this quotient, concentration near \(x=0\) and escape through
\(|x|\to\infty\) lead to different model spaces.

Near \(\{x=0\}\), the dilation
\(\delta_r(x,y)=(rx,r^2y)\) gives the homogeneous Grushin energy
\begin{equation*}
 \mathcal E_G(\phi)=\int_{\mathbb R^N}
 \bigl(|\nabla_x\phi|^2+|x|^2|\nabla_y\phi|^2\bigr)\,dz.
\end{equation*}
Let \(D_G^{1,2}\) be its energy completion and set
\begin{equation*}
 S_G=\inf_{\phi\in D_G^{1,2}\setminus\{0\}}
 \frac{\mathcal E_G(\phi)}{\|\phi\|_q^2}.
\end{equation*}
The underlying geometry goes back to Baouendi \cite{Baouendi1967}.
Local Poincar\'e and Sobolev estimates were developed in
\cite{FranchiLanconelli1983,Jerison1986,FranchiGutierrezWheeden1994}.
Sharp homogeneous inequalities and critical equations are studied in
\cite{Beckner2001,Monti2006,MontiMorbidelli2006,DouSunWangZhu2022}.
For critical \(p\)-Grushin problems, bifurcation and multiplicity
results are given in \cite{MalanchiniMolicaBisciSecchi2025}.
Attainment of \(S_G\) and the two-sided decay of its positive
extremals follow from
\cite[Theorem~1.1]{GandalLoiudiceTyagi2025}.
The energy can also be written using the smooth fields
\(X_i=\partial_{x_i}\), \(Y_{ij}=x_i\partial_{y_j}\):
\begin{equation*}
 [X_i,Y_{ij}]=\partial_{y_j},\qquad
 \sum_i|X_iu|^2+\sum_{i,j}|Y_{ij}u|^2=|\nabla_Gu|^2.
\end{equation*}
Their maximal pointwise homogeneous dimension is \(Q\), attained on
\(\{x=0\}\).  The homogeneous H\"ormander theory in
\cite[Theorem~1.4]{ChenChenLi2026} also applies to this sharp inequality.

At infinity, the coordinates
\((x,y)=(X_k+\xi,b_k+|X_k|\eta)\), \(|X_k|\to\infty\),
give the Euclidean model
\begin{equation*}
 S_E(a)=\inf_{\psi\in H^1(\mathbb R^N)\setminus\{0\}}
 \frac{\displaystyle\int_{\mathbb R^N}
 (|\nabla\psi|^2+a|\psi|^2)\,d\zeta}{\|\psi\|_q^2}.
\end{equation*}
Here \(Q>N\), so \(q\) is Euclidean subcritical and \(S_E(a)\) is
attained.  Its scaling law is
\begin{equation*}
 S_E(a)=a^{n/Q}S_E(1).
\end{equation*}
The mass weight balances both limits.  With
\(\omega(x)=(1+|x|^2)^{n/2}\), its density is \(\omega^{q/Q}\), and
\begin{equation*}
 \omega(rx)\longrightarrow1\quad(r\downarrow0),\qquad
 R^{-n}\omega(Re+\xi)\longrightarrow1\quad(R\to\infty)
\end{equation*}
locally uniformly for \(e\in\mathbb S^{m-1}\).
Thus a bounded sequence has three possible profile spaces:
\(H_a\) after a \(y\)-translation, \(D_G^{1,2}\) near the degeneracy
set, and \(H^1(\mathbb R^N)\) at spatial infinity.

For critical equations, Struwe's theorem \cite{Struwe1984} identifies
the critical profiles of Palais--Smale sequences.
Bhakta and Sandeep \cite{BhaktaSandeep2012} prove a corresponding
decomposition in hyperbolic space and relate its equation to
cylindrically symmetric Grushin problems.
On the Heisenberg group, Sobolev profile decomposition is established
in \cite{Benameur2008}, and Struwe-type global compactness in
\cite{PalatucciPiccininiTemperini2025}.
The Morrey-type improved inequality in
\cite{BasakChakrabortyRanaRoychowdhury2026} yields a Palais--Smale
decomposition for perturbed CR Yamabe equations.
We decompose arbitrary bounded sequences for the weighted embedding
and use the two tangent constants to determine the compactness threshold.

We use \(F,G,E\) for the finite, singular and escape types.
Section 4 defines their frames and cutoff realizations.
Compactness modulo translations means subsequential strong convergence
after a sequence of translations in \(y\).

\begin{theorem}
\label{thm:profile-intro}
Let \((u_k)\) be bounded in \(H_a\).  After passage to a subsequence, there
are \(J_*\in\mathbb N_0\cup\{\infty\}\), nonzero profiles \(p^j\), pairwise
orthogonal frames \(P_k^j\), and types \(\tau_j\in\{F,G,E\}\), indexed by
\(1\leq j\leq J_*\).  Profiles of types \(F,G,E\) belong to
\(H_a,D_G^{1,2},H^1(\mathbb R^N)\), respectively.
There are synchronized realizations \(R_k^jp^j\)
such that, for every \(J\geq1\), writing \(J^\sharp=J\) when
\(J_*=\infty\) and \(J^\sharp=\min\{J,J_*\}\) otherwise,
\begin{equation*}
 u_k=\sum_{j=1}^{J^\sharp} R_k^jp^j+w_k^J.
\end{equation*}
An empty sum is zero; after finite stopping the remainder is fixed.
If \(\mathcal E_{\tau_j}\) and \(\mathcal M_{\tau_j}\) denote the
corresponding model functionals, then
\begin{align}
 \mathcal E_a(u_k)
 &=\sum_{j=1}^{J^\sharp}\mathcal E_{\tau_j}(p^j)
 +\mathcal E_a(w_k^J)+o_k(1),
 \label{eq:energy-splitting-intro}\\
 \mathcal M(u_k)
 &=\sum_{j=1}^{J^\sharp}\mathcal M_{\tau_j}(p^j)
 +\mathcal M(w_k^J)+o_k(1).
 \label{eq:mass-splitting-intro}
\end{align}
Moreover,
\begin{equation}
 \sum_{j=1}^{J_*}\mathcal E_{\tau_j}(p^j)
 \leq\limsup_{k\to\infty}\mathcal E_a(u_k),\qquad
 \lim_{J\to\infty}\limsup_{k\to\infty}\mathcal M(w_k^J)=0.
 \label{eq:profile-summability-intro}
\end{equation}
For every fixed \(J\), the remainder is weakly zero in \(H_a\) and in every
one of the first \(J^\sharp\) selected frames.  In an escape frame the latter means
weak convergence in \(H^1(B_R)\) for every fixed \(R>0\).

Let \(L=\sup_k\mathcal E_a(u_k)\), and put
\begin{equation*}
 \eta_J=\limsup_{k\to\infty}
 \max\{\Theta_F(w_k^J),\Theta_G(w_k^J),\Theta_E(w_k^J)\}.
\end{equation*}
There are \(C>0\) and \(\gamma>0\), depending only on \(m,n,a\), such that
\begin{equation}
 \eta_J\leq\frac{CL}{J+1},\qquad
 \limsup_{k\to\infty}\mathcal M(w_k^J)
 \leq CL^{q/2}(J+1)^{-\gamma}.
 \label{eq:grushin-remainder-rate-intro}
\end{equation}
For finite extraction, the terminal gauges and mass vanish, and these
bounds hold for every \(J\geq J_*\).
\end{theorem}

We normalize minimizing sequences by \(\mathcal M(u_k)=1\).

\begin{theorem}
\label{thm:threshold-intro}
For every \(a>0\),
\begin{equation*}
 0<S_a\leq\min\{S_G,S_E(a)\}.
\end{equation*}
If
\begin{equation}
 S_a<\min\{S_G,S_E(a)\},
 \label{eq:strict-binding}
\end{equation}
then \(S_a\) is attained.  Every normalized minimizing sequence is strongly
compact in \(H_a\), after a sequence of \(y\)-translations.

If \eqref{eq:strict-binding} fails, every normalized minimizing sequence has,
after a subsequence, exactly one of the following forms:
\begin{enumerate}[label=\textup{\(\roman*\)}]
\item it is compact modulo \(y\)-translations and \(S_a\) is attained;
\item \(S_a=S_G\), and it has one singular profile of full mass;
\item \(S_a=S_E(a)\), and it has one escape profile of full mass.
\end{enumerate}
In alternative \textup{(ii)}, let \(\Phi\) be the singular profile.
It is a normalized extremal for \(S_G\), and one realization satisfies
\begin{equation*}
 \mathcal E_a\bigl(u_k-D_{r_k,c_k}(\chi_{L_k}\Phi)\bigr)
 \longrightarrow0.
\end{equation*}
In alternative \textup{(iii)}, the escape frame can be chosen so that its
profile is \(\sigma Z_a\) for some \(\sigma\in\{-1,1\}\), and
\begin{equation*}
 \left\|A_{X_k,b_k}^{-1}u_k-\sigma Z_a\right\|_{X_k,a}
 \longrightarrow0.
\end{equation*}
The last convergence is global in the moving escape norm.
\end{theorem}

For the parameter dependence, set
\begin{equation*}
 a_0=\left(\frac{S_G}{S_E(1)}\right)^{Q/n},\qquad
 \mathcal B=\{a>0:S_a<\min\{S_G,S_E(a)\}\}.
\end{equation*}

\begin{theorem}
\label{thm:connected-binding-intro}
The function \(a\mapsto S_a/a^{n/Q}\) is nonincreasing on
\((0,\infty)\).  The set \(\mathcal B\) is empty or is an open interval
\begin{equation*}
 \mathcal B=(a_-,a_+),\qquad 0\leq a_-<a_0<a_+\leq\infty.
\end{equation*}
Assume \(\mathcal B\ne\varnothing\) in the endpoint assertions below.
If \(a_->0\), then \(S_a=S_E(a)\) for \(0<a\leq a_-\).
If \(a_+<\infty\), then \(S_a=S_G\) for \(a\geq a_+\).
At finite positive endpoints, put
\begin{equation*}
 \ell_-=S_E'(a_-)-(S_a)'_+(a_-),\qquad
 \ell_+=(S_a)'_-(a_+).
\end{equation*}
The binding gaps satisfy
\begin{align}
 S_E(a)-S_a&=\ell_-(a-a_-)+o(a-a_-),&&a\downarrow a_-,
 \label{eq:left-endpoint-asymptotic}\\
 S_G-S_a&=\ell_+(a_+-a)+o(a_+-a),&&a\uparrow a_+.
 \label{eq:right-endpoint-asymptotic}
\end{align}
Moreover,
\begin{equation*}
 \ell_->0\ \Longleftrightarrow\ S_{a_-}\text{ is attained},\qquad
 \ell_+>0\ \Longleftrightarrow\ S_{a_+}\text{ is attained}.
\end{equation*}
If an endpoint is not attained, normalized interior minimizers approaching
it have, after a subsequence, a single escape profile at \(a_-\), or a single singular profile
at \(a_+\), with a strong remainder in \(H_{a_-}\) or \(H_{a_+}\).
\end{theorem}

Corollary \ref{cor:exterior-nonattainment} gives the exterior branches.
An explicit binding interval is available in the following dimensions.

\begin{theorem}
\label{thm:explicit-binding}
Let \(n=1\), and let \(m\geq6\) be even.  Define
\begin{equation*}
 \underline a_m=
 \frac{m^{2m+4}}{4(m+2)^{m+2}(m-1)^{m+1}},\qquad
 \overline a_m=\frac{m(m-2)}{4(m+2)}.
\end{equation*}
Then
\begin{equation}
 a_0\leq\underline a_m<\overline a_m,\qquad
 (\underline a_m,\overline a_m)\subset\mathcal B.
 \label{eq:explicit-binding-interval}
\end{equation}
In particular, when \((m,n)=(6,1)\),
\begin{equation*}
 \left(\frac{43046721}{80000000},\frac34\right)\subset\mathcal B.
\end{equation*}
Thus \(S_{2/3}\) is attained in this model, and every normalized minimizing
sequence is compact modulo \(y\)-translations.
\end{theorem}

For general \(m,n\), Proposition \ref{thm:binding-curvature-intro}
gives second-order sufficient tests.  The singular test requires \(Q>4\).

For critical points, define the model actions and the constants
\begin{equation*}
 I_\tau(p)=\frac12\mathcal E_\tau(p)-\frac1q\mathcal M_\tau(p),\qquad
 c_\tau=\left(\frac12-\frac1q\right)S_\tau^{q/(q-2)},
 \quad \tau\in\{F,G,E\},
\end{equation*}
where \(I_F=I_a\), \(S_F=S_a\), and \(S_E=S_E(a)\).

\begin{theorem}
\label{thm:ps-global-intro}
Let \(a>0\), \(I_a(u_k)\to c\), and \(I_a'(u_k)\to0\) in
\(H_a^*\).  Then \(c\geq0\).  After a subsequence, there are finitely many pairwise
orthogonal profiles \(p^1,\ldots,p^J\), of types
\(\tau_j\in\{F,G,E\}\), such that
\begin{equation}
 u_k=\sum_{j=1}^J R_k^jp^j+o_{H_a}(1),\qquad
 I_{\tau_j}'(p^j)=0.
 \label{eq:global-ps-decomposition}
\end{equation}
Every profile is nonzero, and
\begin{align}
 \lim_k\mathcal E_a(u_k)
 &=\sum_{j=1}^J\mathcal E_{\tau_j}(p^j)
 =\sum_{j=1}^J\mathcal M_{\tau_j}(p^j)
 =\lim_k\mathcal M(u_k),
 \label{eq:global-ps-energy}\\
 c&=\sum_{j=1}^J I_{\tau_j}(p^j),\qquad
 \sum_{j=1}^J c_{\tau_j}\leq c.
 \label{eq:global-ps-action}
\end{align}
The sum is empty precisely when \(c=0\), in which case \(u_k\to0\)
in \(H_a\).
\end{theorem}

In particular, with \(c_*:=\min\{c_F,c_G,c_E\}>0\),
\begin{equation*}
 J\leq\left\lfloor\frac{c}{c_*}\right\rfloor.
\end{equation*}
The action bound gives the following compactness threshold.

\begin{corollary}
\label{thm:stability-ps-intro}
For \(a>0\), every Palais--Smale sequence at a level
\begin{equation}
 0<c<\min\{2c_F,c_G,c_E\}
 \label{eq:ps-threshold-intro}
\end{equation}
is compact in \(H_a\) modulo one sequence of \(y\)-translations.
If \(a\in\mathcal B\), each of the levels \(2c_F\), \(c_G\), and
\(c_E\) admits a Palais--Smale sequence that is not compact modulo
\(y\)-translations.  Their minimum is the first loss of compactness
modulo these translations.
\end{corollary}

Theorem \ref{thm:perturbed-decomposition} extends the decomposition,
the minimizing-sequence alternative and the Palais--Smale conclusions
to coefficients with the same tangent limits.  The parameter interval
and endpoint criteria concern the prototype coefficients above.

\begin{remark}
At an attained finite positive endpoint, tangent recovery minimizing
sequences are noncompact modulo \(y\)-translations, whereas interior
minimizers approaching that endpoint are compact by Corollary
\ref{cor:interior-endpoint-compactness}.
The estimate \eqref{eq:grushin-remainder-rate-intro} concerns the
remaining critical mass.  The deficit modulus in Corollary
\ref{cor:deficit-intro} is qualitative.
\end{remark}

\begin{remark}
Singular profiles may lie outside \(L^2\), and escape pullbacks are
bounded only in local \(H^1\).  We prove the splitting formulas on
compact cores and synchronize their cutoff radii.  Near \(x=0\), mean-zero interpolation
does not control the transverse averages.  Their recurrence has
\(L^q\) contraction factor \(2^{-m/q}<1\), which gives summability
over the annuli.
\end{remark}

Section 2 gives the extraction criterion.  Sections 3--5 verify it
for the Grushin embedding.  Sections 6 and 7 treat sharp constants,
binding and endpoint attainment.  Section 8 proves global
Palais--Smale compactness.  Section 9 treats variable coefficients
and weights.  The appendix contains the cutoff and weighted
Brezis--Lieb lemmas.

\section{A finite multi-tangent extraction criterion}

Let \(H\) be a Hilbert space with quadratic energy \(E\) and mass
\(M(u)=\int\rho|u|^q\), \(q>2\).  For a finite index set
\(\mathcal T\), let \(H_\tau\), \(\tau\in\mathcal T\), be Hilbert model
spaces with energies \(E_\tau\), local weak topologies, and critical
masses \(M_\tau\).  Each \(p\in H_\tau\) has compactly supported approximations satisfying
\begin{equation*}
 E_\tau(p-p_L)+M_\tau(p-p_L)\longrightarrow0.
\end{equation*}
All spaces are separable.  On energy-bounded sets, local weak convergence
is determined by a countable family of local tests.
A frame \(g=(g_k)\) of type \(\tau\) has linear local pullbacks \(T_k^g\) and linear
realizations \(P_{k,L}^g p\in H\) depending only on \(p_L\).
Pullbacks of energy-bounded sequences are locally bounded and have
subsequential local weak limits \(p\in H_\tau\), with
\begin{equation*}
 E_\tau(p)\leq\liminf_k E(v_k)
 \quad\hbox{if }T_k^gv_k\rightharpoonup p\hbox{ locally}.
\end{equation*}
One family contains the identity frame with model \(H\).

\begin{theorem}
\label{thm:three-family-criterion}
Assume the following properties.
\begin{enumerate}[label=\textup{\(\roman*\)}]
\item For every frame and fixed compact cutoff \(h_L\) of
\(h\in H_\tau\),
\begin{align*}
 E(P_{k,L}^gh)&\longrightarrow E_\tau(h_L),&
 M(P_{k,L}^gh)&\longrightarrow M_\tau(h_L),
 \\
 T_k^gP_{k,L}^gh&\longrightarrow h_L\quad\hbox{locally},&
 \langle u_k,P_{k,L}^gh\rangle_H
 &\longrightarrow\langle p,h_L\rangle_{H_\tau}
\end{align*}
whenever \(T_k^gu_k\rightharpoonup p\) locally.  Such convergence has an
almost-everywhere subsequence and satisfies the local weighted Brezis--Lieb
formula.
\item Every pair of frames has an equivalent-or-orthogonal subsequence.
Equivalent frames have the same type and a unitary map
\(U^{h\leftarrow g}\) satisfying
\begin{equation}
 T_k^gv_k\rightharpoonup p\quad\Longrightarrow\quad
 T_k^hv_k\rightharpoonup U^{h\leftarrow g}p
 \label{eq:criterion-equivalent-covariance}
\end{equation}
locally on the same subsequence.  These maps satisfy the cocycle rule.  For orthogonal frames, relative
pullbacks vanish in both directions, and energy and mass decouple.
\item If \(T_k^gv_k\rightharpoonup p\) locally, then, for every fixed
compact cutoff \(p_L\), the local mass formula is
\begin{equation}
 M(v_k)-M(v_k-P_{k,L}^gp)
 \longrightarrow M_\tau(p)-M_\tau(p-p_L).
 \label{eq:criterion-fixed-core-mass}
\end{equation}

\item There are degree-two gauges \(\Theta_\tau\), and constants \(C>0\),
\(0<\gamma<q/2\), such that
\begin{equation*}
 0\leq\Theta_\tau(v)\leq CE(v),
 \qquad v\in H,\quad\tau\in\mathcal T,
\end{equation*}
and
\begin{equation*}
 \left|\Theta_\tau(v)^{1/2}-\Theta_\tau(w)^{1/2}\right|
 \leq C E(v-w)^{1/2},\qquad v,w\in H,
\end{equation*}
and
\begin{equation}
 M(v)\leq C E(v)^{q/2-\gamma}
 \left(\max_{\tau\in\mathcal T}\Theta_\tau(v)\right)^\gamma.
 \label{eq:criterion-refined}
\end{equation}
\item Every bounded sequence with
\(\limsup_k\max_\tau\Theta_\tau(v_k)\geq\theta>0\) has a weak profile
\(p\) in some frame satisfying
\begin{equation*}
 E_\tau(p)\geq c\theta,
\end{equation*}
where \(c>0\) is independent of \(\theta\) and the sequence's energy bound.
\end{enumerate}
Then every bounded sequence admits, along a subsequence,
\(J_*\in\mathbb N_0\cup\{\infty\}\)
nonzero profiles in pairwise orthogonal frames, with energy and mass
splittings.  Put \(J^\sharp=\min\{J,J_*\}\), with \(J^\sharp=J\) if
\(J_*=\infty\).  The same realization of each profile is used at every extraction
level.  Define
\begin{equation*}
 w_k^J=u_k-\sum_{j=1}^{J^\sharp}P_k^jp^j.
\end{equation*}
Then, for every fixed \(J\),
\begin{align*}
 E(u_k)&=\sum_{j=1}^{J^\sharp}E_{\tau_j}(p^j)+E(w_k^J)+o_k(1),\\
 M(u_k)&=\sum_{j=1}^{J^\sharp}M_{\tau_j}(p^j)+M(w_k^J)+o_k(1).
\end{align*}
The synchronized sum agrees with successive subtraction up to
\(o_k(1)\) in energy and mass.  The remainder is weakly zero in
each selected frame and satisfies, for \(A=\sup_kE(u_k)\),
\begin{equation}
 \eta_J\leq\frac{CA}{J+1},\qquad
 \limsup_{k\to\infty}M(w_k^J)
 \leq CA^{q/2}(J+1)^{-\gamma},
 \label{eq:criterion-remainder-rate}
\end{equation}
where
\(
\eta_J=\limsup_k\max_\tau\Theta_\tau(w_k^J)
\).
The constants depend only on \textup{(i)}--\textup{(v)} and \(\mathcal T\).
After finite stopping the terminal remainder is fixed, with zero limiting
gauges and mass.  The estimate holds for all \(J\geq1\).
\end{theorem}

\begin{proof}
Write \(\|v\|_H=E(v)^{1/2}\).  By \textup{(iv)},
\begin{equation*}
 M(v)\leq CE(v)^{q/2},\qquad
 M(v-w)^{1/q}\leq C\|v-w\|_H.
\end{equation*}
Suppose that \(T_k^gv_k\rightharpoonup p\) locally.  For fixed
\(L,L'\), \textup{(i)} and polarization give
\begin{align*}
 \lim_k E(P_{k,L}^gp-P_{k,L'}^gp)
 &=E_\tau(p_L-p_{L'}),\\
 E(v_k-P_{k,L}^gp)
 &=E(v_k)-E_\tau(p)+E_\tau(p-p_L)+o_k(1),\\
 M(v_k-P_{k,L}^gp)
 &=M(v_k)-M_\tau(p)+M_\tau(p-p_L)+o_k(1).
\end{align*}
The mass identity uses \textup{(iii)}.  Here \(L\) is fixed.

Let \(\lambda_1^h,\lambda_2^h,\ldots\) determine local weak
convergence in \(g\) and each previously selected frame.  Choose \(L_s\uparrow\infty\) such that
\begin{equation*}
 E_\tau(p-p_{L_s})^{1/2}+M_\tau(p-p_{L_s})^{1/q}
 \leq2^{-s}.
\end{equation*}
Increase \(k_s>k_{s-1}\) until, for \(k\geq k_s\),
\begin{align*}
 &|E(P_{k,L_s}^gp)-E_\tau(p_{L_s})|
   +|M(P_{k,L_s}^gp)-M_\tau(p_{L_s})|\leq2^{-s},\\
 &|E(v_k)-E(v_k-P_{k,L_s}^gp)
       -E_\tau(p)+E_\tau(p-p_{L_s})|\leq2^{-s},\\
 &|M(v_k)-M(v_k-P_{k,L_s}^gp)
       -M_\tau(p)+M_\tau(p-p_{L_s})|\leq2^{-s},\\
 &|\lambda_\ell^g(T_k^gv_k-p)|
   +|\lambda_\ell^g(T_k^gP_{k,L_s}^gp-p_{L_s})|
       \leq2^{-s},\qquad \ell\leq s,\\
 &|\lambda_\ell^h(T_k^hP_{k,L_s}^gp)|\leq2^{-s},
       \qquad \ell\leq s,\quad h\perp g.
\end{align*}
Here \(h\) ranges over the old orthogonal frames.
The finite conditions at stage \(s\) follow from \textup{(i)}--\textup{(iii)}.  Set
\begin{equation*}
 P_k^gp=P_{k,L_s}^gp\quad(k_s\leq k<k_{s+1}),\qquad
 r_k=v_k-P_k^gp.
\end{equation*}
Then
\begin{align*}
 E(P_k^gp)&\longrightarrow E_\tau(p),&
 M(P_k^gp)&\longrightarrow M_\tau(p),\\
 E(v_k)&=E_\tau(p)+E(r_k)+o_k(1),&
 M(v_k)&=M_\tau(p)+M(r_k)+o_k(1).
\end{align*}
Also,
\begin{equation*}
 T_k^gP_k^gp\rightharpoonup p,\qquad
 T_k^gr_k\rightharpoonup0,\qquad
 T_k^hP_k^gp\rightharpoonup0\quad(h\perp g\text{ old}).
\end{equation*}
For fixed \(L\), the self-energy limit and \textup{(i)} give
\begin{align*}
 \lim_k E(P_k^gp-P_{k,L}^gp)
 &=E_\tau(p)+E_\tau(p_L)-2\langle p,p_L\rangle_{H_\tau}\\
 &=E_\tau(p-p_L).
\end{align*}
Let \(h\perp g\), and let \(\lambda\) be a local weak test in \(h\).
The local pullback bound and \textup{(ii)} give
\begin{align*}
 \limsup_k|\lambda(T_k^hP_k^gp)|
 &\leq\limsup_k|\lambda(T_k^hP_{k,L}^gp)|
       +C_{h,\lambda}\limsup_k\|P_k^gp-P_{k,L}^gp\|_H\\
 &\leq C_{h,\lambda}E_\tau(p-p_L)^{1/2}.
\end{align*}
Letting \(L\to\infty\) gives
\(T_k^hP_k^gp\rightharpoonup0\).
The frame \(h\) need not be chosen before the cutoffs.

Pass to a subsequence such that \(u_k\rightharpoonup u\) in \(H\).
If \(u\ne0\), put
\begin{equation*}
 p^1=u,\qquad \iota=1,\qquad w_k^1=u_k-u.
\end{equation*}
If \(u=0\), put \(\iota=0\) and \(w_k^0=u_k\).  Retain the identity
frame in either case.  Then \(w_k^\iota\rightharpoonup0\), and
\begin{align*}
 E(u_k)&=\sum_{i=1}^{\iota}E_{\tau_i}(p^i)+E(w_k^\iota)+o_k(1),\\
 M(u_k)&=\sum_{i=1}^{\iota}M_{\tau_i}(p^i)+M(w_k^\iota)+o_k(1).
\end{align*}
Use polarization for energy, and \textup{(i)} with the Brezis--Lieb
identity in the fixed measure \(\rho\,dz\) for mass.
The exact realization \(u\) has self-energy \(E(u)\) and identity
pullback \(u\), so the fixed-core tail argument applies to it.

For any bounded sequence \(v=(v_k)\), define
\begin{equation*}
 d(v)=\sup\left\{E_\tau(p):
 T_k^gv_k\rightharpoonup p\text{ locally along a subsequence}
 \right\}.
\end{equation*}
The zero profile is allowed.  Lower semicontinuity gives
\(0\leq d(v)\leq\limsup_kE(v_k)\).
Set \(\eta(v)=\limsup_k\max_\tau\Theta_\tau(v_k)\).
For \(0<\theta<\eta(v)\), \textup{(v)} gives
\(d(v)\geq c\theta\).  Thus
\begin{equation}
 \eta(v)\leq C d(v).
 \label{eq:criterion-gauge-defect}
\end{equation}
In particular, \(d(v)=0\) implies \(\eta(v)=0\), and
\textup{(iv)} then gives \(M(v_k)\to0\).

Suppose that the profile \(p\) has been subtracted from \(v\) in
\(g\), leaving \(r\).  Let \(T_k^hr_k\rightharpoonup p'\ne0\).
If \(h\) were equivalent to \(g\),
\eqref{eq:criterion-equivalent-covariance} and
\(T_k^gr_k\rightharpoonup0\) would give \(p'=0\).
Hence \(h\perp g\).  The tail estimate above yields
\begin{equation*}
 T_k^hv_k=T_k^hr_k+T_k^hP_k^gp\rightharpoonup p'.
\end{equation*}
Every nonzero profile of \(r\) is therefore a profile of \(v\), so
\begin{equation}
 d(r)\leq d(v).
 \label{eq:criterion-defect-monotonicity}
\end{equation}
Restriction to a further subsequence can only decrease \(d\).

At level \(j\geq\iota\), put \(d_j=d(w^j)\).
If \(d_j=0\), stop and put \(J_*=j\).  Otherwise select a subsequence
and a profile \(p^{j+1}\) such that
\begin{equation*}
 T_k^{g^{j+1}}w_k^j\rightharpoonup p^{j+1},\qquad
 E_{\tau_{j+1}}(p^{j+1})\geq\frac12d_j.
\end{equation*}
The frame \(g^{j+1}\) is orthogonal to every old frame.  Indeed,
equivalence with any old frame would give a zero profile by
\eqref{eq:criterion-equivalent-covariance}.
Apply the one-profile construction and set
\begin{equation*}
 w_k^{j+1}=w_k^j-\widetilde P_k^{j+1}p^{j+1}.
\end{equation*}
At each finite level \(J\), summing the subtraction identities gives
\begin{align*}
 E(u_k)&=\sum_{i=1}^JE_{\tau_i}(p^i)+E(w_k^J)+o_k(1),\\
 M(u_k)&=\sum_{i=1}^JM_{\tau_i}(p^i)+M(w_k^J)+o_k(1).
\end{align*}
Each remainder is weakly zero in the identity frame and in all
selected frames.  If the iteration never stops, set \(J_*=\infty\).

Write the nested extraction subsequences as \(K_{j+1}\subset K_j\).
At stage \(s\), consider the first \(s\) available profiles and the
first \(s\) local tests in each selected frame.  Approximate these
profiles by compact cores with energy- and mass-norm errors at most
\(2^{-s}\).  Choose \(k_s\in K_s\), \(k_s>k_{s-1}\), so that all
of the following errors are at most \(2^{-s}\): the first \(s\)
subtraction identities, the self-energy and self-mass limits of these
cores, their pairwise orthogonality limits, and these local weak tests.
If extraction stops, use the terminal subsequence in place of \(K_s\).
Each stage has finitely many conditions.  For fixed \(J,\ell\),
they hold for every \(s\geq\max\{J,\ell\}\) on the common subsequence.

Let \(P_k^jp^j\) be the resulting synchronized realization, and let
\(\widetilde P_k^jp^j\) be its realization at the original subtraction
step.  Both have self-energy limit \(E_{\tau_j}(p^j)\) and pullback
limit \(p^j\).  For a fixed core \(p_L^j\),
\begin{align*}
 \lim_kE(P_k^jp^j-P_{k,L}^jp^j)
 &=E_{\tau_j}(p^j-p_L^j),\\
 \lim_kE(\widetilde P_k^jp^j-P_{k,L}^jp^j)
 &=E_{\tau_j}(p^j-p_L^j).
\end{align*}
Consequently,
\begin{equation*}
 \limsup_k\|P_k^jp^j-\widetilde P_k^jp^j\|_H
 \leq2E_{\tau_j}(p^j-p_L^j)^{1/2}\longrightarrow0
 \quad(L\to\infty).
\end{equation*}
For fixed \(J\), the remainders differ by \(o_H(1)\), hence by
\(o(1)\) in the mass norm.  Their masses converge together because
\begin{equation*}
 |M(v)-M(w)|\leq C\bigl(M(v)+M(w)\bigr)^{(q-1)/q}
 M(v-w)^{1/q}.
\end{equation*}
Both splittings persist.  Local pullback bounds preserve the weak
limits, and \textup{(iv)} preserves the limiting gauges.  We henceforth use these synchronized remainders,
with sums up to \(J^\sharp\) after finite stopping.

Each \(d_j\) is the defect at its extraction step.  The common
subsequence can only decrease it; the greedy energy lower bound remains
valid.  An \(o_H(1)\) change preserves every local weak profile, so
\(d_j\) also bounds the synchronized remainder's defect.
Let \(A=\sup_kE(u_k)\).
Equations \eqref{eq:criterion-defect-monotonicity} and the greedy
choice give, before stopping,
\begin{equation*}
 A\geq\sum_{j=\iota}^{J-1}E_{\tau_{j+1}}(p^{j+1})
 \geq\frac12\sum_{j=\iota}^{J-1}d_j
 \geq\frac{J-\iota}{2}d_J,\qquad J>\iota.
\end{equation*}
Also \(d_\iota\leq A\).  Since \(\iota\in\{0,1\}\),
\eqref{eq:criterion-gauge-defect} gives
\begin{equation*}
 \eta_J\leq\frac{CA}{J+1},\qquad J\geq1.
\end{equation*}
The energy splitting gives \(\limsup_kE(w_k^J)\leq A\).
Applying \eqref{eq:criterion-refined}, we obtain
\begin{equation*}
 \limsup_kM(w_k^J)
 \leq CA^{q/2-\gamma}\eta_J^\gamma
 \leq CA^{q/2}(J+1)^{-\gamma}.
\end{equation*}
After stopping, the terminal defect, gauges and mass are zero.
Thus \eqref{eq:criterion-remainder-rate} holds for all later \(J\).
Finally, the nonnegative remainders in the two splittings give
\begin{equation*}
 \sum_jE_{\tau_j}(p^j)\leq\limsup_kE(u_k),\qquad
 \sum_jM_{\tau_j}(p^j)\leq\limsup_kM(u_k)<\infty.
\end{equation*}
\end{proof}

\section{Geometry, spaces, and model inequalities}

\subsection{Grushin boxes}

For $z_0=(x_0,y_0)$ and $r>0$ define
\begin{equation*}
 \mathcal Q(z_0,r)=
 \{(x,y):|x-x_0|<r, |y-y_0|<r(r+|x_0|)\}.
\end{equation*}
Its Lebesgue measure satisfies
\begin{equation*}
 |\mathcal Q(z_0,r)|=c_{m,n}r^N(r+|x_0|)^n.
\end{equation*}
The boxes are comparable with control balls, uniformly in \(z_0,r\).

\begin{lemma}
\label{lem:local-sobolev}
There is $C=C(m,n)$ such that, for $u\in C_c^\infty(2\mathcal Q(z_0,r))$,
\begin{equation*}
\left(\int_{\mathcal Q(z_0,r)}|u|^q\right)^{2/q}
 \leq C\,r^2|\mathcal Q(z_0,r)|^{-2/Q}
 \int_{2\mathcal Q(z_0,r)}|\nabla_Gu|^2.
\end{equation*}
\end{lemma}

\begin{proof}
If $|x_0|\leq4r$, then
\begin{equation*}
 r^2|\mathcal Q(z_0,r)|^{-2/Q}
 =c_{m,n}^{-2/Q}\left(\frac r{r+|x_0|}\right)^{2n/Q}
 \geq c_{m,n}^{-2/Q}5^{-2n/Q}.
\end{equation*}
The homogeneous Grushin Sobolev inequality
\cite{Monti2006,ChenChenLi2026}, applied to the zero extension of $u$,
gives the estimate.

Suppose $|x_0|>4r$.  Put
\begin{equation*}
 v(\xi,\eta)=u(x_0+r\xi,y_0+r(r+|x_0|)\eta),
 \qquad J=r^N(r+|x_0|)^n.
\end{equation*}
On the fixed doubled box supporting $v$,
\begin{equation*}
 \frac12|x_0|\leq|x_0+r\xi|\leq\frac32|x_0|,
 \qquad
 c\leq\frac{|x_0+r\xi|^2}{(r+|x_0|)^2}\leq C.
\end{equation*}
Consequently,
\begin{equation*}
 \int|\nabla v|^2
 \leq C r^2J^{-1}\int_{2\mathcal Q(z_0,r)}|\nabla_Gu|^2.
\end{equation*}
For $N\geq3$, $q<2N/(N-2)$; Euclidean Sobolev and H\"older give
\(\|v\|_q^2\leq C\int|\nabla v|^2\).
For $N=2$, use any fixed Sobolev exponent above $q$.  Rescaling gives
\begin{equation*}
 \left(\int_{\mathcal Q(z_0,r)}|u|^q\right)^{2/q}
 \leq C r^2|\mathcal Q(z_0,r)|^{2/q-1}
 \int_{2\mathcal Q(z_0,r)}|\nabla_Gu|^2.
\end{equation*}
Use $2/q-1=-2/Q$.
\end{proof}

\begin{lemma}
\label{lem:uniform-local-rellich}
Let \(z_k=(x_k,y_k)\), \(r_k>0\), and define on
\((-2,2)^m\times(-2,2)^n\)
\begin{equation*}
 w_k(\xi,\eta)=v_k\bigl(x_k+r_k\xi,
 y_k+r_k(r_k+|x_k|)\eta\bigr).
\end{equation*}
Put
\begin{equation*}
 a_k(\xi)=\frac{|x_k+r_k\xi|^2}{(r_k+|x_k|)^2}.
\end{equation*}
Assume
\begin{equation*}
 \int_{(-2,2)^N}
 (|\nabla_\xi w_k|^2+a_k|\nabla_\eta w_k|^2+|w_k|^2)\leq C_0.
\end{equation*}
Suppose, for one fixed \(C>2\), that one of the following holds:
\begin{equation*}
 |x_k|\leq2r_k,\qquad
 2r_k<|x_k|\leq Cr_k,\qquad
 |x_k|/r_k\longrightarrow\infty.
\end{equation*}
Then a subsequence converges strongly in
\(L^2((-1,1)^N)\).  The compactness constant depends only on
\(m,n,C,C_0\).  After a singular or escape pullback, the same assertion
holds on every fixed model box.
\end{lemma}

\begin{proof}
Put
\(c_k=x_k/(r_k+|x_k|)\) and
\(\lambda_k=r_k/(r_k+|x_k|)\).  The horizontal fields are
\begin{equation*}
 \partial_{\xi_i},\qquad
 (c_{k,i}+\lambda_k\xi_i)\partial_{\eta_j}.
\end{equation*}
In the first two regimes, \(\lambda_k\geq\lambda_0(C)>0\).  Choose \(\chi\) supported in
the doubled box and equal to one on
\((-3/2,3/2)^N\).  For \(F\in C_c^\infty(\mathbb R^m)\) and \(\tau\in\mathbb R^n\),
the translated harmonic oscillator estimate is
\begin{equation*}
 \lambda_k|\tau|\int_{\mathbb R^m}|F|^2
 \leq C\int_{\mathbb R^m}
 \left(|\nabla_\xi F|^2+|c_k+\lambda_k\xi|^2
 |\tau|^2|F|^2\right).
\end{equation*}
For \(\tau\ne0\), translate by \(c_k/\lambda_k\) and scale by
\((\lambda_k|\tau|)^{1/2}\), then use the first eigenvalue of
\(-\Delta_\xi+|\xi|^2\).  The case \(\tau=0\) is immediate.  For \(\chi w_k\), the product rule gives
\begin{equation*}
 \int\left(|\nabla_\xi(\chi w_k)|^2
 +|c_k+\lambda_k\xi|^2|\nabla_\eta(\chi w_k)|^2\right)
 \leq C\int_{(-2,2)^N}
 \left(|\nabla_\xi w_k|^2+a_k|\nabla_\eta w_k|^2+|w_k|^2\right)
 \leq C.
\end{equation*}
The coefficient is bounded on \(\mathrm{supp}\,\chi\).  Apply the oscillator
estimate to the \(\eta\)-Fourier transform of \(\chi w_k\).  Plancherel gives
\begin{equation*}
 \int_{\mathbb R^n}\int_{\mathbb R^m}
 (1+|\tau|)|\widehat{\chi w_k}(\xi,\tau)|^2\,d\xi d\tau
 +\int|\nabla_\xi(\chi w_k)|^2\leq C.
\end{equation*}
Let \(f_{k,R}\) be the inverse Fourier transform of
\(\mathbf1_{\{|\tau|\leq R\}}\widehat{\chi w_k}\) in \(\eta\).
For each fixed \(R\),
\begin{equation*}
 \|f_{k,R}\|_2^2+\|\nabla_\xi f_{k,R}\|_2^2
 +\|\nabla_\eta f_{k,R}\|_2^2\leq C(1+R^2),
 \qquad
 \|\chi w_k-f_{k,R}\|_2^2\leq C/R.
\end{equation*}
Apply Rellich on the inner box for \(R=1,2,\ldots\).  A diagonal
subsequence makes every fixed \(R\) sequence \(f_{k,R}\) Cauchy there.  Hence
\begin{equation*}
 \limsup_{k,l\to\infty}
 \|w_k-w_l\|_{L^2((-1,1)^N)}\leq 2C^{1/2}R^{-1/2}.
\end{equation*}
Let \(R\to\infty\).

In the third regime, \(\lambda_k\to0\) and, after a subsequence,
\(c_k\to e\in\mathbb S^{m-1}\).  Hence
\begin{equation*}
 a_k(\xi)\longrightarrow1
\end{equation*}
uniformly on the doubled box.  Eventually \(a_k\geq1/2\), so the
assumed bound gives the \(H^1\) bound required by Rellich.
For singular pullbacks, \(x_k=0\); for escape pullbacks, \(r_k=1\)
and \(|x_k|\to\infty\).  Apply these estimates on boxes of radii
\(1,2,3,\ldots\), and take a nested diagonal subsequence.
\end{proof}

\begin{proposition}
\label{prop:embedding}
There is $C=C(m,n,a)>0$ such that
\begin{equation*}
 \mathcal M(u)^{2/q}\leq C\mathcal E_a(u),
 \qquad u\in H_a.
\end{equation*}
Consequently $S_a>0$.
\end{proposition}

\begin{proof}
First let \(u\in C_c^\infty\).
Cover $\{|x|\leq4\}$ by boxes $\mathcal Q(z_j,1)$ with bounded
overlap.  Choose $\sum_j\chi_j^2=1$ on this set, with $\chi_j$
supported in the doubled box and $|\nabla_G\chi_j|\leq C$.  The
weight is bounded on the supports.  Write
\(e_j=\mathcal E_a(\chi_ju)\).  The product rule gives
\begin{equation*}
 \sum_j e_j\leq C\int(|\nabla_Gu|^2+|u|^2)
 \leq C_a\mathcal E_a(u).
\end{equation*}
Lemma \ref{lem:local-sobolev} on the doubled boxes gives
\begin{equation*}
 \sum_j\mathcal E_a(\chi_ju)\leq C_a\mathcal E_a(u),\qquad
 \sum_j\int|\chi_ju|^q\leq C_a\mathcal E_a(u)^{q/2}.
\end{equation*}
At most \(K\) terms are nonzero at each point.  H\"older gives
\begin{equation*}
 1\leq K^{1-2/q}\left(\sum_j|\chi_j|^q\right)^{2/q}.
\end{equation*}
Multiplying by \(|u|^q\) and integrating proves
\begin{equation*}
 \int_{|x|\leq4}\omega^{q/Q}|u|^q
 \leq C\mathcal E_a(u)^{q/2}.
\end{equation*}

For $\ell\geq2$, put
\begin{equation*}
 \mathcal R_\ell=\{2^\ell<|x|<2^{\ell+1}\}.
\end{equation*}
Cover $\mathcal R_\ell$ by boxes of $x$ radius one and $y$ radius
$2^\ell$, with a quadratic partition on the doubled boxes and uniformly
bounded overlap.  The partition derivatives in \(x\) and \(y\) are bounded by
\(C\) and \(C2^{-\ell}\).  Thus its Grushin derivatives are bounded
by \(C\).  On each box,
\begin{equation*}
 \omega(x)^{q/Q}\asymp2^{nq\ell/Q}.
\end{equation*}
For \(f=\chi_{\ell j}u\) and
\(v(\xi,\eta)=2^{n\ell/2}f(x_{\ell j}+\xi,
y_{\ell j}+2^\ell\eta)\), then
\begin{equation*}
 \int\left(|\nabla_\xi v|^2+|\nabla_\eta v|^2+a|v|^2\right)
 \leq C_a\mathcal E_a(f).
\end{equation*}
Here \(|x|\asymp2^\ell\) on the support of \(f\).  Since
\begin{equation*}
 -\frac{nq}{2}+n+\frac{2n}{Q-2}=0.
\end{equation*}
the local $H^1\to L^q$ estimate gives
\begin{equation*}
 \int\omega^{q/Q}|\chi_{\ell j}u|^q
 \leq C_a\mathcal E_a(\chi_{\ell j}u)^{q/2}.
\end{equation*}
For the quadratic
partition \((\chi_{\ell j})_j\), the product rule and bounded overlap give
\begin{equation*}
 \sum_j\int\left(|\nabla_G(\chi_{\ell j}u)|^2
 +a|\chi_{\ell j}u|^2\right)
 \leq C_a\int_{\{2^{\ell-1}<|x|<2^{\ell+2}\}}
 (|\nabla_Gu|^2+a|u|^2).
\end{equation*}
Consequently,
\begin{equation*}
 \int_{\mathcal R_\ell}\omega^{q/Q}|u|^q
 \leq C\left(
 \int_{\{2^{\ell-1}<|x|<2^{\ell+2}\}}
 (|\nabla_Gu|^2+a|u|^2)\right)^{q/2}.
\end{equation*}
Sum in $\ell$.  Since $q/2>1$,
\begin{equation*}
 \sum_\ell E_\ell^{q/2}\leq
 \left(\sum_\ell E_\ell\right)^{q/2}.
\end{equation*}
Bounded overlap proves the estimate for smooth \(u\).  Choose \(u_h\in C_c^\infty\) with
\(\mathcal E_a(u_h-u)\to0\).  Applying the estimate to \(u_h-u_l\) gives a weighted \(L^q\) limit.
It equals \(u\) by \(L^2\) convergence and an almost-everywhere subsequence.
Taking the infimum over \(\mathcal M(u)=1\) gives \(S_a>0\).
\end{proof}

\subsection{The singular model}

The homogeneous inequality
\begin{equation}
 S_G\|\phi\|_q^2\leq\mathcal E_G(\phi)
 \label{eq:singular-sharp-inequality}
\end{equation}
is invariant under $y$ translation and
\begin{equation*}
 \phi(x,y)\mapsto r^{-(Q-2)/2}\phi(x/r,y/r^2).
\end{equation*}
For every \(m,n\geq1\), the constant in
\eqref{eq:singular-sharp-inequality} is attained
\cite[Theorem~1.1]{GandalLoiudiceTyagi2025}.  Every positive normalized
extremal \(\Phi\) satisfies
\begin{equation}
 \Phi(z)\asymp d_G(z,0)^{2-Q}\quad\hbox{as }d_G(z,0)\to\infty,
 \label{eq:singular-extremal-decay}
\end{equation}
where \(d_G\) is a homogeneous control distance.

\begin{proposition}
\label{prop:automatic-singular-moments}
If \(Q>4\), every positive normalized extremal \(\Phi\) for \(S_G\) obeys
\begin{equation*}
 \Phi\in L^2(\mathbb R^N),\qquad
 \int_{\mathbb R^N}(1+|x|^2)^{s+1}|\Phi|^q<\infty,
 \qquad s=\frac n{Q-2}.
\end{equation*}
\end{proposition}

\begin{proof}
Local regularity bounds \(\Phi\) on compact sets.  At infinity,
\(d_G\asymp\rho=(|x|^4+4|y|^2)^{1/4}\) and \(|x|\leq C\rho\).
The dilation \((x,y)\mapsto(rx,r^2y)\) gives
\begin{equation*}
 |\{\rho<r\}|=r^Q|\{\rho<1\}|.
\end{equation*}
Thus radial integration for \(\rho\) has density \(Cr^{Q-1}\,dr\).
By \eqref{eq:singular-extremal-decay},
\begin{equation*}
 \int_{\{\rho>1\}}\Phi^2
 \leq C\int_1^\infty r^{-2(Q-2)}r^{Q-1}\,dr<\infty
\end{equation*}
Since \(q(Q-2)=2Q\),
\begin{align*}
 \int_{\{\rho>1\}}(1+|x|^2)^{s+1}|\Phi|^q
 &\leq C\int_1^\infty r^{2(s+1)-2Q}r^{Q-1}\,dr\\
 &=C\int_1^\infty r^{-Q+1+2s}\,dr.
\end{align*}
The last integral converges: \(2s<Q-2\), since
\(2n\leq Q-1<(Q-2)^2\) for \(Q>4\).
The integrals over \(\{\rho\leq1\}\) are finite by local boundedness.
\end{proof}

\begin{remark}
The lower decay bound in \eqref{eq:singular-extremal-decay} also gives
\(\int_{\rho>1}\Phi^2\geq c\int_1^\infty r^{3-Q}\,dr\).
Thus positive extremals fail to lie in \(L^2\) when \(Q=3,4\).
Section 4 realizes these profiles by cutoffs.  The uncut second-order
singular trial requires \(Q>4\).
\end{remark}

For \(n=1\) and even \(m\), two opposite Heisenberg structures give an explicit extremal.

\begin{proposition}
\label{prop:heisenberg-comparison}
Let \(n=1\), let \(m\geq2\) be even, and put
\begin{equation*}
 b_m=\frac{\pi^{(m+2)/2}}
 {2^{m+1}\Gamma((m+2)/2)},\qquad
 F_m(x,y)=\bigl((1+|x|^2)^2+4y^2\bigr)^{-m/4}.
\end{equation*}
Then
\begin{equation}
 S_G=m^2b_m^{2/(m+2)},\qquad
 \Phi_m=b_m^{-1/q}F_m
 \label{eq:explicit-singular-extremal}
\end{equation}
is a normalized positive extremal.  If \(m>2\), then
\begin{equation}
 \int |x|^2\Phi_m^q=1,\qquad
 \int\Phi_m^2=\frac{4m}{m-2}b_m^{2/(m+2)}.
 \label{eq:explicit-singular-moments}
\end{equation}
\end{proposition}

\begin{proof}
Choose a real orthogonal skew-symmetric matrix \(J\) on \(\mathbb R^m\),
so \(J^2=-\mathrm{Id}\).  Define
\begin{equation*}
 H_i^\pm=\partial_{x_i}\pm(Jx)_i\partial_y,\qquad
 E_\pm(u)=\sum_{i=1}^m\int|H_i^\pm u|^2.
\end{equation*}
Reflection \(y\mapsto-y\) interchanges the two structures.  Denote their
common sharp constant by \(S_H\).
Set \(t=2y\) and \(u(x,y)=v(x,2y)\).  Then
\begin{align*}
 E_\pm(u)
 &=\frac12\int\sum_{i=1}^m
 \left|\bigl(\partial_{x_i}\pm2(Jx)_i\partial_t\bigr)v\right|^2
 \,dx\,dt,\\
 \|u\|_q^2&=2^{-2/q}\|v\|_q^2.
\end{align*}
If \(S_H^{(2)}\) denotes the sharp constant for the fields with
coefficient \(2\), the Jacobian \(dx\,dy=\tfrac12dx\,dt\) gives
\begin{equation*}
 S_H=2^{-(1-2/q)}S_H^{(2)}=2^{-2/Q}S_H^{(2)}.
\end{equation*}
The sharp Heisenberg inequality \cite{JerisonLee1988} has the extremal
\(v(x,t)=((1+|x|^2)^2+t^2)^{-m/4}\), whose pullback is \(F_m\).
Averaging cancels the mixed terms:
\begin{equation*}
 \frac12\bigl(E_+(u)+E_-(u)\bigr)=\mathcal E_G(u)
 \geq S_H\|u\|_q^2.
\end{equation*}
Since \(F_m\) is radial in \(x\),
\((Jx)\cdot\nabla_xF_m=0\).  Hence
\(E_+(F_m)=E_-(F_m)=\mathcal E_G(F_m)\), so \(S_G\leq S_H\).  The preceding lower bound gives equality.
To pass from compact functions to \(D_G^{1,2}\), take
\(u_h\in C_c^\infty\) with
\(\mathcal E_G(u_h-u)\to0\).  Since
\begin{equation*}
 E_\pm(u_h-u_l)\leq2\mathcal E_G(u_h-u_l),
 \qquad
 \|u_h-u_l\|_q^2\leq S_G^{-1}\mathcal E_G(u_h-u_l),
\end{equation*}
the identities pass to \(D_G^{1,2}\).

Put \(H=(1+|x|^2)^2+4y^2\).  Direct differentiation gives
\begin{equation*}
 |\nabla_GH|^2=16|x|^2H,\qquad
 \Delta_GH=4m+4(m+4)|x|^2.
\end{equation*}
Put \(\alpha=m/4\).  The chain rule gives
\begin{equation*}
 -\Delta_GF_m
 =\alpha H^{-\alpha-1}\Delta_GH
 -\alpha(\alpha+1)H^{-\alpha-2}|\nabla_GH|^2
 =m^2H^{-\alpha-1}.
\end{equation*}
Since \(\alpha(q-1)=\alpha+1\), this is
\begin{equation*}
 -\Delta_xF_m-|x|^2\partial_y^2F_m=m^2F_m^{q-1}.
\end{equation*}
For \(\alpha>1/2\) and \(A>0\), integration in \(y\) gives
\begin{equation*}
 \int_{\mathbb R}(A^2+4y^2)^{-\alpha}\,dy
 =\frac{\sqrt\pi\,\Gamma(\alpha-1/2)}
 {2\Gamma(\alpha)}A^{1-2\alpha}.
\end{equation*}
For \(\beta>m/2\), polar integration and the beta integral give
\begin{equation*}
 \int_{\mathbb R^m}(1+|x|^2)^{-\beta}\,dx
 =\frac{\pi^{m/2}\Gamma(\beta-m/2)}{\Gamma(\beta)}.
\end{equation*}
Apply the two integration formulas with
\(\alpha=(m+2)/2\) for \(F_m^q\).  For the \(|x|^2\) moment,
the ratio of the radial beta integrals is
\begin{equation*}
 \frac{\Gamma(m/2+1)\Gamma(m/2)}
 {\Gamma(m/2)\Gamma(m/2+1)}=1.
\end{equation*}
For \(F_m^2\), take \(\alpha=m/2\), \(\beta=m-1\), requiring \(m>2\).
The duplication formula yields
\begin{equation*}
 \int F_m^q=b_m,\qquad
 \int |x|^2F_m^q=b_m,
 \qquad
 \int F_m^2=\frac{4m}{m-2}b_m\quad(m>2).
\end{equation*}
Test the equation with \(\chi_L^2F_m\), where \(\chi_L\) is the
cutoff from Appendix \ref{sec:cutoffs}.  The error is bounded by
\begin{equation*}
 2\left(\int_{\{L<\rho<2L\}}|\nabla_GF_m|^2\right)^{1/2}
 \left(\int F_m^2|\nabla_G\chi_L|^2\right)^{1/2}.
\end{equation*}
H\"older bounds the second factor:
\begin{equation*}
 \int F_m^2|\nabla_G\chi_L|^2
 \leq C L^{-2}|\{L<\rho<2L\}|^{1-2/q}
 \left(\int_{\{L<\rho<2L\}}F_m^q\right)^{2/q}
 \leq C\left(\int_{\{L<\rho<2L\}}F_m^q\right)^{2/q}\to0.
\end{equation*}
Thus \(\mathcal E_G(F_m)=m^2b_m\).  Multiplication by
\(b_m^{-1/q}\) then gives
\eqref{eq:explicit-singular-extremal} and
\eqref{eq:explicit-singular-moments}.
\end{proof}

\subsection{The escape model}

For \(\psi\in H^1(\mathbb R^N)\), put
\begin{equation*}
 \mathcal E_E^a(\psi)
 =\int_{\mathbb R^N}(|\nabla\psi|^2+a|\psi|^2)\,d\zeta.
\end{equation*}
Let $Z_a$ be a positive minimizer for \(S_E(a)\), normalized by
\begin{equation*}
 \int_{\mathbb R^N}Z_a^q=1.
\end{equation*}
It solves
\begin{equation*}
 -\Delta Z_a+aZ_a=S_E(a)Z_a^{q-1}.
\end{equation*}
By \cite{GidasNiNirenberg1979,Kwong1989}, $Z_a$ is radial and unique
up to translation.  Centre it at zero.  Every real normalized extremal
has constant sign: its absolute value is an extremal, and the strong
maximum principle makes it strictly positive.  Continuity excludes
a sign change.

\begin{lemma}
\label{lem:euclidean-pohozaev}
One has
\begin{align*}
 Z_a(\zeta)&=a^{N/(2q)}Z_1(\sqrt a\,\zeta),\\
 S_E(a)&=a^{n/Q}S_E(1),\\
 \int|\nabla Z_a|^2&=\frac{N}{Q}S_E(a),\\
 a\int Z_a^2&=\frac{n}{Q}S_E(a).
\end{align*}
\end{lemma}

\begin{proof}
For \(\psi\in H^1\), put
\(\psi_a(\zeta)=a^{N/(2q)}\psi(\sqrt a\,\zeta)\).
The change of variables \(\theta=\sqrt a\,\zeta\) gives
\begin{equation*}
 \int|\psi_a|^q=\int|\psi|^q,\qquad
 \mathcal E_E^a(\psi_a)
 =a^{1-N/2+N/q}\mathcal E_E^1(\psi).
\end{equation*}
The exponent is
\begin{equation*}
 1-\frac N2+\frac Nq=\frac nQ.
\end{equation*}
Taking infima and applying the inverse scaling gives
\(S_E(a)=a^{n/Q}S_E(1)\).  Uniqueness identifies the centered positive
minimizer with the scaled \(Z_1\).

Testing the Euler equation with \(Z_a\) gives
\begin{equation}
 \int|\nabla Z_a|^2+a\int Z_a^2=S_E(a).
 \label{eq:escape-euler-identity}
\end{equation}
The exponential decay
\(|Z_a(\zeta)|+|\nabla Z_a(\zeta)|\leq Ce^{-c|\zeta|}\)
bounds the boundary terms from testing on \(B_R\) with
\(\zeta\cdot\nabla Z_a\) by
\begin{equation*}
 C R^{N+2}\left(e^{-2cR}+e^{-qcR}\right)\longrightarrow0.
\end{equation*}
Integration by parts and \(R\to\infty\) give
\begin{equation}
 \frac{N-2}{2}\int|\nabla Z_a|^2+
 \frac{Na}{2}\int Z_a^2=\frac NqS_E(a).
 \label{eq:escape-pohozaev-identity}
\end{equation}
Subtract \((N-2)/2\) times \eqref{eq:escape-euler-identity} from
\eqref{eq:escape-pohozaev-identity}.  This gives
\begin{equation*}
 a\int Z_a^2
 =\left(\frac Nq-\frac{N-2}{2}\right)S_E(a)
 =\frac nQ S_E(a).
\end{equation*}
The Euler identity then gives
\(\int|\nabla Z_a|^2=(N/Q)S_E(a)\).
\end{proof}

\section{Frames and orthogonality}

For the three models, write
\begin{align*}
 (\mathcal E_F,\mathcal M_F)=(\mathcal E_a,\mathcal M),\quad
 \mathcal M_G(\phi)=\int_{\mathbb R^N}|\phi|^q,\quad
 \mathcal E_E(\psi)=\mathcal E_E^a(\psi),\qquad
 \mathcal M_E(\psi)=\int_{\mathbb R^N}|\psi|^q.
\end{align*}
The singular energy \(\mathcal E_G\) is defined in Section 1.

\subsection{Finite frames}

For $b\in\mathbb R^n$ set
\begin{equation*}
 (T_bu)(x,y)=u(x,y-b).
\end{equation*}
Then $T_b$ is an isometry of $H_a$ and the weighted $L^q$ space.  Two
finite frames $(b_k)$ and $(\widetilde b_k)$ are orthogonal if
\begin{equation}
 |b_k-\widetilde b_k|\longrightarrow\infty.
 \label{eq:finite-frame-orthogonality}
\end{equation}

\subsection{Singular frames}

Let $r>0$ and $c\in\mathbb R^n$.  Define
\begin{equation}
 (D_{r,c}\phi)(x,y)
 =r^{-(Q-2)/2}\phi(x/r,(y-c)/r^2).
 \label{eq:singular-frame-map}
\end{equation}
Then
\begin{align*}
 \mathcal E_G(D_{r,c}\phi)=\mathcal E_G(\phi),\qquad
 \|D_{r,c}\phi\|_q=\|\phi\|_q,\qquad
 \|D_{r,c}\phi\|_2^2=r^2\|\phi\|_2^2.
\end{align*}
A singular frame is a sequence $(r_k,c_k)$ with $r_k\to0$.  Two such frames
are orthogonal if
\begin{equation}
 \left|\log\frac{r_k}{\widetilde r_k}\right|
 +\frac{|c_k-\widetilde c_k|}{r_k^2+\widetilde r_k^2}
 \longrightarrow\infty.
 \label{eq:singular-frame-orthogonality}
\end{equation}

For $\phi\in D_G^{1,2}\setminus L^2$, the expression in
\eqref{eq:singular-frame-map} need not belong to $H_a$.  Choose
$L_k\to\infty$ so slowly that
\[
 r_kL_k\to0,\qquad
 \|\chi_{L_k}\phi-\phi\|_{D_G^{1,2}}\to0,\qquad
 \|\chi_{L_k}\phi-\phi\|_q\to0,
\]
where $\chi_L$ is defined in Appendix \ref{sec:cutoffs}.  Put
\[
 [D_{r_k,c_k}\phi]_{\mathrm{rel}}
 =D_{r_k,c_k}(\chi_{L_k}\phi).
\]
The ambient zeroth-order term vanishes because Sobolev--H\"older gives
\begin{equation*}
 r_k^2\|\chi_{L_k}\phi\|_2^2
 \leq C(r_kL_k)^2\|\phi\|_q^2\longrightarrow0.
\end{equation*}
Thus two admissible realizations differ by $o(1)$ in the ambient
energy and critical mass norms.
We use $D_{r_k,c_k}\phi$ for this realized sequence in all
decomposition formulas.  This convention is needed when $Q\leq4$.

\subsection{Escape frames}

Let $X\in\mathbb R^m\setminus\{0\}$, $R=|X|$, and $b\in\mathbb R^n$.
For $\psi\in H^1(\mathbb R^N)$ define
\begin{equation}
 (A_{X,b}\psi)(x,y)
 =R^{-n/2}\psi(x-X,(y-b)/R).
 \label{eq:escape-frame-map}
\end{equation}
An escape frame satisfies $R_k=|X_k|\to\infty$.  Two escape frames are
orthogonal if
\begin{equation}
 |X_k-\widetilde X_k|
 +\left|\log\frac{R_k}{\widetilde R_k}\right|
 +\frac{|b_k-\widetilde b_k|}{R_k+\widetilde R_k}
 \longrightarrow\infty.
 \label{eq:escape-frame-orthogonality}
\end{equation}

The pullback energy space depends on $X$.  For $h\in C_c^\infty$ put
\begin{equation}
 \|h\|_{X,a}^2=
 \int_{\mathbb R^N}\left(
 |\nabla_\xi h|^2+
 \frac{|X+\xi|^2}{|X|^2}|\nabla_\eta h|^2+a|h|^2\right).
 \label{eq:moving-escape-norm}
\end{equation}
Let $H_{X,a}$ be its completion.  We write
\begin{equation*}
 h_k\rightharpoonup_E h
\end{equation*}
if $\sup_k\|h_k\|_{X_k,a}<\infty$ and, for every $R>0$,
\begin{equation*}
 h_k\rightharpoonup h\quad\hbox{in }H^1(B_R).
\end{equation*}
The limits agree on overlapping balls.  Lower semicontinuity on
$B_R$, followed by $R\to\infty$, gives $h\in H^1(\mathbb R^N)$.

\subsection{The six orthogonality relations}

Frames of different types are orthogonal because
\begin{align*}
 F\perp G&:\quad r_k\longrightarrow0,\\
 F\perp E&:\quad |X_k|\longrightarrow\infty,\\
 G\perp E&:\quad r_k\longrightarrow0,\quad |X_k|\longrightarrow\infty.
\end{align*}
For a general $H^1$ profile, realize \eqref{eq:escape-frame-map} by a cutoff.
Choose
$\vartheta\in C_c^\infty([0,\infty))$ with $\vartheta=1$ on $[0,1]$ and
$\vartheta=0$ on $[2,\infty)$.  Choose $L_k\to\infty$ so slowly that
$L_k/R_k\to0$ and put
\[
 [A_{X_k,b_k}\psi]_{\mathrm{rel}}
 =A_{X_k,b_k}\bigl(\vartheta(|\cdot|/L_k)\psi\bigr).
\]
The cutoff converges to $\psi$ in $H^1$ and $L^q$, and
$|X_k+\xi|/R_k\to1$ uniformly on its support.  Thus admissible
realizations differ by $o(1)$ in energy and mass norm.
We write $A_{X_k,b_k}\psi$ for a realization in the profile formulas.  No cutoff is needed if
$\int|\xi|^2|\nabla_\eta\psi|^2<\infty$, as for $Z_a$.

\begin{lemma}
\label{lem:one-frame-limits}
Let $\phi\in C_c^\infty(\mathbb R^N)$ and $r_k\to0$.  Then
\begin{align*}
 \mathcal E_a(D_{r_k,c_k}\phi)&\longrightarrow\mathcal E_G(\phi),\\
 \mathcal M(D_{r_k,c_k}\phi)&\longrightarrow\|\phi\|_q^q.
\end{align*}
Let $\psi\in C_c^\infty(\mathbb R^N)$, $R_k=|X_k|\to\infty$, and
$X_k/R_k\to e$.  Then
\begin{align*}
 \mathcal E_a(A_{X_k,b_k}\psi)&\longrightarrow
 \mathcal E_E^a(\psi),\\
 \mathcal M(A_{X_k,b_k}\psi)&\longrightarrow\|\psi\|_q^q.
\end{align*}
The conclusions extend to $D_G^{1,2}$ and $H^1$, respectively.
\end{lemma}

\begin{proof}
The singular change of variables gives
\begin{align*}
 \mathcal E_a(D_{r_k,c_k}\phi)
 &=\mathcal E_G(\phi)+ar_k^2\|\phi\|_2^2,\\
 \mathcal M(D_{r_k,c_k}\phi)
 &=\int(1+r_k^2|\xi|^2)^{n/(Q-2)}|\phi(\xi,\eta)|^q
 \,d\xi\,d\eta.
\end{align*}
Let \(K\) contain the support of \(\phi\).  Then
\begin{equation*}
 \left|\mathcal M(D_{r_k,c_k}\phi)-\|\phi\|_q^q\right|
 \leq\sup_K\left|(1+r_k^2|\xi|^2)^{n/(Q-2)}-1\right|
 \|\phi\|_q^q\longrightarrow0.
\end{equation*}

For the escape frame put \(x=X_k+\xi\), \(y=b_k+R_k\eta\).
The two transformed coefficients are
\begin{equation*}
 c_k(\xi)=\frac{|X_k+\xi|^2}{R_k^2},\qquad
 w_k(\xi)=R_k^{-2n/(Q-2)}
 (1+|X_k+\xi|^2)^{n/(Q-2)}.
\end{equation*}
Here \(-nq/2+n=-2n/(Q-2)\).  Hence
\begin{align*}
 \mathcal E_a(A_{X_k,b_k}\psi)
 &=\int\left(|\nabla_\xi\psi|^2+c_k|\nabla_\eta\psi|^2
 +a|\psi|^2\right),\\
 \mathcal M(A_{X_k,b_k}\psi)&=\int w_k|\psi|^q.
\end{align*}
On every fixed compact set \(K\),
\begin{equation*}
 \sup_K|c_k-1|+\sup_K|w_k-1|\longrightarrow0.
\end{equation*}
For \(K\) containing the support of \(\psi\), this gives
\begin{align*}
 \left|\mathcal E_a(A_{X_k,b_k}\psi)-\mathcal E_E^a(\psi)\right|
 &\leq\sup_K|c_k-1|\,\|\nabla_\eta\psi\|_2^2,\\
 \left|\mathcal M(A_{X_k,b_k}\psi)-\|\psi\|_q^q\right|
 &\leq\sup_K|w_k-1|\,\|\psi\|_q^q.
\end{align*}
For a general singular profile use \(\phi_k=\chi_{L_k}\phi\).
Its support satisfies \(|\xi|\leq CL_k\).  The realization conditions give
\begin{align*}
 \mathcal E_a(D_{r_k,c_k}\phi_k)
 &=\mathcal E_G(\phi_k)+ar_k^2\|\phi_k\|_2^2,\\
 r_k^2\|\phi_k\|_2^2
 &\leq C(r_kL_k)^2\|\phi\|_q^2\longrightarrow0,\\
 \left|\mathcal M(D_{r_k,c_k}\phi_k)-\|\phi_k\|_q^q\right|
 &\leq\sup_{|\xi|\leq CL_k}
 \left|(1+r_k^2|\xi|^2)^{n/(Q-2)}-1\right|\|\phi_k\|_q^q
 \longrightarrow0.
\end{align*}
Now \(\phi_k\to\phi\) in the two model norms.

For a general escape profile put
\(\psi_k=\vartheta(|\cdot|/L_k)\psi\).
Since \(L_k/R_k\to0\),
\begin{equation*}
 \sup_{|\xi|\leq2L_k}|c_k-1|
 +\sup_{|\xi|\leq2L_k}|w_k-1|\longrightarrow0.
\end{equation*}
Apply the same estimates to \(\psi_k\).  Its model norms converge to
those of \(\psi\), so the errors tend to zero.
\end{proof}

\begin{lemma}
\label{lem:relative-frames}
Let two frames be orthogonal, of the same or different families.  The
relative pullback of a fixed compactly supported profile converges weakly
to zero in $H_a$ for a finite receiving frame and in $D_G^{1,2}$ for a
singular receiving frame.  For an escape receiving frame, the pullbacks
are bounded in the moving norms \eqref{eq:moving-escape-norm} and converge
weakly to zero in $H^1(B_R)$ for every fixed $R>0$.
\end{lemma}

\begin{proof}
For \(F\)--\(F\), the relative map is \(T_{b_k-\widetilde b_k}\).
By \eqref{eq:finite-frame-orthogonality}, if \(f,g\) are compactly
supported, their translated supports are disjoint for large \(k\).  Thus
\begin{equation*}
 \langle T_{b_k-\widetilde b_k}f,g\rangle_{H_a}=0
 \quad\hbox{for large }k.
\end{equation*}
For the other pairs, direct substitution gives
\begin{align}
 D_{r,c}^{-1}T_bf(\xi,\eta)
 &=r^{(Q-2)/2}f(r\xi,c+r^2\eta-b),
 \label{eq:relative-gf-formula}\\
 T_b^{-1}D_{r,c}\phi(x,y)
 &=r^{-(Q-2)/2}\phi\left(\frac{x}{r},
 \frac{y+b-c}{r^2}\right),
 \notag\\
 A_{X,b}^{-1}T_c f(\xi,\eta)
 &=R^{n/2}f(X+\xi,b+R\eta-c),
 \notag\\
 T_c^{-1}A_{X,b}\psi(x,y)
 &=R^{-n/2}\psi\left(x-X,\frac{y+c-b}{R}\right),
 \notag\\
 D_{r,c}^{-1}A_{X,b}\psi(\xi,\eta)
 &=r^{(Q-2)/2}R^{-n/2}
 \psi\left(r\xi-X,\frac{c+r^2\eta-b}{R}\right),
 \notag\\
 A_{X,b}^{-1}D_{r,c}\phi(\xi,\eta)
 &=R^{n/2}r^{-(Q-2)/2}
 \phi\left(\frac{X+\xi}{r},\frac{b+R\eta-c}{r^2}\right).
 \label{eq:relative-eg-formula}
\end{align}
Here \(R=|X|\).  For equal model types,
\begin{align*}
 D_{r,c}^{-1}D_{\widetilde r,\widetilde c}\phi(\xi,\eta)
 &=\left(\frac r{\widetilde r}\right)^{(Q-2)/2}
 \phi\left(\frac r{\widetilde r}\xi,
 \frac{c-\widetilde c}{\widetilde r^2}
 +\frac{r^2}{\widetilde r^2}\eta\right),\\
 A_{X,b}^{-1}A_{\widetilde X,\widetilde b}\psi(\xi,\eta)
 &=\left(\frac R{\widetilde R}\right)^{n/2}
 \psi\left(\xi+X-\widetilde X,
 \frac{b-\widetilde b}{\widetilde R}
 +\frac R{\widetilde R}\eta\right).
\end{align*}
All receiving pullbacks have uniform energy bounds.
For a singular receiver this follows from
\(\mathcal E_G(D_{r,c}^{-1}v)=\mathcal E_G(v)\leq\mathcal E_a(v)\).
For an escape receiver the change of variables gives
\begin{equation*}
 \|A_{X,b}^{-1}v\|_{X,a}^2=\mathcal E_a(v).
\end{equation*}
If \(R_k>2M\), then \(|X_k+\xi|/R_k\geq1/2\) on \(B_M\).
Consequently,
\begin{equation*}
 \|h\|_{H^1(B_M)}^2\leq C_a\|h\|_{X_k,a}^2.
\end{equation*}
It remains to test against compactly supported functions.

For \(G\)--\(G\), first suppose
\(\alpha_k=r_k/\widetilde r_k\to0\), and put
\(\beta_k=(c_k-\widetilde c_k)/\widetilde r_k^2\).
The pullback of the smaller profile in the larger frame is
\(D_{\alpha_k,\beta_k}f\).
Its support lies in \(\mathcal Q((0,\beta_k),C\alpha_k)\).
For compact smooth \(f,g\), Cauchy--Schwarz gives
\begin{equation*}
 |\langle D_{\alpha_k,\beta_k}f,g\rangle_{D_G^{1,2}}|
 \leq C\|f\|_{D_G^{1,2}}
 \left(\int_{\mathcal Q((0,\beta_k),C\alpha_k)}
 |\nabla_Gg|^2\right)^{1/2}.
\end{equation*}
The box measure is \(C\alpha_k^Q\to0\).
Absolute continuity of the integral makes the right-hand side tend to
zero uniformly in \(\beta_k\).
For the reverse pullback, unitarity gives
\begin{equation*}
 \langle D_{\alpha_k,\beta_k}^{-1}f,g\rangle_{D_G^{1,2}}
 =\langle f,D_{\alpha_k,\beta_k}g\rangle_{D_G^{1,2}}
 \longrightarrow0.
\end{equation*}
This also treats \(\alpha_k\to\infty\).
If \(\alpha_k\) stays in a compact subset of \((0,\infty)\),
\eqref{eq:singular-frame-orthogonality} gives
\begin{equation*}
 \frac{|c_k-\widetilde c_k|}{r_k^2+\widetilde r_k^2}\longrightarrow\infty.
\end{equation*}
The two fixed supports are then disjoint in the rescaled \(y\) variable.
Every subsequence has one of these further subsequences.
Thus the original sequence has weak limit zero.

For \(E\)--\(E\), first suppose \(|X_k-\widetilde X_k|\to\infty\).
The relative support leaves every bounded \(\xi\) set.
If \(|X_k-\widetilde X_k|\) is bounded, then
\begin{equation*}
 \left|R_k-\widetilde R_k\right|\leq|X_k-\widetilde X_k|,
 \qquad R_k/\widetilde R_k\longrightarrow1.
\end{equation*}
Equation \eqref{eq:escape-frame-orthogonality} now implies
\(|b_k-\widetilde b_k|/\widetilde R_k\to\infty\).
The relative support leaves every bounded \(\eta\) set.
The pullbacks vanish on each fixed ball for large \(k\).
The moving-norm bound gives local weak convergence.
Apply these alternatives to every subsequence to obtain the full limit.

For \(F\)--\(G\), let \(f,g\in C_c^\infty\).
If \(g\) is supported in \(|\xi|\leq C\), then
\(D_{r_k,c_k}g\) is supported in \(|x|\leq Cr_k\).
Thus
\begin{align*}
 |\langle T_{b_k}f,D_{r_k,c_k}g\rangle_{H_a}|
 &\leq
 \left(\int_{\{|x|\leq Cr_k\}}|\nabla_G f|^2\right)^{1/2}
 \|g\|_{D_G^{1,2}}+ar_k\|f\|_2\|g\|_2\\
 &\longrightarrow0.
\end{align*}
The first term tends to zero by absolute continuity.
The bound is independent of the relative \(y\) displacement.
The same estimate without the \(aL^2\) term gives
\begin{equation*}
 \langle D_{r_k,c_k}^{-1}T_{b_k}f,g\rangle_{D_G^{1,2}}
 =\langle T_{b_k}f,D_{r_k,c_k}g\rangle_{D_G^{1,2}}
 \longrightarrow0.
\end{equation*}
For \(F\)--\(E\), the original supports satisfy \(|x|\leq C\)
and \(|x-X_k|\leq C\).
Their energy pairing is zero for large \(k\).
The escape pullback of the finite profile vanishes on each fixed ball.
For \(G\)--\(E\), the original supports satisfy \(|x|\leq Cr_k\)
and \(|x-X_k|\leq C\).
In the singular receiver, the escape support has
\(|\xi|\geq(R_k-C)/r_k\to\infty\).
In the escape receiver, the singular support has
\(|X_k+\xi|\leq Cr_k\), and also leaves each fixed ball.
The uniform receiving bounds finish these cases.

To pass from smooth cores to model profiles, fix a compact core \(f_L\).
The realization estimates of Lemma \ref{lem:one-frame-limits} give
\begin{equation*}
 \limsup_k\|P_kf-P_{k,L}f\|_{H_a}
 \leq C\|f-f_L\|_{H_*}.
\end{equation*}
Here \(H_*\) is the source model space.
A fixed local test of any receiving pullback has error at most
\(C\|f-f_L\|_{H_*}\), by the receiving bounds above.
First let \(k\to\infty\) for the fixed core, then \(L\to\infty\).
\end{proof}

Let $P_k,Q_k$ be orthogonal frames and $f,g\in C_c^\infty$.
Lemma \ref{lem:relative-frames}, the coefficient limits in Lemma
\ref{lem:one-frame-limits}, and the vanishing singular $L^2$ term give
\begin{equation*}
 \langle P_kf,Q_kg\rangle_{H_a}\longrightarrow0.
\end{equation*}
For the mass interactions, put
\begin{equation*}
 I_k=\int\omega^{q/Q}|P_kf|^{q-1}|Q_kg|,\qquad
 J_k=\int\omega^{q/Q}|Q_kg|^{q-1}|P_kf|.
\end{equation*}
Lemma \ref{lem:one-frame-limits} bounds both weighted $L^q$ norms.  For $F$--$F$, $F$--$E$, $G$--$E$, and
$E$--$E$, the support separation proved above gives $I_k+J_k=0$ for
large $k$.

For \(F\)--\(G\), take \(P_kf\) to be the finite profile and
\(Q_kg\) to be the singular profile.  Put
\begin{equation*}
 B_k=\int_{\mathrm{supp}(Q_kg)}\omega^{q/Q}|P_kf|^q.
\end{equation*}
The singular support lies in \(|x|\leq Cr_k\).  Hence
\begin{equation*}
 B_k\leq\int_{|x|\leq Cr_k}\omega^{q/Q}|f|^q\longrightarrow0.
\end{equation*}
This bound is independent of the finite profile's \(y\) translation.
For \(G\)--\(G\), comparable scales give disjoint supports.
Otherwise exchange the frames so that \(Q_kg\) is the smaller profile.
In the larger frame its support lies in
\begin{equation*}
 \mathcal Q((0,\beta_k),C\alpha_k),\qquad
 \alpha_k=r_k/\widetilde r_k\to0,\qquad
 \beta_k=(c_k-\widetilde c_k)/\widetilde r_k^2.
\end{equation*}
The transformed weight is uniformly bounded on the fixed larger support.
The box measure is \(O(\alpha_k^Q)\).
Absolute continuity therefore gives \(B_k\to0\), uniformly in \(\beta_k\).
In both cases, H\"older's inequality gives
\begin{equation*}
 I_k+J_k\leq
 B_k^{(q-1)/q}\mathcal M(Q_kg)^{1/q}
 +B_k^{1/q}\mathcal M(Q_kg)^{(q-1)/q}\longrightarrow0.
\end{equation*}
Finally,
\begin{equation*}
 \bigl||s+t|^q-|s|^q-|t|^q\bigr|
 \leq C_q(|s|^{q-1}|t|+|t|^{q-1}|s|)
\end{equation*}
gives compact-profile mass decoupling in all six cases.
For nonsmooth compact cores, choose smooth approximations on a fixed
larger support with model errors at most \(\delta\).
The transformed coefficients are uniformly bounded there.
Their realized energy and weighted \(L^q\) errors are therefore at most
\(C\delta+o_k(1)\).
Cauchy--Schwarz controls the energy pairing errors.
The power inequality and H\"older's inequality control the mass errors
by \(C(\delta+\delta^{q-1})+o_k(1)\).
First let \(k\to\infty\), then \(\delta\downarrow0\).
The same argument permits a full finite profile in place of its compact
core, since translation preserves both approximation norms.

\begin{lemma}
\label{lem:synchronized-realizations}
Let a finite number of pairwise orthogonal frames and fixed profiles be
given.  Their singular and escape cutoff radii can be chosen so that
\begin{enumerate}[label=\textup{\(\roman*\)}]
\item every cutoff profile converges to the original profile in its model
energy and critical norm;
\item every coefficient convergence in Lemma \ref{lem:one-frame-limits} is
uniform on the corresponding support;
\item every pairwise energy and mass interaction tends to zero.
\end{enumerate}
For a countable family, one choice works at every fixed extraction level.
\end{lemma}

\begin{proof}
Use the block construction of Theorem \ref{thm:three-family-criterion}.
Number the profiles by \(j\).
At stage \(s\), take compact cores \(p_s^j\), \(j\leq s\), with
\begin{equation*}
 \|p_s^j-p^j\|_{H_{\tau_j}}
 +\mathcal M_{\tau_j}(p_s^j-p^j)^{1/q}\leq2^{-s}.
\end{equation*}
Their cutoff radii can be chosen at least \(s\) and increasing in \(s\).
For finite profiles use the full translated profile.
Its compact-core error is unchanged by translation.

For these finitely many cores, choose \(k_s>k_{s-1}\).
For every \(k\geq k_s\), require each self-energy and self-mass error,
each pairwise energy and mass error, and the first \(s\) local-test
errors to be at most \(2^{-s}\).
The fixed-core limits proved above permit this choice.
For a singular core of radius \(L\), also require \(r_kL<2^{-s}\).
For an escape core require \(L/|X_k|<2^{-s}\).
Include the three errors
\begin{align*}
 &\sup_{|\xi|\leq CL}
 \left|(1+r_k^2|\xi|^2)^{n/(Q-2)}-1\right|,\\
 &\sup_{|\xi|\leq2L}
 \left|\frac{|X_k+\xi|^2}{|X_k|^2}-1\right|,\\
 &\sup_{|\xi|\leq2L}
 \left||X_k|^{-2n/(Q-2)}
 (1+|X_k+\xi|^2)^{n/(Q-2)}-1\right|
\end{align*}
in the same finite list.
For fixed \(s\), all these errors tend to zero with \(k\).
Use the stage-\(s\) cores on \(k_s\leq k<k_{s+1}\).

Fix \(j\).  For large \(s\), its core is present and its model error is
at most \(2^{-s}\).
Lemma \ref{lem:one-frame-limits} therefore gives its energy and mass limits.
The displayed coefficient errors tend to zero on its selected support.
Fix also \(i\ne j\).
Their interaction errors are included once \(s\geq\max\{i,j\}\),
and hence tend to zero.  This proves \textup{(i)}--\textup{(iii)}.
The same construction works for a countable list because each stage
contains only finitely many requirements.

Let \(P_kp\) and \(\widetilde P_kp\) be two such realizations.
Write their model cutoff cores as \(p_k\) and \(\widetilde p_k\).
Then
\begin{equation*}
 \|p_k-\widetilde p_k\|_{H_*}
 +\mathcal M_*(p_k-\widetilde p_k)^{1/q}\longrightarrow0.
\end{equation*}
For a singular profile put
\(L_k^*=\max\{L_k,\widetilde L_k\}\).
The extra zeroth-order term satisfies
\begin{equation*}
 r_k^2\|p_k-\widetilde p_k\|_2^2
 \leq C(r_kL_k^*)^2\|p\|_q^2\longrightarrow0.
\end{equation*}
The singular weight tends to one uniformly on this union of supports.
For an escape profile, both transformed coefficients tend to one
uniformly on the union because \(L_k^*/R_k\to0\).
Consequently, in either case,
\begin{align*}
 \mathcal E_a(P_kp-\widetilde P_kp)
 &\leq C\mathcal E_*(p_k-\widetilde p_k)+o_k(1)\longrightarrow0,\\
 \mathcal M(P_kp-\widetilde P_kp)
 &\leq C\mathcal M_*(p_k-\widetilde p_k)\longrightarrow0.
\end{align*}
Finite profiles have the same full translation in both realizations.
Finally, include the first \(s\) tests from each of the first \(s\)
old orthogonal frames.
Every fixed test is then included for all sufficiently large stages.
Its limit remains zero.
\end{proof}

\begin{lemma}
\label{lem:pairwise-decoupling}
Let $P_k$ and $Q_k$ be two orthogonal frames.  Let $f$ and $g$ be fixed
profiles of the corresponding types.  Then
\begin{equation}
 \langle P_kf,Q_kg\rangle_{H_a}\longrightarrow0.
 \label{eq:pairwise-energy}
\end{equation}
Moreover,
\begin{equation}
 \int\omega^{q/Q}\left(
 |P_kf+Q_kg|^q-|P_kf|^q-|Q_kg|^q\right)\longrightarrow0.
 \label{eq:pairwise-mass}
\end{equation}
\end{lemma}

\begin{proof}
Let \(\tau_1,\tau_2\) be the two model types.
For \(\delta>0\), choose compact cores \(f_0,g_0\) with
\begin{equation*}
 \|f-f_0\|_{H_{\tau_1}}+\|g-g_0\|_{H_{\tau_2}}\leq\delta.
\end{equation*}
Lemma \ref{lem:one-frame-limits} and the realization estimates give
\begin{equation*}
 \limsup_k\left(
 \|P_kf-P_kf_0\|_{H_a}+\|Q_kg-Q_kg_0\|_{H_a}\right)
 \leq C\delta.
\end{equation*}
All four realized norms are bounded by a constant depending on the
fixed profile norms and \(\delta\leq1\).
Hence
\begin{align*}
 &\left|\langle P_kf,Q_kg\rangle_{H_a}
 -\langle P_kf_0,Q_kg_0\rangle_{H_a}\right|\\
 &\qquad\leq
 \|P_kf-P_kf_0\|_{H_a}\|Q_kg\|_{H_a}
 +\|P_kf_0\|_{H_a}\|Q_kg-Q_kg_0\|_{H_a}.
\end{align*}
The compact-core pairing tends to zero.
Take \(\limsup_k\), then \(\delta\downarrow0\), to obtain \eqref{eq:pairwise-energy}.

Write \(\|v\|_{q,\omega}=\mathcal M(v)^{1/q}\).
The global Sobolev inequality bounds the two approximation errors in
this norm by \(C\delta\) in the limit superior.
For any weighted \(L^q\) functions \(F,G\), the power inequality and
H\"older's inequality give
\begin{equation*}
 \left|\mathcal M(F)-\mathcal M(G)\right|
 \leq C_q\left(\|F\|_{q,\omega}+\|G\|_{q,\omega}\right)^{q-1}
 \|F-G\|_{q,\omega}.
\end{equation*}
Apply this estimate to the sum and to its two individual terms.
The mass interaction changes by at most \(C\delta\) in the limit superior
when \((f,g)\) is replaced by \((f_0,g_0)\).
The compact-core interaction tends to zero.
Let \(\delta\downarrow0\) to obtain \eqref{eq:pairwise-mass}.
\end{proof}

\begin{lemma}
\label{lem:one-profile-subtraction}
Let \((u_k)\) be bounded in \(H_a\).  Suppose that the pullback by a finite,
singular, or escape frame converges weakly to a nonzero profile \(p\) in the
corresponding model topology.  After a subsequence, there are synchronized
cutoff radii \(L_k\to\infty\) and realizations \(R_kp\in H_a\) such that
\begin{align}
 \mathcal E_a(u_k)
 &=\mathcal E_a(u_k-R_kp)+\mathcal E_*(p)+o_k(1),
 \label{eq:one-profile-energy-splitting}\\
 \mathcal M(u_k)
 &=\mathcal M(u_k-R_kp)+\mathcal M_*(p)+o_k(1),
 \label{eq:one-profile-mass-splitting}
\end{align}
and the pullback of \(u_k-R_kp\) is weakly zero in that frame.  If the
original sequence was weakly zero in finitely many old orthogonal frames,
then the new remainder is weakly zero in those frames as well.
\end{lemma}

\begin{proof}
Fix \(L\).  Put \(p_L=p\) in a finite frame,
\(p_L=\chi_Lp\) in a singular frame, and
\(p_L=\vartheta(|\cdot|/L)p\) in an escape frame.
Write \(v_k\) for the receiving pullback of \(u_k\).
In the finite case, global weak convergence gives the pairing limit.
In the singular case,
\begin{equation*}
 \langle u_k,D_{r_k,c_k}p_L\rangle_{H_a}
 =\langle v_k,p_L\rangle_{D_G^{1,2}}
 +a\int u_kD_{r_k,c_k}p_L.
\end{equation*}
The last term satisfies
\begin{equation*}
 \left|a\int u_kD_{r_k,c_k}p_L\right|
 \leq ar_k\|u_k\|_2\|p_L\|_2\longrightarrow0.
\end{equation*}
In an escape frame, let \(K\) contain the support of \(p_L\).
Then
\begin{align*}
 \langle u_k,A_{X_k,b_k}p_L\rangle_{H_a}
 &=\int_K\left(\nabla_\xi v_k\cdot\nabla_\xi p_L
 +c_k\nabla_\eta v_k\cdot\nabla_\eta p_L+av_kp_L\right),\\
 c_k(\xi)&=|X_k+\xi|^2/|X_k|^2.
\end{align*}
The local \(H^1\) norms of \(v_k\) are bounded, and
\(\sup_K|c_k-1|\to0\).
The coefficient error is at most
\(\sup_K|c_k-1|\,\|\nabla_\eta v_k\|_{L^2(K)}
\|\nabla_\eta p_L\|_2\to0\).
Hence
\begin{equation}
 \langle u_k,P_{k,L}p\rangle_{H_a}
 \longrightarrow\langle p,p_L\rangle_{H_*}.
 \label{eq:fixed-cutoff-pairing}
\end{equation}
They also apply to every fixed compact model test.
Polarization and Lemma \ref{lem:one-frame-limits} yield
\begin{align}
 \mathcal E_a(u_k-P_{k,L}p)
 =\mathcal E_a(u_k)-2\langle p,p_L\rangle_{H_*}
 +\mathcal E_*(p_L)+o_k(1).
 \label{eq:fixed-cutoff-energy-subtraction}
\end{align}

Apply Lemma \ref{lem:uniform-local-rellich} on model boxes \(K_s\)
increasing to the whole space.
Choose nested subsequences with strong \(L^2(K_s)\) convergence.
Choose the \(k\)-th diagonal term far enough along its subsequence that
\begin{equation*}
 \sum_{s=1}^k\|v_k-p\|_{L^2(K_s)}^2\leq2^{-k}.
\end{equation*}
For each fixed \(s\), the series
\(\sum_k\|v_k-p\|_{L^2(K_s)}^2\) is finite.
Hence \(v_k\to p\) almost everywhere on every \(K_s\).
For a finite frame, apply Lemma \ref{lem:weighted-bl} in the fixed measure.
For the other two frames, take a compact set \(K\) containing
the support of \(p_L\).
The pulled-back density \(w_k\) satisfies
\begin{equation*}
 \sup_K|w_k-1|\longrightarrow0,\qquad
 \sup_k\int_K\left(|v_k|^q+|v_k-p_L|^q\right)<\infty.
\end{equation*}
Apply the Brezis--Lieb lemma twice on \(K\):
\begin{align*}
 \int_K\left(|v_k|^q-|v_k-p|^q\right)&\longrightarrow\int_K|p|^q,\\
 \int_K\left(|v_k-p_L|^q-|v_k-p|^q\right)
 &\longrightarrow\int_K|p-p_L|^q.
\end{align*}
Subtract these limits.
Replacing \(1\) by \(w_k\) changes their difference by at most
\begin{equation*}
 \sup_K|w_k-1|\int_K\left(|v_k|^q+|v_k-p_L|^q\right)
 \longrightarrow0.
\end{equation*}
Outside \(K\), the subtraction changes nothing.
In all three cases,
\begin{equation}
 \mathcal M(u_k)-\mathcal M(u_k-P_{k,L}p)
 =\mathcal M_*(p)-\mathcal M_*(p-p_L)+o_k(1).
 \label{eq:fixed-cutoff-mass-subtraction}
\end{equation}
Moreover,
\begin{align*}
 &\left|\mathcal M_*(p)-\mathcal M_*(p-p_L)-\mathcal M_*(p_L)\right|\\
 &\qquad\leq C\left(
 \|p_L\|_q^{q-1}\|p-p_L\|_q
 +\|p-p_L\|_q^{q-1}\|p_L\|_q\right)
 \longrightarrow0\quad\hbox{as }L\to\infty.
\end{align*}
The norms use the corresponding model measure.
For a finite profile the error is zero.

Choose \(L_s\to\infty\) with model energy- and mass-norm errors at most
\(2^{-s}\).
For fixed \(s\), choose \(k_s>k_{s-1}\) so that, for \(k\geq k_s\),
the errors in \eqref{eq:fixed-cutoff-pairing},
\eqref{eq:fixed-cutoff-energy-subtraction}, and
\eqref{eq:fixed-cutoff-mass-subtraction},
the self-energy and self-mass limits, and the first \(s\) local weak tests
relative to their fixed-\(L_s\) limits are at most \(2^{-s}\).
Impose also the coefficient and scale constraints in
Lemma \ref{lem:synchronized-realizations}.
This is the block construction of Theorem \ref{thm:three-family-criterion}.
For each of the finitely many old orthogonal pullbacks \(S_k^i\),
Lemma \ref{lem:relative-frames} gives
\begin{equation*}
 S_k^iP_{k,L_s}p\rightharpoonup0\quad\hbox{locally}
 \quad\hbox{for fixed }s.
\end{equation*}
Include its first \(s\) tests in the same finite list before choosing
\(k_s\).
Use \(L_k=L_s\) on \(k_s\leq k<k_{s+1}\).
Then
\begin{equation*}
 \langle p,p_{L_k}\rangle_{H_*}\to\mathcal E_*(p),\qquad
 \mathcal E_*(p_{L_k})\to\mathcal E_*(p),\qquad
 \mathcal M_*(p-p_{L_k})\to0.
\end{equation*}
The two subtraction formulas give
\eqref{eq:one-profile-energy-splitting} and
\eqref{eq:one-profile-mass-splitting}.
In the new frame, each fixed weak test of the remainder has limit
that of \(p-p_{L_s}\), which tends to zero as \(s\to\infty\).
In every old frame, both the original sequence and the selected
realization have zero weak limit.
For a finite profile use the full translation throughout;
no cutoff selection is needed.
\end{proof}

\section{A three-scale vanishing lemma}

\subsection{Scale-invariant local quantities}

Fix $\varepsilon=1/40$ throughout this section.  For $u\in H_a$ define
\begin{align}
 \Theta_F(u)
 &=\sup_{b\in\mathbb R^n}
 \int_{\{|x|<2\varepsilon^{-1},\ |y-b|<2\varepsilon^{-1}\}}|u|^2,
 \label{eq:theta-f}\\
 \Theta_G(u)
 &=\sup_{0<r<16\varepsilon}\sup_{c\in\mathbb R^n}
 r^{-2}\int_{\mathcal Q((0,c),r)}|u|^2,
 \label{eq:theta-g}\\
 \Theta_E(u)
 &=\sup_{|X|>(2\varepsilon)^{-1}}\sup_{b\in\mathbb R^n}
 \int_{\{|x-X|<2\varepsilon^{-1},\
 |y-b|<2\varepsilon^{-1}|X|\}}|u|^2.
 \label{eq:theta-e}
\end{align}
Each gauge is homogeneous of degree two.  The factor \(r^{-2}\) in
\eqref{eq:theta-g} comes from the singular scaling; the escape change
of variables preserves local \(L^2\) mass.

For the finite and escape gauges,
\begin{equation*}
 \Theta_F(u),\Theta_E(u)\leq\|u\|_2^2\leq a^{-1}\mathcal E_a(u),
\end{equation*}
and H\"older's inequality and the homogeneous Sobolev inequality give
\begin{equation*}
 r^{-2}\int_{\mathcal Q((0,c),r)}|u|^2
 \leq r^{-2}|\mathcal Q((0,c),r)|^{2/Q}\|u\|_q^2
 \leq C\mathcal E_G(u)\leq C\mathcal E_a(u).
\end{equation*}
The triangle inequality for the local $L^2$ norms gives
\begin{equation*}
 \left|\Theta_\tau(u)^{1/2}-\Theta_\tau(v)^{1/2}\right|
 \leq \Theta_\tau(u-v)^{1/2}
 \leq C\mathcal E_a(u-v)^{1/2},\qquad \tau\in\{F,G,E\}.
\end{equation*}

We first estimate functions with zero weighted mean in \(\xi\).

\begin{lemma}
\label{lem:discrete-interpolation}
Let \(K_\xi\subset\mathbb R^m\) be bounded, connected, and Lipschitz.
Fix \(0\leq\zeta_K\in C_c^\infty(K_\xi)\) with \(\int\zeta_K=1\).
Suppose \(h\in H^1(K_\xi\times\mathbb R^n)\) and
\begin{equation*}
 \int_{K_\xi}\zeta_K(\xi)h(\xi,\eta)\,d\xi=0
 \quad\text{for almost every }\eta.
\end{equation*}
For \(0<t<1\), let
\begin{equation*}
 E=\int(|\nabla_\xi h|^2+|\nabla_\eta h|^2+at^2|h|^2),\qquad
 \mu=\sup_{B\in\mathcal B_0}\int_{K_\xi\times2B}|h|^2,
\end{equation*}
where \(\mathcal B_0\) is the unit dyadic cube family in \(\mathbb R^n\).
If \(E\leq A\), then, for some \(\delta>0\),
\begin{equation}
 \|h\|_q^q\leq C_A E\mu^\delta.
 \label{eq:mean-zero-discrete}
\end{equation}
The constants depend on the fixed set and averaging weight, and are
independent of \(t\).
\end{lemma}

\begin{proof}
Since \(Q>N\), we have \(q<2N/(N-2)\) when \(N\geq3\).
Choose \(p>q\) below this exponent.  For \(N=2\), choose a finite
\(p>q\).  Take \(p\) sufficiently close to \(q\) that
\begin{equation*}
 \frac1q=\frac{1-\theta}{2}+\frac\theta p,\qquad
 \delta=\frac{(1-\theta)q}{2}>0,\qquad
 \sigma=\frac{\theta q}{2}-1>0.
\end{equation*}
The weighted-mean Poincar\'e inequality on \(K_\xi\) gives
\begin{equation*}
 \int_{K_\xi}|h(\xi,\eta)|^2\,d\xi
 \leq C\int_{K_\xi}|\nabla_\xi h(\xi,\eta)|^2\,d\xi.
\end{equation*}
For \(B\in\mathcal B_0\), put
\begin{equation*}
 E_B=\int_{K_\xi\times2B}
 (|\nabla_\xi h|^2+|\nabla_\eta h|^2+at^2|h|^2).
\end{equation*}
After integration over \(2B\), local Sobolev gives
\begin{equation*}
 \|h\|_{L^p(K_\xi\times B)}^2
 \leq C\int_{K_\xi\times2B}
 (|\nabla_\xi h|^2+|\nabla_\eta h|^2+|h|^2)
 \leq CE_B.
\end{equation*}
Uniformly in \(t\), interpolation gives
\begin{align*}
 \int_{K_\xi\times B}|h|^q
 &\leq
 \left(\int_{K_\xi\times B}|h|^2\right)^{(1-\theta)q/2}
 \left(\int_{K_\xi\times B}|h|^p\right)^{\theta q/p}\\
 &\leq C\mu^\delta E_B^{1+\sigma}.
\end{align*}
The cubes partition \(\mathbb R^n\), with bounded overlap of their doubles:
\begin{equation*}
 \sum_BE_B\leq CE,\qquad
 \sum_BE_B^{1+\sigma}\leq A^\sigma\sum_BE_B\leq CA^\sigma E.
\end{equation*}
Summing proves \eqref{eq:mean-zero-discrete}.
\end{proof}

\subsection{A contracting recurrence for transverse averages}

\begin{lemma}
\label{lem:vertical-high-frequency}
Set \(Tf(\eta)=2^{-(Q-2)/2}f(\eta/4)\).
Let \(g_j,d_j\in H^1(\mathbb R^n)\) satisfy
\begin{equation}
 g_{j+1}=T(g_j+d_j),\qquad j\geq0.
 \label{eq:abstract-terminal-recurrence}
\end{equation}
Put
\begin{align*}
 E_{\mathrm{rec}}&=\|g_0\|_{H^1}^2+
 \sum_{j\geq0}\|d_j\|_{H^1}^2,\\
 \nu_{\mathrm{rec}}&=\sup_{B\in\mathcal B_0}\int_{2B}|g_0|^2+
 \sup_{j\geq0}\sup_{B\in\mathcal B_0}\int_{2B}|d_j|^2.
\end{align*}
If \(E_{\mathrm{rec}}\leq L\), there is \(\delta_1>0\) such that
\begin{equation}
 \sum_{j\geq0}\|g_j\|_q^q
 \leq C_L\nu_{\mathrm{rec}}^{\delta_1}E_{\mathrm{rec}}.
 \label{eq:vertical-high-frequency-bound}
\end{equation}
The estimate is uniform over finite truncations of the recurrence.
\end{lemma}

\begin{proof}
If \(n\geq3\), then \(Q>n\) implies \(q<2n/(n-2)\).
If \(n=1,2\), every finite local Sobolev exponent is available.
Choose \(p>q\) sufficiently close to \(q\) that
\begin{equation*}
 \frac1q=\frac{1-\theta}{2}+\frac\theta p,\qquad
 \delta_1=\frac{(1-\theta)q}{2}>0,\qquad
 \rho=\frac{\theta q}{2}=\frac q2-\delta_1>1.
\end{equation*}
Let \((h_\ell)\subset H^1(\mathbb R^n)\) satisfy
\begin{equation*}
 \sum_\ell\|h_\ell\|_{H^1}^2\leq L,\qquad
 \sup_{\ell,B}\int_{2B}|h_\ell|^2\leq\nu.
\end{equation*}
Set
\begin{equation*}
 e_{\ell B}=\int_{2B}(|\nabla h_\ell|^2+|h_\ell|^2).
\end{equation*}
Local Sobolev interpolation gives
\begin{equation*}
 \int_B|h_\ell|^q\leq C\nu^{\delta_1}e_{\ell B}^{\rho}.
\end{equation*}
Since \(e_{\ell B}\leq L\) and the doubled cubes have bounded overlap,
\begin{align*}
 \sum_\ell\|h_\ell\|_q^q
 &\leq C\nu^{\delta_1}\sum_{\ell,B}e_{\ell B}^{\rho}\\
 &\leq C\nu^{\delta_1}L^{\rho-1}
 \sum_{\ell,B}e_{\ell B}
 \leq C_L\nu^{\delta_1}\sum_\ell\|h_\ell\|_{H^1}^2.
\end{align*}
For \(g_0,d_0,d_1,\ldots\), the two bounds are
\(E_{\mathrm{rec}}\) and \(\nu_{\mathrm{rec}}\).  Hence
\begin{equation*}
 \|g_0\|_q^q+\sum_{j\geq0}\|d_j\|_q^q
 \leq C_L\nu_{\mathrm{rec}}^{\delta_1}E_{\mathrm{rec}}.
\end{equation*}
A change of variables gives
\begin{equation*}
 \|Tf\|_q^q
 =2^{-q(Q-2)/2}4^n\|f\|_q^q
 =2^{-Q+2n}\|f\|_q^q
 =2^{-m}\|f\|_q^q.
\end{equation*}
Put \(\lambda=2^{-m/q}<1\).  Iterating
\eqref{eq:abstract-terminal-recurrence}, we obtain
\begin{equation*}
 g_j=T^jg_0+\sum_{\ell=0}^{j-1}T^{j-\ell}d_\ell,\qquad
 \|g_j\|_q\leq\lambda^j\|g_0\|_q+
 \sum_{\ell=0}^{j-1}\lambda^{j-\ell}\|d_\ell\|_q.
\end{equation*}
Discrete Young's inequality on \(0\leq j\leq J\) yields
\begin{equation*}
 \left(\sum_{j=0}^J\|g_j\|_q^q\right)^{1/q}
 \leq\frac1{1-\lambda}
 \left(\|g_0\|_q^q+\sum_{\ell=0}^{J-1}\|d_\ell\|_q^q\right)^{1/q}.
\end{equation*}
Raise to the power \(q\) and let \(J\to\infty\) to obtain
\eqref{eq:vertical-high-frequency-bound}.
For \(m=1\), no \(L^2\) contraction is used: only the defects
\(d_j\) and \(g_0\) have summed squared \(H^1\) norms.
\end{proof}

\subsection{Annular control near the degeneracy set}

\begin{lemma}
\label{lem:singular-morrey}
For every \(L>0\), there are \(C_L>0\) and \(\gamma>0\) such that
\begin{equation}
 \int_{|x|<2\varepsilon}|v|^q
 \leq C_L\Theta_G(v)^\gamma\mathcal E_a(v),
 \qquad \mathcal E_a(v)\leq L.
 \label{eq:singular-morrey}
\end{equation}
\end{lemma}

\begin{proof}
First let \(v\in C_c^\infty\).  Put \(t_j=2^{1-j}\varepsilon\) and
\begin{align*}
 \mathcal A_j&=\{t_j/2<|x|<t_j\},\qquad
 \Omega=\{1/2<|\xi|<1\},\\
 \Omega^*&=\{1/4<|\xi|<2\},\qquad
 \Omega^{**}=\{1/8<|\xi|<4\}.
\end{align*}
Define
\begin{equation}
 h_j(\xi,\eta)=t_j^{(Q-2)/2}v(t_j\xi,t_j^2\eta),\qquad
 E_j=\int_{\Omega^{**}\times\mathbb R^n}
 (|\nabla_\xi h_j|^2+|\nabla_\eta h_j|^2+at_j^2|h_j|^2).
 \label{eq:annular-rescaling}
\end{equation}
The substitution in \eqref{eq:annular-rescaling} gives
\begin{align*}
 \int_{\Omega^{**}\times\mathbb R^n}|\nabla_\xi h_j|^2
 &=\int_{t_j\Omega^{**}\times\mathbb R^n}|\nabla_xv|^2,\\
 \int_{\Omega^{**}\times\mathbb R^n}|\nabla_\eta h_j|^2
 &=t_j^2\int_{t_j\Omega^{**}\times\mathbb R^n}|\nabla_yv|^2,\\
 at_j^2\int_{\Omega^{**}\times\mathbb R^n}|h_j|^2
 &=a\int_{t_j\Omega^{**}\times\mathbb R^n}|v|^2.
\end{align*}
Since \(t_j/8<|x|<4t_j\) on \(t_j\Omega^{**}\), bounded overlap gives
\begin{equation}
 \sum_jE_j\leq C\mathcal E_a(v),\qquad
 \int_{\Omega\times\mathbb R^n}|h_j|^q=\int_{\mathcal A_j}|v|^q.
 \label{eq:annular-rescaled-mass}
\end{equation}
For \(m=1\), work on one half-line at a time, intersecting all three
annuli with it.  The resulting domains are connected.

Choose \(0\leq\zeta\in C_c^\infty(\Omega^*)\) with
\(\int\zeta=1\).  On each half-line this is a separate normalized
weight.  Set
\begin{equation*}
 g_j(\eta)=\int\zeta(\xi)h_j(\xi,\eta)\,d\xi,\qquad
 k_j=h_j-g_j,\qquad
 \mu=\sup_{j,B\in\mathcal B_0}
 \int_{\Omega^{**}\times2B}|h_j|^2.
\end{equation*}
There are \(M_n\) points \(\eta_{B,s}\), with \(M_n\) depending
only on \(n\), such that
\begin{equation*}
 t_j\Omega^{**}\times t_j^2(2B)
 \subset\mathcal Q((0,t_j^2\eta_{B,1}),4t_j)\cup\ldots\cup
 \mathcal Q((0,t_j^2\eta_{B,M_n}),4t_j).
\end{equation*}
The \(x\) radius is \(4t_j\); the transverse cube of side
\(2t_j^2\) requires only a dimension-dependent number of balls
of radius \(16t_j^2\).
Since \(4t_j\leq8\varepsilon<16\varepsilon\),
\begin{align*}
 \int_{\Omega^{**}\times2B}|h_j|^2
 &=t_j^{-2}\int_{t_j\Omega^{**}\times t_j^2(2B)}|v|^2\\
 &\leq16M_n\Theta_G(v).
\end{align*}
Thus \(\mu\leq C\Theta_G(v)\leq C_L\).

The weighted mean of \(k_j\) is zero.  Poincar\'e's inequality and
Jensen's inequality give
\begin{align*}
 \int_{\Omega^{**}}|k_j|^2
 &\leq C\int_{\Omega^{**}}|\nabla_\xi h_j|^2,\\
 |\nabla_\eta g_j|^2
 &\leq\int\zeta|\nabla_\eta h_j|^2,\\
 \int_{\Omega^{**}}|\nabla_\eta k_j|^2
 &\leq C\int_{\Omega^{**}}|\nabla_\eta h_j|^2.
\end{align*}
Jensen also gives
\begin{equation*}
 \sup_B\int_{\Omega^{**}\times2B}|k_j|^2\leq C\mu.
\end{equation*}
After integration in \(\eta\), the energy of \(k_j\) is at most
\(CE_j\leq CL\).  Lemma \ref{lem:discrete-interpolation} yields
\begin{equation}
 \sum_j\|k_j\|_{L^q(\Omega\times\mathbb R^n)}^q
 \leq C_L\mu^{\delta_0}\sum_jE_j
 \label{eq:annular-mean-zero-bound}
\end{equation}
for some \(\delta_0>0\).

Put \(\zeta_1(\xi)=2^m\zeta(2\xi)\) and
\begin{equation*}
 d_j(\eta)=\int(\zeta_1-\zeta)(\xi)h_j(\xi,\eta)\,d\xi.
\end{equation*}
Both weights are supported in \(\Omega^{**}\), and
\(\int(\zeta_1-\zeta)=0\).  Hence
\begin{equation*}
 d_j=\int(\zeta_1-\zeta)(h_j-g_j).
\end{equation*}
Poincar\'e controls this expression; Jensen controls its
\(\eta\)-derivatives and local mass:
\begin{equation}
 \|d_j\|_2^2\leq C\int|\nabla_\xi h_j|^2,
 \qquad \|\nabla_\eta d_j\|_2^2\leq C\int|\nabla_\eta h_j|^2,
 \qquad \sup_B\int_{2B}|d_j|^2\leq C\mu.
 \label{eq:annular-mean-defect-energy}
\end{equation}
The right-hand integrals are over \(\Omega^{**}\times\mathbb R^n\).
At \(j=0\), Jensen gives
\begin{equation*}
 \|g_0\|_2^2\leq C\int|h_0|^2
 \leq\frac C{at_0^2}E_0,\qquad
 \|\nabla_\eta g_0\|_2^2\leq CE_0.
\end{equation*}
Since \(t_0=2\varepsilon\) is fixed,
\eqref{eq:annular-mean-defect-energy} gives
\begin{equation*}
 \|g_0\|_{H^1}^2+\sum_j\|d_j\|_{H^1}^2
 \leq C_{a,\varepsilon}\sum_jE_j,\qquad
 \sup_B\int_{2B}|g_0|^2\leq C\mu.
\end{equation*}
Moreover,
\begin{align*}
 h_{j+1}(\xi,\eta)
 &=2^{-(Q-2)/2}h_j(\xi/2,\eta/4),\\
 g_{j+1}(\eta)
 &=2^{-(Q-2)/2}\int2^m\zeta(2\xi)
 h_j(\xi,\eta/4)\,d\xi.
\end{align*}
Therefore
\begin{equation}
 g_{j+1}(\eta)=2^{-(Q-2)/2}(g_j+d_j)(\eta/4).
 \label{eq:exact-mean-chain}
\end{equation}
Apply Lemma \ref{lem:vertical-high-frequency} to
\eqref{eq:exact-mean-chain}, with
\(E_{\mathrm{rec}}\leq C_{a,\varepsilon}L\) and
\(\nu_{\mathrm{rec}}\leq C\mu\).  We obtain
\begin{equation*}
 \sum_j\|g_j\|_q^q\leq C_L\mu^{\delta_1}\sum_jE_j.
\end{equation*}
Since \(|h_j|^q\leq2^{q-1}(|k_j|^q+|g_j|^q)\),
\eqref{eq:annular-mean-zero-bound} implies
\begin{equation*}
 \sum_j\int_{\Omega\times\mathbb R^n}|h_j|^q
 \leq C_L(\mu^{\delta_0}+\mu^{\delta_1})\sum_jE_j.
\end{equation*}
Take \(\gamma=\min\{\delta_0,\delta_1\}\).
The bound \(\mu\leq C_L\) gives
\begin{equation}
 \sum_j\int_{\Omega\times\mathbb R^n}|h_j|^q
 \leq C_L\Theta_G(v)^\gamma\mathcal E_a(v).
 \label{eq:global-annular-bridge}
\end{equation}
Add the half-line estimates when \(m=1\).  The annuli partition
\(0<|x|<2\varepsilon\), so \eqref{eq:annular-rescaled-mass} proves
\eqref{eq:singular-morrey} for smooth \(v\).

Now take \(v_\ell\in C_c^\infty\) with \(v_\ell\to v\) in \(H_a\).
For large \(\ell\), \(\mathcal E_a(v_\ell)\leq L+1\).
Proposition \ref{prop:embedding} gives
\begin{equation*}
 \int\omega^{q/Q}|v_\ell-v|^q\longrightarrow0,
 \qquad
 \|v_\ell-v\|_{L^q(\{|x|<2\varepsilon\})}\longrightarrow0.
\end{equation*}
The gauge estimate gives
\begin{equation*}
 |\Theta_G(v_\ell)^{1/2}-\Theta_G(v)^{1/2}|
 \leq C\mathcal E_a(v_\ell-v)^{1/2}\longrightarrow0.
\end{equation*}
Pass to the limit in \eqref{eq:global-annular-bridge}, absorbing
the dependence on \(L+1\) into \(C_L\).
\end{proof}

\subsection{Inverse vanishing and profile extraction}

\begin{proposition}
\label{prop:refined}
Let $(u_k)$ be bounded in $H_a$.  If
\begin{equation*}
 \Theta_F(u_k)+\Theta_G(u_k)+\Theta_E(u_k)\longrightarrow0,
\end{equation*}
then
\begin{equation}
 \mathcal M(u_k)\longrightarrow0.
 \label{eq:mass-vanishing}
\end{equation}
There are \(C>0\) and \(0<\gamma<q/2\), depending only on \(m,n,a\), such that
\begin{equation}
 \mathcal M(u)
 \leq C\mathcal E_a(u)^{q/2-\gamma}
 \left(\max\{\Theta_F(u),\Theta_G(u),\Theta_E(u)\}\right)^\gamma.
 \label{eq:refined-global}
\end{equation}
\end{proposition}

\begin{proof}
First assume \(\mathcal E_a(u)\leq1\).  All constants below depend
only on \(m,n,a\), since \(\varepsilon=1/40\) is fixed.
Lemma \ref{lem:singular-morrey} gives
\begin{equation}
 \int_{|x|<2\varepsilon}|u|^q\leq C\Theta_G(u)^\gamma.
 \label{eq:singular-vanishing-bound}
\end{equation}
The weight is bounded on this slab.

For the remaining regions, use the following Euclidean estimate.
Let \(2<p<2_N^*\), with any finite upper exponent when \(N=2\).
Suppose that \(h_i\in H^1(2Q_i)\), where \(Q_i\) are unit cubes, and
\begin{equation*}
 e_i=\int_{2Q_i}(|\nabla h_i|^2+|h_i|^2),\qquad
 \sum_ie_i\leq L,\qquad
 \sup_i\int_{2Q_i}|h_i|^2\leq\tau.
\end{equation*}
Choose \(s>p\) below the local Sobolev exponent.  Take \(s\) close
enough to \(p\) that
\begin{equation*}
 \frac1p=\frac{1-\theta}{2}+\frac\theta s,\qquad
 \gamma_p=\frac{(1-\theta)p}{2}>0,\qquad
 \rho_p=\frac{\theta p}{2}>1.
\end{equation*}
Sobolev interpolation gives
\begin{equation*}
 \int_{Q_i}|h_i|^p\leq C\tau^{\gamma_p}e_i^{\rho_p}.
\end{equation*}
Since \(e_i\leq L\), summation yields
\begin{equation}
 \sum_i\int_{Q_i}|h_i|^p
 \leq C\tau^{\gamma_p}L^{\rho_p-1}\sum_ie_i
 \leq C_L\tau^{\gamma_p}.
 \label{eq:euclidean-lions}
\end{equation}
The estimate also holds on fixed rescalings of products of balls,
uniformly under translation.

Cover \(2\varepsilon<|x|<\varepsilon^{-1}\) by finitely many
\(x\) balls of radius at most \(\varepsilon/4\).  Choose the centers
in this strip.  Their doubles lie in
\(\varepsilon<|x|<2\varepsilon^{-1}\).  In \(y\), use unit cubes.
After fixed rescalings in \(x\), denote the restrictions of \(u\)
by \(h_i\).  On the doubled sets,
\begin{equation*}
 \varepsilon^2|\nabla_yu|^2\leq|x|^2|\nabla_yu|^2,
 \qquad |u|^2\leq a^{-1}a|u|^2.
\end{equation*}
Bounded overlap therefore gives
\begin{equation*}
 \sum_i\int_{2Q_i}(|\nabla h_i|^2+|h_i|^2)
 \leq C\mathcal E_a(u),\qquad
 \sup_i\int_{2Q_i}|h_i|^2\leq C\Theta_F(u).
\end{equation*}
For the second bound, cover each doubled set by a dimension-dependent
number of boxes from \eqref{eq:theta-f}.  The weight is bounded here.
Apply \eqref{eq:euclidean-lions} with \(p=q\).  We obtain
\begin{equation}
 \int_{2\varepsilon<|x|<\varepsilon^{-1}}
 \omega^{q/Q}|u|^q\leq C\Theta_F(u)^\gamma.
 \label{eq:finite-vanishing-bound}
\end{equation}

For \(\ell\geq5\), put \(R_\ell=2^\ell\).
Cover \(R_\ell<|x|<2R_\ell\) by \(x\) balls of radius one and
\(y\) balls of radius \(R_\ell\).  Choose centers
\((X_{\ell j},b_{\ell j})\) with
\(R_\ell\leq|X_{\ell j}|\leq2R_\ell\), and use a covering whose
doubles have bounded overlap.  Write \(Q_{\ell j}\) for the fixed reference product of balls.  Set
\begin{equation*}
 h_{\ell j}(\xi,\eta)
 =R_\ell^{n/2}u(X_{\ell j}+\xi,b_{\ell j}+R_\ell\eta).
\end{equation*}
On each doubled reference box, \(|X_{\ell j}+\xi|\asymp R_\ell\).
The change of variables gives
\begin{align*}
 \int|h_{\ell j}|^2&=\int|u|^2,\\
 \int|\nabla_\xi h_{\ell j}|^2&=\int|\nabla_xu|^2,\\
 \int|\nabla_\eta h_{\ell j}|^2
 &=R_\ell^2\int|\nabla_yu|^2
 \leq C\int|x|^2|\nabla_yu|^2.
\end{align*}
Integrals are over the corresponding boxes.  Bounded overlap in
physical space, including the sum in \(\ell\), gives
\begin{equation*}
 \sum_{\ell,j}\int_{2Q_{ell j}}
 (|\nabla h_{\ell j}|^2+|h_{\ell j}|^2)
 \leq C\mathcal E_a(u).
\end{equation*}
Moreover, \(|X_{\ell j}|\geq32>(2\varepsilon)^{-1}\).
The doubled physical boxes are allowed by \eqref{eq:theta-e}, so
\begin{equation*}
 \sup_{\ell,j}\int_{2Q_{ell j}}|h_{\ell j}|^2
 \leq C\Theta_E(u).
\end{equation*}
The transformed mass coefficient satisfies
\begin{equation*}
 R_\ell^{-nq/2+n}
 (1+|X_{\ell j}+\xi|^2)^{n/(Q-2)}\asymp1,
 \qquad -\frac{nq}{2}+n+\frac{2n}{Q-2}=0.
\end{equation*}
Apply \eqref{eq:euclidean-lions} once to the full family
\((h_{\ell j})\), with \(p=q\).  Since the annuli with \(\ell\geq5\)
cover \(|x|>\varepsilon^{-1}=40\), we get
\begin{equation}
 \int_{|x|>\varepsilon^{-1}}\omega^{q/Q}|u|^q
 \leq C\Theta_E(u)^\gamma.
 \label{eq:escape-vanishing-bound}
\end{equation}

The gauges are bounded on \(\mathcal E_a(u)\leq1\).  Taking the
smallest exponent in \eqref{eq:singular-vanishing-bound},
\eqref{eq:finite-vanishing-bound}, and \eqref{eq:escape-vanishing-bound}
and summing gives
\begin{equation}
 \mathcal M(u)\leq C
 \left(\max\{\Theta_F(u),\Theta_G(u),\Theta_E(u)\}\right)^\gamma.
 \label{eq:normalized-refined-global}
\end{equation}
The interpolation exponents satisfy \(0<\gamma<q/2\).
For \(u\ne0\), put \(v=u/\mathcal E_a(u)^{1/2}\).
Then
\begin{equation*}
 \mathcal M(v)=\mathcal E_a(u)^{-q/2}\mathcal M(u),\qquad
 \Theta_\tau(v)=\mathcal E_a(u)^{-1}\Theta_\tau(u).
\end{equation*}
Substitution into \eqref{eq:normalized-refined-global} proves
\eqref{eq:refined-global}, hence \eqref{eq:mass-vanishing}.
\end{proof}

A lower bound for one gauge gives a nonzero profile with a uniform energy lower bound.

\begin{lemma}
\label{lem:inverse-vanishing}
Let $\mathcal E_a(u_k)\leq L$ and
$\limsup_k\mathcal M(u_k)>0$.  After a subsequence, one of the following
holds.
\begin{enumerate}[label=\textup{\(\roman*\)}]
\item There are $b_k$ such that
$T_{-b_k}u_k\rightharpoonup U\ne0$ in $H_a$.
\item There are $r_k\to0$ and $c_k$ such that
$D_{r_k,c_k}^{-1}u_k\rightharpoonup\Phi\ne0$ in
$D_G^{1,2}$.
\item There are $|X_k|\to\infty$ and $b_k$ such that
$A_{X_k,b_k}^{-1}u_k\rightharpoonup\Psi\ne0$ in $H^1_{\rm loc}$.
\end{enumerate}
In \textup{(iii)}, $\Psi\in H^1(\mathbb R^N)$.  If
$\limsup_k\mathcal M(u_k)\geq\eta>0$, then the selected profile satisfies
\begin{equation}
 \mathcal E_*(p)\geq c(L,\eta)>0.
 \label{eq:quantitative-profile-energy}
\end{equation}
If the maximum of the three gauges has limit superior at least
\(\theta>0\), the selected profile satisfies
\begin{equation}
 \mathcal E_*(p)\geq c\theta,
 \label{eq:gauge-profile-energy}
\end{equation}
where \(c>0\) depends only on \(m,n,a\).
\end{lemma}

\begin{proof}
Suppose first that the maximum gauge has limit superior at least
\(\theta>0\).  On a subsequence, one fixed gauge is at least
\(\tau=\theta/2\).  For each \(k\), choose a box carrying at least
\(\tau/2\) of that gauge.

Assume first that \(\Theta_G(u_k)\geq\tau\).  Choose
\(0<r_k<16\varepsilon\) and \(c_k\) such that
\begin{equation*}
 r_k^{-2}\int_{\mathcal Q((0,c_k),r_k)}|u_k|^2\geq\frac\tau2.
\end{equation*}
Set \(v_k=D_{r_k,c_k}^{-1}u_k\).  Scaling gives
\begin{equation*}
 \mathcal E_G(v_k)=\mathcal E_G(u_k)\leq L,\qquad
 \int_{\mathcal Q((0,0),1)}|v_k|^2\geq\frac\tau2.
\end{equation*}
For every fixed bounded set \(K\), Sobolev--H\"older gives
\begin{equation*}
 \int_K|v_k|^2\leq |K|^{2/Q}\|v_k\|_q^2\leq C_KL.
\end{equation*}
Weak compactness and Lemma \ref{lem:uniform-local-rellich} on
increasing boxes give
\begin{equation*}
 v_k\rightharpoonup\Phi\quad\hbox{in }D_G^{1,2},\qquad
 v_k\longrightarrow\Phi\quad\hbox{in }L^2_{\rm loc}.
\end{equation*}
Consequently,
\begin{equation*}
 \frac\tau2\leq\int_{\mathcal Q((0,0),1)}|\Phi|^2
 \leq |\mathcal Q((0,0),1)|^{2/Q}\|\Phi\|_q^2
 \leq C\mathcal E_G(\Phi).
\end{equation*}
Pass to a further subsequence with \(r_k\to r_\infty\).
If \(r_\infty=0\), this proves \textup{(ii)}.
If \(r_\infty>0\), the sequence \(T_{-c_k}u_k\) is bounded in
\(H_a\).  Take its weak limit \(U\).  On compact sets,
\begin{equation*}
 T_{-c_k}u_k(x,y)
 =r_k^{-(Q-2)/2}v_k(x/r_k,y/r_k^2)
 \longrightarrow D_{r_\infty,0}\Phi
 \quad\hbox{in }L^2_{\rm loc}.
\end{equation*}
Here the convergence follows from local \(L^2\) convergence of
\(v_k\), followed by continuity of fixed dilations in \(L^2\).
Hence \(U=D_{r_\infty,0}\Phi\), and
\begin{equation*}
 \frac\tau2\leq r_\infty^{-2}
 \int_{\mathcal Q((0,0),r_\infty)}|U|^2
 \leq r_\infty^{-2}|\mathcal Q((0,0),r_\infty)|^{2/Q}\|U\|_q^2
 \leq C\mathcal E_a(U).
\end{equation*}
The constant is independent of \(r_\infty\), since
\(|\mathcal Q((0,0),r_\infty)|^{2/Q}=Cr_\infty^2\).

Assume next that \(\Theta_E(u_k)\geq\tau\), and choose
\(X_k,b_k\) as above.  If \((X_k)\) is bounded, the selected boxes,
after translation by \(T_{-b_k}\), lie in one fixed bounded set
\(K\).  Weak compactness in \(H_a\) and local Rellich compactness
give
\begin{equation*}
 T_{-b_k}u_k\rightharpoonup U\quad\hbox{in }H_a,\qquad
 T_{-b_k}u_k\longrightarrow U\quad\hbox{in }L^2(K).
\end{equation*}
Since \(\int_K|T_{-b_k}u_k|^2\geq\tau/2\),
\begin{equation*}
 \|U\|_2^2\geq\frac\tau2,\qquad
 \mathcal E_a(U)\geq\frac{a\tau}{2}.
\end{equation*}

Otherwise pass to \(R_k=|X_k|\to\infty\), and set
\(v_k=A_{X_k,b_k}^{-1}u_k\).  The selected box becomes
\begin{equation*}
 K_0=\{(\xi,\eta):|\xi|<2\varepsilon^{-1},\
 |\eta|<2\varepsilon^{-1}\},\qquad
 \int_{K_0}|v_k|^2\geq\frac\tau2.
\end{equation*}
For a fixed ball \(B_R\), put
\(c_k(\xi)=|X_k+\xi|^2/R_k^2\).  Then
\(c_k\to1\) uniformly on \(B_R\).  Formula
\eqref{eq:moving-escape-norm} gives, for all large \(k\),
\begin{equation*}
 \int_{B_R}(|\nabla v_k|^2+a|v_k|^2)\leq2L.
\end{equation*}
A diagonal subsequence satisfies
\begin{equation*}
 v_k\rightharpoonup\Psi\quad\hbox{in }H^1(B_R),\qquad
 v_k\longrightarrow\Psi\quad\hbox{in }L^2(B_R)
\end{equation*}
for every fixed \(R\).  Take \(K_0\subset B_R\).
Then \(\int_{K_0}|\Psi|^2\geq\tau/2\).
For any fixed \(R\), uniform convergence of \(c_k\) and lower
semicontinuity give
\begin{equation*}
 \int_{B_R}(|\nabla\Psi|^2+a|\Psi|^2)
 \leq\liminf_k\int_{B_R}
 (|\nabla_\xi v_k|^2+c_k|\nabla_\eta v_k|^2+a|v_k|^2)
 \leq L.
\end{equation*}
Let \(R\to\infty\).  We obtain \(\Psi\in H^1(\mathbb R^N)\) and
\begin{equation*}
 \frac{a\tau}{2}\leq\mathcal E_E^a(\Psi)\leq L.
\end{equation*}

Finally, assume that \(\Theta_F(u_k)\geq\tau\).
Choose \(b_k\) from \eqref{eq:theta-f}.  The translated sequence
is bounded in \(H_a\).  On the fixed defining box, take a strong
\(L^2\) subsequence and the corresponding weak \(H_a\) limit \(U\).
Then
\begin{equation*}
 \|U\|_2^2\geq\frac\tau2,\qquad
 \mathcal E_a(U)\geq\frac{a\tau}{2}.
\end{equation*}
Thus every case gives \(\mathcal E_*(p)\geq c\theta\), with \(c\)
depending only on \(m,n,a\).  This proves
\eqref{eq:gauge-profile-energy}.

If \(\limsup_k\mathcal M(u_k)\geq\eta>0\),
\eqref{eq:refined-global} implies
\begin{equation*}
 \limsup_k\max\{\Theta_F(u_k),\Theta_G(u_k),\Theta_E(u_k)\}
 \geq\left(\frac\eta{CL^{q/2-\gamma}}\right)^{1/\gamma}.
\end{equation*}
Apply the gauge argument with this positive lower bound.
It gives \eqref{eq:quantitative-profile-energy}, for example with
\begin{equation*}
 c(L,\eta)=c\left(\frac\eta{CL^{q/2-\gamma}}\right)^{1/\gamma}.
\end{equation*}

\end{proof}

\begin{proof}[Proof of Theorem \ref{thm:profile-intro}]
Use Theorem \ref{thm:three-family-criterion} with ambient space
\(H_a\) and model spaces \(H_a\), \(D_G^{1,2}\), and
\(H^1(\mathbb R^N)\).  Their energies are
\(\mathcal E_a\), \(\mathcal E_G\), and \(\mathcal E_E^a\).
Their masses are \(\mathcal M\), \(\|\cdot\|_q^q\), and
\(\|\cdot\|_q^q\).  The constant finite frame is the identity frame.
All three spaces are separable.  Compact smooth functions are dense;
for the singular space use Appendix \ref{sec:cutoffs}.

For the local bounds, finite pullbacks preserve \(\mathcal E_a\).  Singular pullbacks satisfy
\begin{equation*}
 \mathcal E_G(D_{r_k,c_k}^{-1}u_k)\leq\mathcal E_a(u_k),\qquad
 \|D_{r_k,c_k}^{-1}u_k\|_{L^2(K)}^2
 \leq C_K\mathcal E_a(u_k)
\end{equation*}
for every fixed bounded \(K\).  The second bound is Sobolev--H\"older.  Escape pullbacks have bounded
local \(H^1\) norms by the coefficient limit in
\eqref{eq:moving-escape-norm}.  In the singular and escape
cases, lower semicontinuity on fixed boxes, followed by increasing
boxes, gives
\begin{equation*}
 \mathcal E_*(p)\leq\liminf_k\mathcal E_a(u_k).
\end{equation*}
Lemma \ref{lem:uniform-local-rellich} gives local strong \(L^2\)
convergence along a subsequence.  A diagonal over fixed boxes gives
almost-everywhere convergence.

For \textup{(i)}, fix a compact model cutoff \(h_L\).
The frame pullback of its realization is exactly \(h_L\).
Lemma \ref{lem:one-frame-limits} gives convergence of its energy and
mass.  The singular and escape mass densities are
\begin{equation*}
 (1+r_k^2|\xi|^2)^{n/(Q-2)},\qquad
 R_k^{-nq/2+n}(1+|X_k+\xi|^2)^{n/(Q-2)}.
\end{equation*}
Both converge uniformly to one on the fixed support, with bounded
reciprocals.  If the pullbacks of \(u_k\) converge
weakly to \(p\), the fixed-cutoff pairing in Lemma
\ref{lem:one-profile-subtraction} gives
\begin{equation*}
 \langle u_k,P_{k,L}h\rangle_{H_a}
 \longrightarrow\langle p,h_L\rangle_{H_*}.
\end{equation*}
For a singular frame, the extra pairing is
\(ar_k^2\int(D_{r_k,c_k}^{-1}u_k)h_L\), which tends to zero by the
local \(L^2\) bound.  For an escape frame, the gradient coefficient
converges uniformly on the fixed support.  Finite-frame pairings are
preserved exactly.  Local almost-everywhere convergence and these
mass-density bounds give the local weighted Brezis--Lieb formula;
see Lemma \ref{lem:weighted-bl}.  

For \textup{(ii)}, take a pair of frames.  If they have the same type,
pass to a subsequence on which their relative parameters either
converge or satisfy the corresponding orthogonality condition.
When the parameters converge, the relative formulas in Section 4
have limits of the following forms:
\begin{equation*}
 f(x,y+d),\qquad
 \alpha^{(Q-2)/2}\phi(\alpha\xi,\beta+\alpha^2\eta),\qquad
 \psi(\xi+d,\eta+e),\qquad \alpha>0.
\end{equation*}
These maps preserve model energy and mass.  For finite and singular
frames, compact smooth tests converge in the model norm on a fixed
compact set; \(\alpha\) stays bounded above and away from zero.
Changing variables proves \eqref{eq:criterion-equivalent-covariance}.

For equivalent escape frames, bounded relative \(x\)-displacement
implies \(R_k/\widetilde R_k\to1\).  Their exact relative map is
\begin{equation*}
 A_{X_k,b_k}^{-1}A_{\widetilde X_k,\widetilde b_k}\psi(\xi,\eta)
 =\left(\frac{R_k}{\widetilde R_k}\right)^{n/2}
 \psi\left(\xi+X_k-\widetilde X_k,
 \frac{b_k-\widetilde b_k}{\widetilde R_k}
 +\frac{R_k}{\widetilde R_k}\eta\right).
\end{equation*}
The affine maps converge with their derivatives to a translation.
Changing variables produces strongly convergent tests on a fixed
compact set.  The local \(H^1\) bounds give covariance.
The limiting maps extend by density to unitary maps of the model
spaces.  Composition of the exact relative maps gives the
cocycle rule in the limit.

If the parameters diverge, Lemma \ref{lem:relative-frames} gives
vanishing relative pullbacks in both directions.  Different types
are orthogonal by
\eqref{eq:relative-gf-formula}--\eqref{eq:relative-eg-formula}.
Lemma \ref{lem:pairwise-decoupling} gives energy and mass decoupling
for these pairs.  First fix its cutoff radii, then choose synchronized radii after
the finitely many tests.

For \textup{(iii)}, fix \(p_L\) and first let \(k\to\infty\).
In a finite frame, apply Lemma \ref{lem:weighted-bl} to the pullbacks
and to the pullbacks minus \(p_L\).  Their weak limits are \(p\)
and \(p-p_L\).  Subtract the two mass identities.  For singular
and escape frames, use the fixed-cutoff identity in Lemma
\ref{lem:one-profile-subtraction}.  In all three cases,
\begin{equation*}
 \mathcal M(u_k)-\mathcal M(u_k-P_{k,L}p)
 \longrightarrow\mathcal M_*(p)-\mathcal M_*(p-p_L).
\end{equation*}
This is \eqref{eq:criterion-fixed-core-mass}.  Now let \(L\to\infty\)
to remove the model cutoff errors.
The block selection in Theorem \ref{thm:three-family-criterion} fixes
one admissible realization for every extraction level.

For \textup{(iv)}, the estimates at the start of this section give
\begin{equation*}
 0\leq\Theta_\tau(u)\leq C\mathcal E_a(u),\qquad
 |\Theta_\tau(u)^{1/2}-\Theta_\tau(v)^{1/2}|
 \leq C\mathcal E_a(u-v)^{1/2}.
\end{equation*}
Each gauge has degree two.  Proposition \ref{prop:refined} supplies
\eqref{eq:criterion-refined}.  For \textup{(v)}, apply the first part
of Lemma \ref{lem:inverse-vanishing} to a positive gauge limit superior.
The extracted profile satisfies
\(\mathcal E_*(p)\geq c\theta\), with \(c\) independent of the
sequence's energy bound.

Theorem \ref{thm:three-family-criterion} gives
\eqref{eq:energy-splitting-intro} and
\eqref{eq:mass-splitting-intro}.  Its profile bound and remainder
estimate give \eqref{eq:profile-summability-intro} and
\eqref{eq:grushin-remainder-rate-intro}.
\end{proof}

\section{Sharp constants and the threshold alternative}

\subsection{Upper bounds}

\begin{lemma}
\label{lem:recovery}
For every $a>0$,
\begin{equation*}
 S_a\leq S_G,\qquad S_a\leq S_E(a).
\end{equation*}
\end{lemma}

\begin{proof}
Let \(0\ne\phi\in C_c^\infty\).  Put \(s=n/(Q-2)\).
The change of variables \(x=r\xi\), \(y=r^2\eta\) gives
\begin{align*}
 \mathcal E_a(D_{r,0}\phi)
 &=\mathcal E_G(\phi)+ar^2\|\phi\|_2^2,\\
 \mathcal M(D_{r,0}\phi)
 &=\int(1+r^2|\xi|^2)^s|\phi(\xi,\eta)|^q\,d\xi\,d\eta
 \longrightarrow\|\phi\|_q^q.
\end{align*}
The mass coefficient converges uniformly on the support of \(\phi\).
By Lemma \ref{lem:one-frame-limits},
\begin{equation*}
 S_a\leq\lim_{r\downarrow0}
 \frac{\mathcal E_a(D_{r,0}\phi)}
 {\mathcal M(D_{r,0}\phi)^{2/q}}
 =\frac{\mathcal E_G(\phi)}{\|\phi\|_q^2}.
\end{equation*}
Density in \(D_G^{1,2}\) gives \(S_a\leq S_G\).

Let \(0\ne\psi\in C_c^\infty\), and fix \(e\in\mathbb S^{m-1}\).
In the coordinates \(x=Re+\xi\), \(y=R\eta\),
\begin{align*}
 \mathcal E_a(A_{Re,0}\psi)
 &=\int\left(|\nabla_\xi\psi|^2
 +|e+\xi/R|^2|\nabla_\eta\psi|^2+a|\psi|^2\right),\\
 \mathcal M(A_{Re,0}\psi)
 &=\int\left(R^{-2}+|e+\xi/R|^2\right)^s|\psi|^q.
\end{align*}
Uniform convergence of both coefficients on the support of \(\psi\) gives
\begin{equation*}
 S_a\leq\lim_{R\to\infty}
 \frac{\mathcal E_a(A_{Re,0}\psi)}
 {\mathcal M(A_{Re,0}\psi)^{2/q}}
 =\frac{\mathcal E_E^a(\psi)}{\|\psi\|_q^2}.
\end{equation*}
Density in \(H^1(\mathbb R^N)\) gives \(S_a\leq S_E(a)\).
\end{proof}

\subsection{Binding tests and explicit nonemptiness}

Let \(s=n/(Q-2)\) and define
\begin{align*}
 \kappa_E^{\mathrm{test}}={}&\frac{mn}{N(N+2)}
 \int|\zeta|^2|\nabla Z_1|^2
 -\frac nQ S_E(1)\frac{m+2(s-1)}N\int|\zeta|^2Z_1^q,\\
 a_E^{\mathrm{test}}={}&\frac{Q\kappa_E^{\mathrm{test}}}{nS_E(1)}.
\end{align*}

\begin{proposition}
\label{prop:recovery-curvature}
Let \(e\in\mathbb S^{m-1}\).  Then
\begin{equation*}
 \frac{\mathcal E_a(A_{Re,0}Z_a)}
 {\mathcal M(A_{Re,0}Z_a)^{2/q}}
 =S_E(a)+R^{-2}\kappa_E(a)+o(R^{-2}),
\end{equation*}
where
\begin{align}
 \kappa_E(a)={}&\frac{mn}{N(N+2)}
 \int |\zeta|^2|\nabla Z_a|^2
 -\frac nQ S_E(a)\left(1+
 \frac{m+2(s-1)}N\int |\zeta|^2Z_a^q\right)
 \notag\\
 ={}&a^{-N/Q}\left(\kappa_E^{\mathrm{test}}-\frac nQ S_E(1)a\right).
 \label{eq:escape-curvature-coefficient}
\end{align}
Consequently, \(\kappa_E(a)<0\) implies \(S_a<S_E(a)\).

If \(Q>4\), let \(\Phi\) be a positive normalized extremal for
\(S_G\).  Its required moments are finite by Proposition
\ref{prop:automatic-singular-moments}.  Then
\begin{equation}
 \frac{\mathcal E_a(D_{r,0}\Phi)}
 {\mathcal M(D_{r,0}\Phi)^{2/q}}
 =S_G+r^2\kappa_G(a)+o(r^2),
 \label{eq:singular-curvature-expansion}
\end{equation}
where
\begin{equation*}
 \kappa_G(a)=a\int\Phi^2-\frac nQ S_G\int |x|^2|\Phi|^q.
\end{equation*}
Consequently, \(\kappa_G(a)<0\) implies \(S_a<S_G\).
\end{proposition}

\begin{proof}
Write \(\zeta=(\xi,\eta)\in\mathbb R^m\times\mathbb R^n\).
Exponential decay gives all polynomial moments of \(Z_a\).
The energy is
\begin{align*}
 \mathcal E_a(A_{Re,0}Z_a)
 &=\int\left(|\nabla_\xi Z_a|^2
 +\left(1+\frac{2e\cdot\xi}{R}+\frac{|\xi|^2}{R^2}\right)
 |\nabla_\eta Z_a|^2+aZ_a^2\right)\\
 &=S_E(a)+\frac1{R^2}\int|\xi|^2|\nabla_\eta Z_a|^2.
\end{align*}
The term of order \(R^{-1}\) is odd in \(\xi\).
For distinct coordinates \(i,j\), spherical integration gives
\begin{equation*}
 \frac1{|\mathbb S^{N-1}|}
 \int_{\mathbb S^{N-1}}\theta_i^2\theta_j^2\,d\theta
 =\frac1{N(N+2)}.
\end{equation*}
Radial symmetry and summation over the \(mn\) coordinate pairs give
\begin{equation*}
 \int|\xi|^2|\nabla_\eta Z_a|^2
 =\frac{mn}{N(N+2)}\int|\zeta|^2|\nabla Z_a|^2.
\end{equation*}

The mass is
\begin{equation*}
 \mathcal M(A_{Re,0}Z_a)
 =\int\left(1+\frac{2e\cdot\xi}{R}
 +\frac{1+|\xi|^2}{R^2}\right)^s Z_a^q.
\end{equation*}
On \(|\zeta|\leq R^{1/2}\), put
\begin{equation*}
 h_R=\frac{2e\cdot\xi}{R}+\frac{1+|\xi|^2}{R^2}.
\end{equation*}
Then \(|h_R|\leq CR^{-1/2}\), and
\begin{equation*}
 (1+h_R)^s=1+sh_R+\frac{s(s-1)}2h_R^2+O(|h_R|^3).
\end{equation*}
The remainder satisfies
\begin{equation*}
 \int_{|\zeta|\leq R^{1/2}}|h_R|^3Z_a^q
 \leq\frac C{R^3}\int(1+|\zeta|^6)Z_a^q=O(R^{-3}).
\end{equation*}
On the complement, exponential decay gives
\begin{equation*}
 \int_{|\zeta|>R^{1/2}}
 \left[1+(1+|\zeta|^2)^s+|\zeta|^4\right]Z_a^q
 =o(R^{-2}).
\end{equation*}
The same polynomial bounds the mass coefficient and Taylor terms
uniformly for \(R\geq1\).  Odd terms vanish, so
\begin{align*}
 \mathcal M(A_{Re,0}Z_a)
 &=1+\frac s{R^2}\int(1+|\xi|^2)Z_a^q
 +\frac{2s(s-1)}{R^2}\int(e\cdot\xi)^2Z_a^q+o(R^{-2})\\
 &=1+\frac s{R^2}\left(1+
 \frac{m+2(s-1)}N\int|\zeta|^2Z_a^q\right)+o(R^{-2}).
\end{align*}
Here we used
\begin{equation*}
 \int|\xi|^2Z_a^q=\frac mN\int|\zeta|^2Z_a^q,
 \qquad
 \int(e\cdot\xi)^2Z_a^q=\frac1N\int|\zeta|^2Z_a^q.
\end{equation*}
Apply \((1+t)^{-2/q}=1-(2/q)t+o(t)\).
Since \(2s/q=n/Q\), the coefficient of \(R^{-2}\) is the first
line of \eqref{eq:escape-curvature-coefficient}.
The substitution \(\vartheta=\sqrt a\,\zeta\) gives
\begin{align*}
 \int|\zeta|^2|\nabla Z_a|^2
 &=a^{-N/Q}\int|\vartheta|^2|\nabla Z_1|^2,\\
 \int|\zeta|^2Z_a^q
 &=a^{-1}\int|\vartheta|^2Z_1^q.
\end{align*}
Using \(S_E(a)=a^{n/Q}S_E(1)\) and \(n/Q-1=-N/Q\)
gives the second line.
For \(\kappa_E(a)<0\), the quotient lies below \(S_E(a)\)
when \(R\) is large.

For the singular test, scaling gives
\begin{equation*}
 \mathcal E_a(D_{r,0}\Phi)=S_G+ar^2\int\Phi^2.
\end{equation*}
Here \(0<s\leq1\).  For \(0<r\leq1\),
\begin{equation*}
 0\leq\frac{(1+r^2|x|^2)^s-1}{r^2}\leq s|x|^2.
\end{equation*}
Proposition \ref{prop:automatic-singular-moments} makes the right side
integrable against \(|\Phi|^q\,dz\).  Dominated convergence gives
\begin{equation*}
 \mathcal M(D_{r,0}\Phi)
 =1+sr^2\int|x|^2|\Phi|^q+o(r^2).
\end{equation*}
Hence
\begin{align*}
 \frac{\mathcal E_a(D_{r,0}\Phi)}
 {\mathcal M(D_{r,0}\Phi)^{2/q}}
 &=\left(S_G+ar^2\int\Phi^2\right)
 \left(1-\frac{2s}{q}r^2\int|x|^2|\Phi|^q+o(r^2)\right)\\
 &=S_G+r^2\left(a\int\Phi^2
 -\frac nQ S_G\int|x|^2|\Phi|^q\right)+o(r^2).
\end{align*}
This gives \eqref{eq:singular-curvature-expansion}.
For \(\kappa_G(a)<0\), the quotient lies below \(S_G\) when \(r>0\)
is small.
\end{proof}

\begin{proposition}
\label{thm:binding-curvature-intro}
The binding set is nonempty if
\begin{equation}
 a_0>a_E^{\mathrm{test}}.
 \label{eq:escape-curvature-crossing}
\end{equation}
If \(Q>4\), choose a positive normalized extremal \(\Phi\) for \(S_G\)
and define
\begin{equation*}
 a_G^{\mathrm{test}}=\frac{nS_G}{Q}
 \frac{\int |x|^2\Phi^q}{\int\Phi^2}.
\end{equation*}
Then the binding set is also nonempty if
\begin{equation}
 a_0<a_G^{\mathrm{test}}.
 \label{eq:singular-curvature-crossing}
\end{equation}
\end{proposition}

\begin{proof}
Under \eqref{eq:escape-curvature-crossing}, choose
\begin{equation*}
 \max\{0,a_E^{\mathrm{test}}\}<a<a_0.
\end{equation*}
Then
\begin{equation*}
 \kappa_E(a)=\frac nQ S_E(1)a^{-N/Q}
 (a_E^{\mathrm{test}}-a)<0,
 \qquad S_E(a)<S_E(a_0)=S_G.
\end{equation*}
Proposition \ref{prop:recovery-curvature} gives
\(S_a<S_E(a)<S_G\).

Under \eqref{eq:singular-curvature-crossing}, choose
\(a_0<a<a_G^{\mathrm{test}}\).  Then
\begin{equation*}
 \kappa_G(a)=\left(\int\Phi^2\right)
 (a-a_G^{\mathrm{test}})<0,
 \qquad S_G=S_E(a_0)<S_E(a).
\end{equation*}
It gives \(S_a<S_G<S_E(a)\).  Thus \(a\in\mathcal B\) in either case.
\end{proof}

\begin{proof}[Proof of Theorem \ref{thm:explicit-binding}]
Let \(n=1\), \(m\geq6\) be even, \(N=m+1\), and \(Q=m+2\).
Proposition \ref{prop:heisenberg-comparison} gives
\begin{equation*}
 S_G=m^2b_m^{2/Q},\qquad
 \int|x|^2\Phi_m^q=1,\qquad
 \int\Phi_m^2=\frac{4m}{m-2}b_m^{2/Q}.
\end{equation*}
Thus
\begin{align*}
 \kappa_G(a)
 &=\frac{4ma}{m-2}b_m^{2/Q}
 -\frac{m^2}{m+2}b_m^{2/Q}\\
 &=\frac{4m}{m-2}b_m^{2/Q}
 \left(a-\frac{m(m-2)}{4(m+2)}\right).
\end{align*}
Thus \(S_a<S_G\) for \(0<a<\overline a_m\).

Let
\begin{equation*}
 S_N^{\mathrm{Sob}}=N(N-2)\pi
 \left(\frac{\Gamma(N/2)}{\Gamma(N)}\right)^{2/N},
 \qquad \theta=\frac NQ.
\end{equation*}
For \(v\in H^1(\mathbb R^N)\), put
\(A=\int|\nabla v|^2\) and \(B=\int v^2\).
Since
\begin{equation*}
 \frac1q=\theta\frac{N-2}{2N}+\frac{1-\theta}{2},
\end{equation*}
interpolation and the sharp Sobolev inequality \cite{Talenti1976} give
\begin{equation*}
 \|v\|_q^2\leq
 \|v\|_{2N/(N-2)}^{2\theta}\|v\|_2^{2(1-\theta)}
 \leq(S_N^{\mathrm{Sob}})^{-\theta}A^\theta B^{1-\theta}.
\end{equation*}
For \(A+B>0\), maximizing \(t^\theta(1-t)^{1-\theta}\)
on \([0,1]\) gives \(t=\theta\) and
\begin{equation*}
 A^\theta B^{1-\theta}
 \leq\theta^\theta(1-\theta)^{1-\theta}(A+B),
\end{equation*}
and
\begin{equation*}
 S_E(1)\geq
 \frac{(S_N^{\mathrm{Sob}})^{N/Q}}
 {(N/Q)^{N/Q}(1/Q)^{1/Q}}.
\end{equation*}
Since \(Q=N+1\),
\begin{align*}
 a_0=\left(\frac{S_G}{S_E(1)}\right)^Q
 &\leq\frac{m^{2Q}b_m^2(N/Q)^N(1/Q)}
 {(S_N^{\mathrm{Sob}})^N}\\
 &=\frac{m^{2Q}b_m^2}
 {Q^Q(N-2)^N\pi^N}
 \left(\frac{\Gamma(N)}{\Gamma(N/2)}\right)^2.
\end{align*}
Use
\begin{equation*}
 b_m=\frac{\pi^{Q/2}}{2^{m+1}\Gamma(Q/2)},\qquad
 \frac{\Gamma(N)}{\Gamma(N/2)}
 =\frac{2^{N-1}\Gamma((N+1)/2)}{\sqrt\pi}.
\end{equation*}
Canceling the gamma factors and using \(N-2=m-1\) gives
\begin{equation*}
 a_0\leq
 \frac{m^{2m+4}}
 {4(m+2)^{m+2}(m-1)^{m+1}}=\underline a_m.
\end{equation*}
Moreover,
\begin{equation*}
 \frac{\underline a_m}{\overline a_m}
 =\frac m{m-2}
 \left(1+\frac{m-2}{m^2}\right)^{-(m+1)}.
\end{equation*}
Bernoulli's inequality gives
\begin{align*}
 \left(1+\frac{m-2}{m^2}\right)^{m+1}
 &\geq1+\frac{(m+1)(m-2)}{m^2},\\
 \frac{(m+1)(m-2)}{m^2}-\frac2{m-2}
 &=\frac{m^2(m-5)+4}{m^2(m-2)}>0.
\end{align*}
Thus \(\underline a_m<\overline a_m\).
For \(\underline a_m<a<\overline a_m\),
\begin{equation*}
 a>a_0,\qquad S_a<S_G<S_E(a).
\end{equation*}
This proves \eqref{eq:explicit-binding-interval}.
For \(m=6\),
\begin{equation*}
 \underline a_6=\frac{43046721}{80000000}
 <\frac23<\frac34=\overline a_6.
\end{equation*}
Theorem \ref{thm:threshold-intro} gives attainment and compactness.
\end{proof}

\begin{remark}
The interval \((\underline a_m,\overline a_m)\) is contained in
\(\mathcal B\).  The actual endpoints of \(\mathcal B\) remain
undetermined.  Proposition \ref{prop:recovery-curvature} expands the
specified recovery quotients; it does not give a relaxed expansion
of \(S_a\).
\end{remark}

\subsection{Exclusion of dichotomy}

\begin{lemma}
\label{lem:concavity}
Let $0<\alpha<1$ and $t_j\geq0$.  Then
\begin{equation*}
 \sum_jt_j^\alpha\geq\left(\sum_jt_j\right)^\alpha.
\end{equation*}
If $\sum_jt_j=1$, equality holds only when exactly one $t_j$ is nonzero.
\end{lemma}

\begin{proof}
For \(s,t>0\), put \(\lambda=s/(s+t)\).  Since
\(0<\lambda<1\) and \(0<\alpha<1\),
\begin{equation*}
 s^\alpha+t^\alpha
 =(s+t)^\alpha\left(\lambda^\alpha+(1-\lambda)^\alpha\right)
 >(s+t)^\alpha.
\end{equation*}
The nonstrict inequality also holds if \(s=0\) or \(t=0\).
Iteration gives
\begin{equation*}
 \sum_{j=1}^Jt_j^\alpha
 \geq\left(\sum_{j=1}^Jt_j\right)^\alpha.
\end{equation*}
Monotone convergence as \(J\to\infty\) gives the inequality.
If \(\sum_jt_j=1\) and two terms are positive, combine them first.
Their strict gap remains after adding the other terms.  Equality
therefore requires one term equal to one and all others equal to zero.
\end{proof}

\begin{lemma}
\label{lem:profilewise-model-bounds}
Let \(p\ne0\) be a profile in Theorem \ref{thm:profile-intro}.  Then
\begin{equation*}
 \mathcal E_{\tau}(p)
 \geq S_\tau\mathcal M_\tau(p)^{2/q},
 \qquad
 S_F=S_a,\qquad S_E=S_E(a).
\end{equation*}
Moreover, \(\mathcal M_\tau(p)>0\).
\end{lemma}

\begin{proof}
By definition,
\begin{equation*}
 \frac{\mathcal E_\tau(p)}{\mathcal M_\tau(p)^{2/q}}
 \geq S_\tau.
\end{equation*}
For \(G\) and \(E\), \(\mathcal M_\tau(p)=\int|p|^q>0\).
For \(F\), positivity of \(\omega^{q/Q}\) gives
\(\mathcal M_F(p)>0\).  Multiply by \(\mathcal M_\tau(p)^{2/q}\).
\end{proof}

\begin{proof}[Proof of Theorem \ref{thm:threshold-intro}]
Proposition \ref{prop:embedding} gives
\(\mathcal M(u)^{2/q}\leq C_a\mathcal E_a(u)\).
Hence \(S_a\geq C_a^{-1}>0\).  Lemma \ref{lem:recovery} gives the
upper bounds.

Let \(\mathcal M(u_k)=1\) and \(\mathcal E_a(u_k)\to S_a\).
Apply Theorem \ref{thm:profile-intro}, and put
\begin{equation*}
 t_j=\mathcal M_{\tau_j}(p^j).
\end{equation*}
For each \(J\), mass splitting gives
\begin{equation*}
 1=\sum_{j=1}^Jt_j+\mathcal M(w_k^J)+o_k(1).
\end{equation*}
Let \(k\to\infty\), then \(J\to\infty\).
The remainder mass vanishes, so \(\sum_jt_j=1\).
Energy splitting and Lemma \ref{lem:profilewise-model-bounds} give
\begin{align*}
 S_a
 &\geq\sum_j\mathcal E_{\tau_j}(p^j)\\
 &\geq\sum_{\tau_j=F}S_at_j^{2/q}
 +\sum_{\tau_j=G}S_Gt_j^{2/q}
 +\sum_{\tau_j=E}S_E(a)t_j^{2/q}\\
 &\geq S_a\sum_jt_j^{2/q}\geq S_a.
\end{align*}
By Lemma \ref{lem:concavity}, equality leaves one nonzero profile,
denoted by \(p^1\).  Then
\begin{equation*}
 \mathcal M_{\tau_1}(p^1)=1,\qquad
 \mathcal E_{\tau_1}(p^1)=S_a=S_{\tau_1},\qquad
 \mathcal E_a(w_k^1)
 =\mathcal E_a(u_k)-\mathcal E_{\tau_1}(p^1)+o(1)\longrightarrow0.
\end{equation*}

Under \eqref{eq:strict-binding}, \(\tau_1=G\) would give
\(S_a=S_G\), and \(\tau_1=E\) would give \(S_a=S_E(a)\).
Both contradict binding.  Hence \(\tau_1=F\), and the profile is a normalized minimizer:
\begin{equation*}
 \|u_k-T_{b_k}p^1\|_{H_a}^2=\mathcal E_a(w_k^1)\longrightarrow0.
\end{equation*}

If \(\tau_1=G\), the profile is a normalized extremal for \(S_G\).
Its realization satisfies
\begin{equation*}
 \mathcal E_a\bigl(u_k-D_{r_k,c_k}(\chi_{L_k}p^1)\bigr)
 =\mathcal E_a(w_k^1)\longrightarrow0.
\end{equation*}
If \(\tau_1=E\), Section 3 gives
\begin{equation*}
 p^1(\xi,\eta)=\psi(\xi,\eta)
 =\sigma Z_a(\xi-\xi_0,\eta-\eta_0),\qquad \sigma\in\{-1,1\}.
\end{equation*}
Set
\begin{equation*}
 X'_k=X_k+\xi_0,\qquad b'_k=b_k+R_k\eta_0,\qquad
 R'_k=|X'_k|,\qquad t_k=R_k/R'_k\longrightarrow1.
\end{equation*}
The difference \(A_{X_k,b_k}\psi-\sigma A_{X'_k,b'_k}Z_a\)
is represented in the original escape coordinates by
\begin{equation*}
 h_k(\xi,\eta)=\sigma\bigl[Z_a(\xi-\xi_0,\eta-\eta_0)
 -t_k^{n/2}Z_a(\xi-\xi_0,t_k(\eta-\eta_0))\bigr].
\end{equation*}
For large \(k\), \(1/2\leq t_k\leq2\).
Exponential decay bounds the squared integrand in
\(\|h_k\|_{H^1}^2+\||\xi|\nabla_\eta h_k\|_2^2\)
by \(C(1+|\zeta|^4)e^{-c|\zeta|}\), with fixed \(C,c>0\).
Pointwise convergence to zero and domination give
\begin{equation*}
 \|h_k\|_{H^1}+\||\xi|\nabla_\eta h_k\|_2\longrightarrow0.
\end{equation*}
Since \(|X_k+\xi|^2/R_k^2\leq2+2|\xi|^2/R_k^2\),
\begin{align*}
 \|A_{X_k,b_k}\psi-\sigma A_{X'_k,b'_k}Z_a\|_{H_a}^2
 &=\|h_k\|_{X_k,a}^2\\
 &\leq C_a\|h_k\|_{H^1}^2
 +2R_k^{-2}\||\xi|\nabla_\eta h_k\|_2^2\longrightarrow0.
\end{align*}
Exponential decay and the cutoff construction replace the canonical
realization by the uncut extremal with error \(o_{H_a}(1)\).
After recentering,
\begin{equation*}
 \left\|A_{X'_k,b'_k}^{-1}u_k-\sigma Z_a\right\|_{X'_k,a}^2
 =\|u_k-\sigma A_{X'_k,b'_k}Z_a\|_{H_a}^2\longrightarrow0.
\end{equation*}
The convergence is global in the norm \eqref{eq:moving-escape-norm}.
If \(S_G=S_E(a)=S_a\), two positive tangent masses would make
\(\sum_jt_j^{2/q}>1\), so they cannot coexist.
\end{proof}

\subsection{Dependence on the coefficient a}

\begin{proposition}
\label{prop:phase-properties}
The function $a\mapsto S_a$ is nondecreasing, concave, and locally Lipschitz on
$(0,\infty)$.  Moreover,
\begin{equation}
 \frac{S_a}{a^{n/Q}}=
 \inf_{v\in H_1\setminus\{0\}}
 \frac{\mathcal E_1(v)}
 {\left(\displaystyle\int(a+|x|^2)^{n/(Q-2)}|v|^q\right)^{2/q}},
 \label{eq:parameter-rescaling}
\end{equation}
so \(a\mapsto S_a/S_E(a)\) is nonincreasing on \((0,\infty)\).
If \(S_{a_1}\) is attained and \(a_2>a_1\), then
\(S_{a_2}/S_E(a_2)<S_{a_1}/S_E(a_1)\).
The set
\begin{equation*}
\{a>0:S_a<S_G\}
\end{equation*}
is an interval of the form \((0,a_*)\), with \(a_*\in[a_0,\infty]\),
and it contains \((0,a_0)\).  At every
$a$ where $S_a$ is attained by a normalized minimizer $U_a$, one has
\begin{equation*}
 (S_a)'_+\leq\int U_a^2\leq(S_a)'_-.
\end{equation*}
If \eqref{eq:strict-binding} holds, let $\mathcal F_a$ be the set of
normalized minimizers.  Then
\begin{align}
 (S_a)'_+&=\inf_{U\in\mathcal F_a}\int U^2,
 \label{eq:right-subgradient}\\
 (S_a)'_-&=\sup_{U\in\mathcal F_a}\int U^2.
 \label{eq:left-subgradient}
\end{align}
\end{proposition}

\begin{proof}
For \(\mathcal M(u)=1\),
\begin{equation*}
 \mathcal E_a(u)=\mathcal E_G(u)+a\|u\|_2^2.
\end{equation*}
Hence \(S_b\geq S_a\) for \(b>a\).
For \(0<t<1\),
\begin{align*}
 S_{ta+(1-t)b}
 &=\inf_{\mathcal M(u)=1}
 \left[t\mathcal E_a(u)+(1-t)\mathcal E_b(u)\right]\\
 &\geq tS_a+(1-t)S_b.
\end{align*}
Concavity follows; finiteness and positivity follow from Lemma
\ref{lem:recovery} and Proposition \ref{prop:embedding}.
For \(\alpha\leq a<b\leq\beta\), concavity gives
\begin{equation*}
 0\leq\frac{S_b-S_a}{b-a}
 \leq\frac{S_\alpha-S_{\alpha/2}}{\alpha/2}
 \leq\frac{2S_G}{\alpha}.
\end{equation*}
For the bijection \(u=D_{a^{-1/2},0}v\) of the common underlying space,
scaling gives
\begin{equation*}
 \mathcal E_a(u)=\mathcal E_1(v),\qquad
 \mathcal M(u)=a^{-n/(Q-2)}
 \int(a+|x|^2)^{n/(Q-2)}|v|^q.
\end{equation*}
Since \(2n/[q(Q-2)]=n/Q\), taking infima gives
\eqref{eq:parameter-rescaling}.
For \(v\ne0\) and \(a_2>a_1\),
\begin{equation*}
 \int(a_2+|x|^2)^{n/(Q-2)}|v|^q
 >\int(a_1+|x|^2)^{n/(Q-2)}|v|^q.
\end{equation*}
Each test quotient strictly decreases.  Taking infima gives \(S_{a_2}/a_2^{n/Q}\leq S_{a_1}/a_1^{n/Q}\).
If the infimum at \(a_1\) is attained, its rescaled minimizer
as a test at \(a_2\) gives
\begin{equation*}
 \frac{S_{a_2}}{a_2^{n/Q}}
 <\frac{S_{a_1}}{a_1^{n/Q}}.
\end{equation*}
Division by \(S_E(1)\) gives the assertions for \(S_a/S_E(a)\).
Monotonicity and continuity show that \(\{a:S_a<S_G\}\)
is an initial open interval.  It contains \((0,a_0)\), since
\(S_a\leq S_E(a)<S_G\) there.

Let \(U_a\) be a normalized minimizer, and take \(0<h<a\).
Testing at the neighboring parameters gives
\begin{equation*}
 S_{a+h}\leq S_a+h\|U_a\|_2^2,\qquad
 S_{a-h}\leq S_a-h\|U_a\|_2^2.
\end{equation*}
Thus
\begin{equation*}
 \frac{S_{a+h}-S_a}{h}\leq\|U_a\|_2^2
 \leq\frac{S_a-S_{a-h}}h.
\end{equation*}
Let \(h\downarrow0\) to obtain the derivative bounds.

Under \eqref{eq:strict-binding}, choose a normalized minimizer
\(U_{a+h}\) for small \(|h|\).  Since \(a+h\geq a/2\),
\begin{equation*}
 \|U_{a+h}\|_2^2\leq\frac{S_{a+h}}{a+h}\leq\frac{2S_G}{a}.
\end{equation*}
For every \(U\in\mathcal F_a\), the two test comparisons give
\begin{equation*}
 h\|U_{a+h}\|_2^2\leq S_{a+h}-S_a\leq h\|U\|_2^2.
\end{equation*}
Also,
\begin{align*}
 \mathcal E_a(U_{a+h})&=S_{a+h}-h\|U_{a+h}\|_2^2,\\
 0\leq\mathcal E_a(U_{a+h})-S_a
 &\leq|S_{a+h}-S_a|+\frac{2S_G}{a}|h|\longrightarrow0.
\end{align*}
For any \(h_j\to0\), Theorem \ref{thm:threshold-intro} gives a strong
\(H_a\) limit in \(\mathcal F_a\) after translations and a subsequence.
For \(h_j\downarrow0\), denote it by \(U_+\).  Dividing by
\(h_j>0\) gives
\begin{equation*}
 \|U_{a+h_j}\|_2^2
 \leq\frac{S_{a+h_j}-S_a}{h_j}\leq\|U\|_2^2.
\end{equation*}
Hence
\begin{equation*}
 \inf_{V\in\mathcal F_a}\|V\|_2^2
 \leq\|U_+\|_2^2\leq(S_a)'_+
 \leq\inf_{U\in\mathcal F_a}\|U\|_2^2.
\end{equation*}
This is \eqref{eq:right-subgradient}.
For \(h_j\uparrow0\), write the strong limit as \(U_-\).
Division by \(h_j<0\) reverses the comparison:
\begin{equation*}
 \|U\|_2^2\leq\frac{S_{a+h_j}-S_a}{h_j}
 \leq\|U_{a+h_j}\|_2^2.
\end{equation*}
Consequently,
\begin{equation*}
 \sup_{U\in\mathcal F_a}\|U\|_2^2
 \leq(S_a)'_-\leq\|U_-\|_2^2
 \leq\sup_{V\in\mathcal F_a}\|V\|_2^2.
\end{equation*}
This is \eqref{eq:left-subgradient}.
\end{proof}

\begin{remark}
Translation in $x$ changes both the energy coefficient and the mass
density.  Thus rearrangement alone does not classify minimizers.
\end{remark}

\section{Parameter slopes and endpoint attainment}

Let \(U\geq0\) be a normalized minimizer:
\(\mathcal M(U)=1\), \(\mathcal E_a(U)=S_a\).  Testing the Lagrange
multiplier identity by \(U\) gives
\begin{align*}
 &\int_{\mathbb R^N}\left(
 \nabla_xU\cdot\nabla_x\varphi+
|x|^2\nabla_yU\cdot\nabla_y\varphi+aU\varphi\right)\,dz\\
 &\hspace{30mm}
 =S_a\int_{\mathbb R^N}(1+|x|^2)^{n/(Q-2)}
 U^{q-1}\varphi\,dz,
 \qquad \varphi\in C_c^\infty(\mathbb R^N).
\end{align*}
Equivalently,
\begin{equation*}
 -\Delta_xU-|x|^2\Delta_yU+aU
 =S_a(1+|x|^2)^{n/(Q-2)}U^{q-1}
 \quad\hbox{in }\mathbb R^N
\end{equation*}
in the weak sense.  The Pohozaev identities follow by differentiating
the quotient along the following scalings.

For $U\in H_a$ put
\begin{equation*}
 A_x=\int|\nabla_xU|^2,\qquad
 A_y=\int|x|^2|\nabla_yU|^2,\qquad
 B=a\int U^2.
\end{equation*}
For $\lambda,\mu>0$ define
\begin{equation*}
 U_\lambda(x,y)=U(\lambda x,y),\qquad
 V_\mu(x,y)=U(x,\mu y).
\end{equation*}

\begin{lemma}
\label{lem:scaling-derivatives}
Let $U\in H_a$ and $\mathcal M(U)<\infty$.  The two paths above are
continuous in $H_a$ for $\lambda,\mu>0$.  At $\lambda=\mu=1$,
\begin{align*}
 \frac{d}{d\lambda}\mathcal E_a(U_\lambda)
 &=(2-m)A_x-(m+2)A_y-mB,\\
 \frac{d}{d\lambda}\mathcal M(U_\lambda)
 &=-m\mathcal M(U)-\int x\cdot\nabla_x(\omega^{q/Q})|U|^q,
 \\
 \frac{d}{d\mu}\mathcal E_a(V_\mu)
 &=-nA_x+(2-n)A_y-nB,\\
 \frac{d}{d\mu}\mathcal M(V_\mu)
 &=-n\mathcal M(U).
\end{align*}
\end{lemma}

\begin{proof}
The changes \(\xi=\lambda x\) and \(\eta=\mu y\) give
\begin{align}
 \mathcal E_a(U_\lambda)
 &=\lambda^{2-m}A_x+\lambda^{-m-2}A_y+\lambda^{-m}B,\notag\\
 \mathcal M(U_\lambda)
 &=\lambda^{-m}\int\omega(x/\lambda)^{q/Q}|U(x,y)|^q\,dx\,dy,
 \label{eq:x-mass-exact}\\
 \mathcal E_a(V_\mu)
 &=\mu^{-n}A_x+\mu^{2-n}A_y+\mu^{-n}B,\notag\\
 \mathcal M(V_\mu)&=\mu^{-n}\mathcal M(U).\notag
\end{align}
For \(\lambda,\mu\in[1/2,2]\), these formulas give
\begin{equation*}
 \|(U-U_j)_\lambda\|_{H_a}
 +\|(U-U_j)(x,\mu y)\|_{H_a}
 \leq C\|U-U_j\|_{H_a}.
\end{equation*}
Choose \(U_j\in C_c^\infty\) converging to \(U\) in \(H_a\).
The scaled paths of \(U_j\) are continuous.  For fixed
\(\lambda_0\in[1/2,2]\),
\begin{equation*}
 \|U_\lambda-U_{\lambda_0}\|_{H_a}
 \leq2C\|U-U_j\|_{H_a}
 +\|(U_j)_\lambda-(U_j)_{\lambda_0}\|_{H_a}.
\end{equation*}
Let \(\lambda\to\lambda_0\), then \(j\to\infty\).
Apply the same estimate to \(V_\mu\) and to any compact parameter
interval.
Proposition \ref{prop:embedding} also gives continuity in the weighted
\(L^q\) norm.

Put \(w=\omega^{q/Q}\).  Direct differentiation gives
\begin{equation}
 x\cdot\nabla_xw(x)=\frac{2n}{Q-2}
 \frac{|x|^2}{1+|x|^2}w(x),
 \qquad |x\cdot\nabla_xw(x)|\leq\frac{2n}{Q-2}w(x).
 \label{eq:weight-derivative-bound}
\end{equation}
For \(1/2\leq\lambda\leq2\),
\begin{align*}
 \partial_\lambda\bigl[\lambda^{-m}w(x/\lambda)\bigr]
 &=-\lambda^{-m-1}\left[mw(x/\lambda)
 +(x/\lambda)\cdot\nabla_xw(x/\lambda)\right],\\
 \left|\partial_\lambda\bigl[\lambda^{-m}w(x/\lambda)\bigr]\right|
 &\leq Cw(x).
\end{align*}
Since \(w|U|^q\in L^1\), dominated convergence in \eqref{eq:x-mass-exact} gives
\begin{equation*}
 \left.\frac{d}{d\lambda}\mathcal M(U_\lambda)\right|_{\lambda=1}
 =-m\mathcal M(U)-\int x\cdot\nabla_xw(x)|U|^q.
\end{equation*}
The energy powers and the mass factor \(\mu^{-n}\) give the remaining
derivatives.
\end{proof}

Put
\begin{equation*}
 C=S_a\int\omega^{q/Q}|U|^q.
\end{equation*}

Stationarity under the two scalings gives the following identities.

\begin{proposition}
\label{prop:pohozaev}
Every normalized minimizer satisfies
\begin{align}
 A_x+A_y+B&=C,\label{eq:pohozaev-energy}\\
 \frac{m-2}{2}A_x+\frac{m+2}{2}A_y+\frac m2B
 &=\frac{S_a}{q}\int
 \left(m\omega^{q/Q}+x\cdot\nabla_x(\omega^{q/Q})\right)|U|^q,
 \label{eq:pohozaev-x}\\
 \frac n2A_x+\frac{n-2}{2}A_y+\frac n2B
 &=\frac{nS_a}{q}\int\omega^{q/Q}|U|^q.
 \label{eq:pohozaev-y}
\end{align}
Consequently,
\begin{equation}
 \frac{Q-2}{2}(A_x+A_y)+\frac Q2B
 =\frac{S_a}{q}\int
 \left(Q\omega^{q/Q}+x\cdot\nabla_x(\omega^{q/Q})\right)|U|^q.
 \label{eq:pohozaev-anisotropic}
\end{equation}
\end{proposition}

\begin{proof}
Since \(\mathcal M(U)=1\) and \(\mathcal E_a(U)=S_a\),
\begin{equation*}
 A_x+A_y+B=S_a=C.
\end{equation*}
The quotient along \(U_\lambda\) has a minimum at \(\lambda=1\).
Lemma \ref{lem:scaling-derivatives} gives
\begin{align*}
 0&=\left.\frac{d}{d\lambda}
 \left[\mathcal E_a(U_\lambda)\mathcal M(U_\lambda)^{-2/q}\right]
 \right|_{\lambda=1}\\
 &=(2-m)A_x-(m+2)A_y-mB
 +\frac{2S_a}{q}\left(m+
 \int x\cdot\nabla_x(\omega^{q/Q})|U|^q\right).
\end{align*}
Division by \(-2\) gives \eqref{eq:pohozaev-x}.
The quotient along \(V_\mu\) likewise gives
\begin{equation*}
 0=-nA_x+(2-n)A_y-nB+\frac{2nS_a}{q}.
\end{equation*}
Division by \(-2\) gives \eqref{eq:pohozaev-y}.
Add \eqref{eq:pohozaev-x} to twice \eqref{eq:pohozaev-y}.
The coefficients satisfy
\begin{equation*}
 \frac{m-2}{2}+n=\frac{Q-2}{2},\qquad
 \frac{m+2}{2}+n-2=\frac{Q-2}{2},\qquad
 \frac m2+n=\frac Q2.
\end{equation*}
Since \((m+2n)\omega^{q/Q}=Q\omega^{q/Q}\), this gives
\eqref{eq:pohozaev-anisotropic}.
\end{proof}

\begin{corollary}
\label{cor:mass-identity}
Every normalized minimizer satisfies
\begin{equation}
 B=\frac{S_a}{q}
 \int x\cdot\nabla_x(\omega^{q/Q})|U|^q.
 \label{eq:mass-identity}
\end{equation}
Moreover,
\begin{equation}
 0<a\int U^2<\frac nQ S_a.
 \label{eq:mass-slope-bound}
\end{equation}
\end{corollary}

\begin{proof}
Subtract \((Q-2)/2\) times \eqref{eq:pohozaev-energy}
from \eqref{eq:pohozaev-anisotropic}.  Since \(Q/q=(Q-2)/2\),
the gradient terms cancel, as do the mass terms on the right, leaving
\begin{equation*}
 B=\frac{S_a}{q}\int x\cdot\nabla_x(\omega^{q/Q})|U|^q.
\end{equation*}
By \eqref{eq:weight-derivative-bound},
\begin{equation*}
 B=\frac nQ S_a
 \int\frac{|x|^2}{1+|x|^2}\omega^{q/Q}|U|^q.
\end{equation*}
The measure \(\omega^{q/Q}|U|^q\,dz\) has total mass one.
The set \(\{x=0\}\) has measure zero, and
\(0<|x|^2/(1+|x|^2)<1\) elsewhere.  The last integral lies
strictly between zero and one, giving \eqref{eq:mass-slope-bound}.
\end{proof}

\begin{corollary}
\label{cor:slope-escape}
For every \(a>0\),
\begin{equation*}
 0\leq a(S_a)'_+(a)\leq a(S_a)'_-(a)\leq\frac nQ S_a.
\end{equation*}
At every point of the binding set,
\begin{equation}
 0<a(S_a)'_+(a)\leq a(S_a)'_-(a)<\frac nQ S_a.
 \label{eq:sa-log-slope}
\end{equation}
Consequently $a\mapsto S_a/S_E(a)$ is strictly decreasing on every
connected component of the binding set.
\end{corollary}

\begin{proof}
Monotonicity and concavity of \(S_a\) give
\(0\leq(S_a)'_+\leq(S_a)'_-\).
For \(0<h<a\), monotonicity of \(S_a/a^{n/Q}\) gives
\begin{equation*}
 S_{a-h}\geq S_a\left(1-\frac ha\right)^{n/Q}.
\end{equation*}
Therefore
\begin{equation*}
 (S_a)'_-(a)
 =\lim_{h\downarrow0}\frac{S_a-S_{a-h}}h
 \leq S_a\lim_{h\downarrow0}
 \frac{1-(1-h/a)^{n/Q}}h=\frac{nS_a}{Qa}.
\end{equation*}
Suppose \(a\in\mathcal B\).  Theorem \ref{thm:threshold-intro}
makes \(\mathcal F_a/\mathbb R_y^n\) compact.
The continuous, translation-invariant map \(U\mapsto\|U\|_2^2\)
attains its extrema on this quotient.  Applying
\eqref{eq:mass-slope-bound}, derived from \eqref{eq:mass-identity},
to both extremizers gives
\begin{equation*}
 0<a\inf_{U\in\mathcal F_a}\|U\|_2^2
 \leq a\sup_{U\in\mathcal F_a}\|U\|_2^2<\frac nQ S_a.
\end{equation*}
Equations \eqref{eq:right-subgradient} and \eqref{eq:left-subgradient}
now give \eqref{eq:sa-log-slope}.

Concavity gives \(S_a\in AC_{\rm loc}(0,\infty)\).
Since \(S_E(a)=a^{n/Q}S_E(1)>0\),
\(\log(S_a/S_E(a))\) is locally absolutely continuous.
At almost every point of a binding component,
\begin{equation*}
 \frac{d}{da}\log\frac{S_a}{S_E(a)}
 =\frac{S_a'}{S_a}-\frac{n}{Qa}<0.
\end{equation*}
For \(a_1<a_2\) in that component,
\begin{equation*}
 \log\frac{S_{a_2}}{S_E(a_2)}-
 \log\frac{S_{a_1}}{S_E(a_1)}
 =\int_{a_1}^{a_2}
 \left(\frac{S_a'}{S_a}-\frac{n}{Qa}\right)\,da<0.
\end{equation*}
\end{proof}

\begin{proof}[Proof of Theorem \ref{thm:connected-binding-intro}]
Proposition \ref{prop:phase-properties} gives continuity and
\begin{equation*}
 S_{a_1}\leq S_{a_2},\qquad
 \frac{S_{a_2}}{S_E(a_2)}\leq\frac{S_{a_1}}{S_E(a_1)}
 \quad\hbox{for }a_1<a_2.
\end{equation*}
Thus \(\{a:S_a<S_G\}\) is an initial open interval, and
\(\{a:S_a<S_E(a)\}\) is a terminal open interval.
Their intersection \(\mathcal B\) is empty or one open interval.
Suppose \(a\in\mathcal B\).
If \(a<a_0\), then
\begin{equation*}
 S_{a_0}/S_E(a_0)\leq S_a/S_E(a)<1.
\end{equation*}
If \(a>a_0\), then \(S_{a_0}\leq S_a<S_G\).
Since \(S_E(a_0)=S_G\), both cases give \(a_0\in\mathcal B\).
Openness gives
\(a_-<a_0<a_+\).
For \(0<a\leq a_-\), one has \(a<a_0\), and hence
\(S_a\leq S_E(a)<S_G\).  This parameter lies outside \(\mathcal B\), so
\(S_a=S_E(a)\).  For \(a\geq a_+\), the inequalities
\(S_a\leq S_G<S_E(a)\) similarly give \(S_a=S_G\).
These assertions concern only finite positive parameters.

If \(a_->0\), put \(a=a_-+h\), where \(h\downarrow0\).
The one-sided expansions give
\begin{align*}
 S_{a_-+h}&=S_E(a_-)+h(S_a)'_+(a_-)+o(h),\\
 S_E(a_-+h)&=S_E(a_-)+hS_E'(a_-)+o(h).
\end{align*}
Their difference is \(\ell_-h+o(h)\), which proves
\eqref{eq:left-endpoint-asymptotic}.
If \(a_+<\infty\), put \(a=a_+-h\).  Then
\begin{equation*}
 S_{a_+-h}=S_G-h(S_a)'_-(a_+)+o(h).
\end{equation*}
This proves \eqref{eq:right-endpoint-asymptotic}.
Both gaps are positive in \(\mathcal B\), so
\(\ell_-,\ell_+\geq0\).

Suppose \(a_->0\) and a normalized \(U\) attains \(S_{a_-}\).
Proposition \ref{prop:phase-properties} and
\eqref{eq:mass-slope-bound} give
\begin{equation*}
 (S_a)'_+(a_-)\leq\|U\|_2^2
 <\frac{nS_{a_-}}{Qa_-}=S_E'(a_-).
\end{equation*}
Hence \(\ell_->0\).
If \(a_+<\infty\) is attained, choose a normalized minimizer \(U_+\)
at \(a_+\).  The other subgradient bound gives
\begin{equation*}
 \ell_+=(S_a)'_-(a_+)\geq\|U_+\|_2^2>0.
\end{equation*}

Let \(a_j\in\mathcal B\) tend to an endpoint
\(a_*\in(0,\infty)\), and choose normalized minimizers \(U_j\).
For large \(j\),
\begin{equation*}
 a_j\geq a_*/2,\qquad
 \|U_j\|_2^2\leq\frac{2S_G}{a_*},\qquad
 \mathcal E_{a_*}(U_j)
 =S_{a_j}+(a_*-a_j)\|U_j\|_2^2\longrightarrow S_{a_*}.
\end{equation*}
Thus \((U_j)\) is minimizing at \(a_*\).
Suppose that endpoint is not attained.
Theorem \ref{thm:threshold-intro} excludes a finite profile.
At \(a_*=a_-\), the strict inequality \(S_E(a_-)<S_G\)
leaves only an escape profile.
Since the escape map is an \(L^2\) isometry,
\begin{align*}
 \left|\|U_j\|_2-\|Z_{a_-}\|_2\right|
 &\leq\|U_j-\sigma A_{X_j,b_j}Z_{a_-}\|_2\\
 &\leq a_-^{-1/2}
 \|U_j-\sigma A_{X_j,b_j}Z_{a_-}\|_{H_{a_-}}\longrightarrow0.
\end{align*}
Consequently,
\begin{equation*}
 \|U_j\|_2^2\longrightarrow\|Z_{a_-}\|_2^2
 =\frac{nS_E(a_-)}{Qa_-}=S_E'(a_-).
\end{equation*}
At \(a_*=a_+\), the inequality \(S_G<S_E(a_+)\) excludes an
escape profile.  Write the singular realization as
\(D_{r_j,c_j}(\chi_{L_j}\Phi)\).  Its cutoff choice gives
\begin{equation*}
 \|D_{r_j,c_j}(\chi_{L_j}\Phi)\|_2
 =r_j\|\chi_{L_j}\Phi\|_2\longrightarrow0.
\end{equation*}
The strong remainder gives
\begin{equation*}
 \|U_j\|_2\leq r_j\|\chi_{L_j}\Phi\|_2
 +a_+^{-1/2}\|U_j-D_{r_j,c_j}(\chi_{L_j}\Phi)\|_{H_{a_+}}
 \longrightarrow0.
\end{equation*}

At every \(a_j\),
\begin{equation*}
 (S_a)'_+(a_j)\leq\|U_j\|_2^2\leq(S_a)'_-(a_j).
\end{equation*}
For fixed small \(h>0\) and \(a_-<a_j<a_-+h\), concavity gives
\begin{equation*}
 \frac{S_{a_-+h}-S_{a_j}}{a_-+h-a_j}
 \leq(S_a)'_+(a_j)\leq(S_a)'_-(a_j)
 \leq\frac{S_{a_j}-S_{a_-}}{a_j-a_-}.
\end{equation*}
Let \(j\to\infty\), then \(h\downarrow0\).
Both bounds tend to \((S_a)'_+(a_-)\).
At the right endpoint, use \(a_+-h<a_j<a_+\) and
\begin{equation*}
 \frac{S_{a_+}-S_{a_j}}{a_+-a_j}
 \leq(S_a)'_+(a_j)\leq(S_a)'_-(a_j)
 \leq\frac{S_{a_j}-S_{a_+-h}}{a_j-a_++h}.
\end{equation*}
The same order of limits gives \((S_a)'_-(a_+)\).
Thus a nonattained left endpoint has
\((S_a)'_+(a_-)=S_E'(a_-)\), and a nonattained right endpoint has
\((S_a)'_-(a_+)=0\).  Hence \(\ell_-=0\) or \(\ell_+=0\), respectively.
The attained-endpoint inequalities give both equivalences.
\end{proof}

\begin{corollary}
\label{cor:exterior-nonattainment}
Assume \(\mathcal B\ne\varnothing\).
If \(0<a<a_-\), then \(S_a=S_E(a)\) is not attained, and every
normalized minimizing sequence has, after a subsequence, one escape
profile of full mass.
If \(a_+<a<\infty\), then \(S_a=S_G\) is not attained, and every
normalized minimizing sequence has, after a subsequence, one singular
profile of full mass.  In both cases the realized profile has a strong
remainder in \(H_a\).
\end{corollary}

\begin{proof}
Let \(0<a<a_-\).  Theorem \ref{thm:connected-binding-intro}
gives \(S_b=S_E(b)\) for all \(b\) near \(a\).
If a normalized minimizer \(U\) existed, Proposition
\ref{prop:phase-properties} would give
\begin{equation*}
 S_E'(a)=(S_a)'_+(a)\leq\|U\|_2^2
 \leq(S_a)'_-(a)=S_E'(a).
\end{equation*}
Thus \(a\|U\|_2^2=nS_a/Q\), contrary to
\eqref{eq:mass-slope-bound}.
For \(a_+<a<\infty\), the same theorem gives \(S_b=S_G\)
near \(a\).  The subgradient bounds would give
\begin{equation*}
 0\leq\|U\|_2^2\leq(S_a)'_-(a)=0,
\end{equation*}
a contradiction.

For a normalized minimizing sequence, Theorem
\ref{thm:threshold-intro} excludes a finite profile by nonattainment.
On the left, \(a<a_0\) and \(S_a=S_E(a)<S_G\), which excludes a
singular profile.  On the right, \(a>a_0\) and
\(S_a=S_G<S_E(a)\), which excludes an escape profile.
The remaining profile has mass one and energy \(S_a\), so
\begin{equation*}
 \mathcal E_a(w_k^1)
 =\mathcal E_a(u_k)-S_a+o(1)\longrightarrow0.
\end{equation*}
\end{proof}

\begin{corollary}
\label{cor:interior-endpoint-compactness}
Assume \(\mathcal B\ne\varnothing\), let
\(a_*\in\{a_-,a_+\}\cap(0,\infty)\), and suppose that \(S_{a_*}\)
is attained.  Let \(a_j\in\mathcal B\) tend to \(a_*\), and let \(U_j\)
be any normalized minimizer for \(S_{a_j}\).
After a subsequence, there are \(b_j\in\mathbb R^n\) and a normalized
minimizer \(U_*\) for \(S_{a_*}\) such that
\begin{equation*}
 \|T_{-b_j}U_j-U_*\|_{H_{a_*}}\longrightarrow0.
\end{equation*}
The limit satisfies
\begin{equation*}
 \int U_*^2=
 \begin{cases}
 (S_a)'_+(a_-)=S_E'(a_-)-\ell_-,&a_*=a_-,\\
 (S_a)'_-(a_+)=\ell_+,&a_*=a_+.
 \end{cases}
\end{equation*}
\end{corollary}

\begin{proof}
Since \(a_j\to a_*>0\), for large \(j\),
\begin{equation*}
 \|U_j\|_2^2\leq\frac{2S_G}{a_*},\qquad
 \mathcal E_{a_*}(U_j)
 =S_{a_j}+(a_*-a_j)\|U_j\|_2^2\longrightarrow S_{a_*}.
\end{equation*}
Thus \((U_j)\) is minimizing at \(a_*\).  The subgradient bounds give
\begin{equation*}
 (S_a)'_+(a_j)\leq\|U_j\|_2^2\leq(S_a)'_-(a_j).
\end{equation*}
The concavity bounds in Theorem \ref{thm:connected-binding-intro} give
\begin{equation*}
 \|U_j\|_2^2\longrightarrow
 \begin{cases}
 S_E'(a_-)-\ell_-,&a_*=a_-,\\
 \ell_+,&a_*=a_+.
 \end{cases}
\end{equation*}

At \(a_*=a_-\), \(S_E(a_-)<S_G\).
Theorem \ref{thm:threshold-intro} permits only a finite or an escape
profile.  An escape profile would give
\begin{equation*}
 \|U_j\|_2^2\longrightarrow\|Z_{a_-}\|_2^2=S_E'(a_-).
\end{equation*}
This contradicts \(\ell_->0\) from endpoint attainment.
At \(a_*=a_+\), \(S_G<S_E(a_+)\), and only a finite or a singular
profile is possible.  A singular realization has \(L^2\) norm tending
to zero.  Its strong remainder would give \(\|U_j\|_2^2\to0\),
contrary to \(\ell_+>0\).

The finite profile \(U_*\) has mass one and energy \(S_{a_*}\).
Energy splitting gives
\begin{equation*}
 \|T_{-b_j}U_j-U_*\|_{H_{a_*}}^2
 =\mathcal E_{a_*}(U_j)-\mathcal E_{a_*}(U_*)+o(1)
 \longrightarrow0.
\end{equation*}
In particular, \(\|U_j\|_2^2\to\|U_*\|_2^2\).
The endpoint mass limits give the identity.
\end{proof}

\section{Global compactness of Palais--Smale sequences}

\subsection{Critical profiles of Palais--Smale sequences}

For the two tangent spaces put
\begin{align*}
 I_G(\phi)=\frac12\mathcal E_G(\phi)-\frac1q\int|\phi|^q,\qquad
 I_E(\psi)=\frac12\mathcal E_E^a(\psi)-\frac1q\int|\psi|^q.
\end{align*}

\begin{lemma}
\label{lem:ps-profile-criticality}
Let \((u_k)\) be bounded in \(H_a\), with
\(I_a'(u_k)\to0\) in \(H_a^*\).  Every profile in Theorem
\ref{thm:profile-intro} is a critical point of its model functional
\(I_F=I_a\), \(I_G\), or \(I_E\).  Every nonzero profile of type
\(\tau\) satisfies
\begin{equation}
 \mathcal E_\tau(p)=\mathcal M_\tau(p)
 \geq S_\tau^{q/(q-2)},\qquad
 I_\tau(p)\geq c_\tau.
 \label{eq:critical-profile-quantization}
\end{equation}
\end{lemma}

\begin{proof}
Fix a profile \(p=p^j\), and let \(T_k^j\) be its pullback.
Lemma \ref{lem:relative-frames} and the cutoff choice give
\begin{equation*}
 T_k^ju_k=T_k^jw_k^{j-1}
 +\sum_{i<j}T_k^jR_k^ip^i\rightharpoonup p.
\end{equation*}
This convergence holds in \(H_a\) for type \(F\), in \(D_G^{1,2}\)
for type \(G\), and in \(H^1(B_R)\) for each fixed \(R\) for type \(E\).
For the identity-frame profile, use the weak limit of \(u_k\).

Let \(v\in C_c^\infty\) be a model test function.  Its realization
\(R_kv\) satisfies \(\sup_k\|R_kv\|_{H_a}<\infty\).  Hence
\begin{equation*}
 |\langle I_a'(u_k),R_kv\rangle|
 \leq\|I_a'(u_k)\|_{H_a^*}\|R_kv\|_{H_a}\longrightarrow0.
\end{equation*}
The pullback formulas in Section 4 give
\begin{equation*}
 \langle u_k,R_kv\rangle_{H_a}
 \longrightarrow\langle p,v\rangle_{H_\tau}
\end{equation*}
For type \(G\), the zeroth-order pairing has the factor \(r_k^2\)
and tends to zero on the fixed support of \(v\).
For type \(E\), the coefficients converge uniformly on that support.

Put \(f_k=T_k^ju_k\), and let \(K\) be a compact model set containing
the support of \(v\).  Lemma \ref{lem:uniform-local-rellich} gives,
after a subsequence, \(f_k\to p\) almost everywhere on \(K\).
Write \(\rho_k\) for the pulled-back mass density.  Then
\begin{equation*}
 \rho_k\longrightarrow\rho_\tau\quad\hbox{uniformly on }K,
 \qquad 0<c_K\leq\rho_k\leq C_K,
 \qquad\sup_k\int_K|f_k|^q<\infty.
\end{equation*}
Here \(\rho_F=\omega^{q/Q}\) and \(\rho_G=\rho_E=1\).
Thus
\begin{equation*}
 |f_k|^{q-2}f_k\rightharpoonup|p|^{q-2}p
 \quad\hbox{in }L^{q/(q-1)}(K).
\end{equation*}
Moreover,
\begin{align*}
 &\left|\int_K\rho_k|f_k|^{q-2}f_kv
       -\int_K\rho_\tau|p|^{q-2}pv\right|\\
 &\quad\leq\|\rho_k-\rho_\tau\|_{L^\infty(K)}
       \|f_k\|_{L^q(K)}^{q-1}\|v\|_{L^q(K)}
 +\left|\int_K\rho_\tau
       (|f_k|^{q-2}f_k-|p|^{q-2}p)v\right|
 \longrightarrow0.
\end{align*}
It follows that \(\langle I_\tau'(p),v\rangle=0\).

For \(v\in H_\tau\), choose \(v_\ell\in C_c^\infty\) with
\(v_\ell\to v\) in \(H_\tau\).  The model Sobolev inequality and
H\"older give
\begin{equation*}
 |\langle I_\tau'(p),v-v_\ell\rangle|
 \leq\|p\|_{H_\tau}\|v-v_\ell\|_{H_\tau}
 +\mathcal M_\tau(p)^{(q-1)/q}
  \mathcal M_\tau(v-v_\ell)^{1/q}\longrightarrow0.
\end{equation*}
Hence \(I_\tau'(p)=0\) on \(H_\tau\).
Testing with \(p\) gives
\begin{equation*}
 \mathcal E_\tau(p)=\mathcal M_\tau(p),\qquad
 \mathcal M_\tau(p)\geq S_\tau\mathcal M_\tau(p)^{2/q}.
\end{equation*}
Since \(p\ne0\), \eqref{eq:critical-profile-quantization} follows:
\begin{equation*}
 \mathcal E_\tau(p)=\mathcal M_\tau(p)
 \geq S_\tau^{q/(q-2)},\qquad
 I_\tau(p)=\left(\frac12-\frac1q\right)\mathcal E_\tau(p)
 \geq c_\tau.
\end{equation*}
\end{proof}

\begin{proof}[Proof of Theorem \ref{thm:ps-global-intro}]
Put \(\delta_k=\|I_a'(u_k)\|_{H_a^*}\to0\).  Then
\begin{align*}
 \left(\frac12-\frac1q\right)\|u_k\|_{H_a}^2
 &=I_a(u_k)-\frac1q\langle I_a'(u_k),u_k\rangle\\
 &\leq C+\frac{\delta_k}{q}\|u_k\|_{H_a}.
\end{align*}
Since \(q>2\), \((u_k)\) is bounded in \(H_a\) and
\(\langle I_a'(u_k),u_k\rangle=o(1)\).  Hence
\begin{equation}
 \mathcal E_a(u_k)-\mathcal M(u_k)=o(1),\qquad
 \mathcal E_a(u_k)\longrightarrow\frac{c}{1/2-1/q}.
 \label{eq:ps-total-energy}
\end{equation}
In particular, \(c\geq0\).

Apply Theorem \ref{thm:profile-intro}.  By Lemma
\ref{lem:ps-profile-criticality}, every nonzero profile is critical and
has energy at least
\begin{equation*}
 e_*:=\min_{\tau\in\{F,G,E\}}S_\tau^{q/(q-2)}>0.
\end{equation*}
For every finite \(J\leq J_*\), energy splitting gives
\begin{equation*}
 Je_*\leq\sum_{j=1}^J\mathcal E_{\tau_j}(p^j)
 \leq\lim_k\mathcal E_a(u_k)=\frac{c}{1/2-1/q}.
\end{equation*}
Thus \(J_*<\infty\).  Write \(J=J_*\) and \(w_k=w_k^J\).
The extraction continues whenever a limiting gauge is positive, so
the terminal gauges vanish.  Hence \(\mathcal M(w_k)\to0\).
Subtract the energy and mass splittings:
\begin{align*}
 \mathcal E_a(u_k)-\mathcal M(u_k)
 &=\sum_{j=1}^J
   \bigl(\mathcal E_{\tau_j}(p^j)-\mathcal M_{\tau_j}(p^j)\bigr)
   +\mathcal E_a(w_k)-\mathcal M(w_k)+o(1)\\
 &=\mathcal E_a(w_k)-\mathcal M(w_k)+o(1).
\end{align*}
Equation \eqref{eq:ps-total-energy} now gives
\begin{equation*}
 \mathcal E_a(w_k)=\mathcal M(w_k)+o(1)\longrightarrow0.
\end{equation*}
This proves \eqref{eq:global-ps-decomposition} and
\eqref{eq:global-ps-energy}.  Also,
\begin{equation*}
 c=\left(\frac12-\frac1q\right)
   \sum_{j=1}^J\mathcal E_{\tau_j}(p^j)
  =\sum_{j=1}^JI_{\tau_j}(p^j)
  \geq\sum_{j=1}^Jc_{\tau_j},
\end{equation*}
which is \eqref{eq:global-ps-action}.
If \(c=0\), \eqref{eq:ps-total-energy} gives \(u_k\to0\) in \(H_a\).
If no profiles occur, the terminal remainder is \(u_k\), so the same
argument gives \(u_k\to0\) and \(c=0\).
\end{proof}

\begin{proof}[Proof of the compactness assertion in Corollary
\ref{thm:stability-ps-intro}]
By \eqref{eq:global-ps-action}, a singular profile would give
\(c\geq c_G\), an escape profile would give \(c\geq c_E\), and two
finite profiles would give \(c\geq2c_F\).
Condition \eqref{eq:ps-threshold-intro} and \(c>0\) therefore leave
exactly one finite profile \(p\ne0\).  Equation
\eqref{eq:global-ps-decomposition} gives
\begin{equation*}
 u_k=T_{b_k}p+o_{H_a}(1),\qquad
 \|T_{-b_k}u_k-p\|_{H_a}\longrightarrow0.
\end{equation*}
\end{proof}

\begin{proof}[Proof of the sharpness assertion in Corollary
\ref{thm:stability-ps-intro}]
Let \(a\in\mathcal B\), and put \(s=n/(Q-2)\).
Choose normalized extremals \(U\in\mathcal F_a\), \(\Phi\) for
\(S_G\), and \(Z_a\) for \(S_E(a)\).  Set
\begin{equation*}
 W_F=S_a^{1/(q-2)}U,\qquad
 W_G=S_G^{1/(q-2)}\Phi,\qquad
 W_E=S_E(a)^{1/(q-2)}Z_a.
\end{equation*}
For each \(\tau\in\{F,G,E\}\),
\begin{equation*}
 I_\tau'(W_\tau)=0,\qquad
 \mathcal E_\tau(W_\tau)=\mathcal M_\tau(W_\tau)
 =S_\tau^{q/(q-2)},\qquad I_\tau(W_\tau)=c_\tau.
\end{equation*}

Let \(|b_k|\to\infty\), and put
\begin{equation*}
 u_k^F=T_{b_k}W_F+T_{-b_k}W_F.
\end{equation*}
Translation orthogonality and Brezis--Lieb splitting give
\begin{equation*}
 \mathcal E_a(u_k^F)=2\mathcal E_a(W_F)+o(1),\qquad
 \mathcal M(u_k^F)=2\mathcal M(W_F)+o(1),\qquad
 I_a(u_k^F)\longrightarrow2c_F.
\end{equation*}
Write \(f(t)=|t|^{q-2}t\).  The pointwise estimate
\begin{equation*}
 |f(t)-f(v)|\leq C_q(|t|+|v|)^{q-2}|t-v|
\end{equation*}
and H\"older's inequality give
\begin{equation*}
 \|f(v)-f(w)\|_{L^{q/(q-1)}(\omega^{q/Q})}
 \leq C_q\bigl(\|v\|_{L^q(\omega^{q/Q})}
                +\|w\|_{L^q(\omega^{q/Q})}\bigr)^{q-2}
             \|v-w\|_{L^q(\omega^{q/Q})}.
\end{equation*}
Choose \(V_\ell\in C_c^\infty\) with \(V_\ell\to W_F\) in \(H_a\)
and hence in \(L^q(\omega^{q/Q})\).
For fixed \(\ell\), the supports of \(T_{b_k}V_\ell\) and
\(T_{-b_k}V_\ell\) are disjoint for large \(k\).
Their nonlinear interaction vanishes.  The preceding bound gives,
uniformly in \(k\),
\begin{align*}
 &\limsup_{k\to\infty}\bigl\||u_k^F|^{q-2}u_k^F
 -|T_{b_k}W_F|^{q-2}T_{b_k}W_F
 -|T_{-b_k}W_F|^{q-2}T_{-b_k}W_F\bigr\|_
 {L^{q/(q-1)}(\omega^{q/Q})}\\
 &\qquad\leq C\|W_F-V_\ell\|_{L^q(\omega^{q/Q})}.
\end{align*}
Let \(\ell\to\infty\).  Criticality of both translates and the
weighted Sobolev inequality give
\begin{equation*}
 \|I_a'(u_k^F)\|_{H_a^*}
 \leq C\bigl\|f(u_k^F)-f(T_{b_k}W_F)-f(T_{-b_k}W_F)\bigr\|_
 {L^{q/(q-1)}(\omega^{q/Q})}\longrightarrow0.
\end{equation*}
For any \(\beta_k\), at most one of \(b_k-\beta_k\) and
\(-b_k-\beta_k\) remains bounded.  Lemma
\ref{lem:relative-frames} gives, after a subsequence,
\begin{equation*}
 T_{-\beta_k}u_k^F\rightharpoonup W,\qquad
 \mathcal E_a(W)\in\{0,\mathcal E_a(W_F)\},\qquad
 \mathcal E_a(T_{-\beta_k}u_k^F)\longrightarrow2\mathcal E_a(W_F).
\end{equation*}
Thus no such translated sequence converges strongly in \(H_a\).

Let \(r_k\downarrow0\).  Choose \(L_k\to\infty\) slowly enough that
\begin{equation*}
 V_k=\chi_{L_k}W_G,\qquad
 \|V_k-W_G\|_{D_G^{1,2}}+\|V_k-W_G\|_q\longrightarrow0,
 \qquad r_kL_k\longrightarrow0.
\end{equation*}
By Lemma \ref{lem:homogeneous-cutoff}, first fix a radius with small
model error, then choose \(k\) so that \(r_kL\) is small.
Sobolev--H\"older gives
\begin{equation*}
 r_k\|V_k\|_2\leq Cr_kL_k\|W_G\|_q\longrightarrow0.
\end{equation*}
Set \(u_k^G=D_{r_k,0}V_k\).  Continuity of the model derivative
and \(r_kL_k\to0\) give
\begin{equation*}
 \|I_G'(V_k)\|_{(D_G^{1,2})^*}\longrightarrow0,
 \qquad
 \varepsilon_k:=\sup_{|\xi|\leq2L_k}
 |(1+r_k^2|\xi|^2)^s-1|\longrightarrow0.
\end{equation*}
For \(h\in H_a\), put \(\widetilde h=D_{r_k,0}^{-1}h\).
Homogeneous scaling and Sobolev give
\begin{equation*}
 \|\widetilde h\|_{D_G^{1,2}}\leq\|h\|_{H_a},\qquad
 \|\widetilde h\|_q\leq C\|h\|_{H_a},\qquad
 \|h\|_2\leq a^{-1/2}\|h\|_{H_a}.
\end{equation*}
The pulled-back derivative is
\begin{align*}
 \langle I_a'(u_k^G),h\rangle
 &=\langle I_G'(V_k),\widetilde h\rangle
   +ar_k^2\int V_k\widetilde h\\
 &\quad-\int\bigl((1+r_k^2|\xi|^2)^s-1\bigr)
             |V_k|^{q-2}V_k\widetilde h.
\end{align*}
Consequently,
\begin{align*}
 |\langle I_a'(u_k^G),h\rangle|
 &\leq\|I_G'(V_k)\|_{(D_G^{1,2})^*}
       \|\widetilde h\|_{D_G^{1,2}}
       +ar_k\|V_k\|_2\|h\|_2\\
 &\quad+\varepsilon_k\|V_k\|_q^{q-1}\|\widetilde h\|_q
 \leq o(1)\|h\|_{H_a}.
\end{align*}
Also,
\begin{align*}
 \mathcal E_a(u_k^G)&=\mathcal E_G(V_k)+ar_k^2\|V_k\|_2^2
          \longrightarrow\mathcal E_G(W_G),\\
 \left|\mathcal M(u_k^G)-\int|V_k|^q\right|
 &\leq\varepsilon_k\|V_k\|_q^q\longrightarrow0.
\end{align*}
Thus
\begin{equation*}
 I_a(u_k^G)\longrightarrow c_G,\qquad
 I_a'(u_k^G)\longrightarrow0\quad\hbox{in }H_a^*.
\end{equation*}
For every \(\beta_k\),
\begin{equation*}
 \|T_{-\beta_k}u_k^G\|_2\longrightarrow0,\qquad
 \mathcal E_a(T_{-\beta_k}u_k^G)
 \longrightarrow\mathcal E_G(W_G)>0.
\end{equation*}
A strong \(H_a\) limit would be zero by the first limit, contradicting
the second.

Finally, let \(X_k=R_ke\), where \(R_k\to\infty\), and put
\begin{equation*}
 u_k^E=A_{X_k,0}W_E.
\end{equation*}
For \(R=R_k\), define
\begin{equation*}
 b_R(\xi)=\frac{|Re+\xi|^2}{R^2},\qquad
 w_R(\xi)=\left(R^{-2}+|e+\xi/R|^2\right)^{n/(Q-2)},\qquad
 h_R=A_{Re,0}^{-1}h.
\end{equation*}
Both \(b_R\) and \(w_R\) tend pointwise to one.  For \(R\geq1\),
\begin{equation*}
 b_R(\xi)\leq C(1+|\xi|^2),\qquad
 w_R(\xi)\leq C(1+|\xi|^{2s}).
\end{equation*}
Exponential decay of \(W_E\) and its derivatives gives
\begin{align*}
 \mathcal E_a(u_k^E)
 &=\int\bigl(|\nabla_\xi W_E|^2+b_R|\nabla_\eta W_E|^2
                  +aW_E^2\bigr)\longrightarrow\mathcal E_E^a(W_E),\\
 \mathcal M(u_k^E)&=\int w_R|W_E|^q\longrightarrow\int|W_E|^q.
\end{align*}
For \(h\in C_c^\infty\), integration by parts and the equation for
\(W_E\) give
\begin{equation*}
 \langle I_a'(u_k^E),h\rangle
 =\int\left[(1-b_R)\Delta_\eta W_E
             +(1-w_R)W_E^{q-1}\right]h_R.
\end{equation*}
The bracketed source belongs to \(L^2\).  Dominated convergence gives
\begin{equation*}
 \|(1-b_R)\Delta_\eta W_E\|_2
 +\|(1-w_R)W_E^{q-1}\|_2\longrightarrow0.
\end{equation*}
Since \(\|h_R\|_2=\|h\|_2\leq a^{-1/2}\|h\|_{H_a}\), the bound
extends by density from \(C_c^\infty\) to \(H_a\), and
\begin{equation*}
 \|I_a'(u_k^E)\|_{H_a^*}
 \leq a^{-1/2}\bigl(
 \|(1-b_R)\Delta_\eta W_E\|_2
 +\|(1-w_R)W_E^{q-1}\|_2\bigr)\longrightarrow0.
\end{equation*}
\begin{equation*}
 I_a(u_k^E)\longrightarrow c_E,\qquad
 I_a'(u_k^E)\longrightarrow0\quad\hbox{in }H_a^*.
\end{equation*}
For every fixed \(M>0\) and every \(\beta_k\),
\begin{equation*}
 \int_{|x|\leq M}|T_{-\beta_k}u_k^E|^2
 =\int_{|Re+\xi|\leq M}|W_E|^2\longrightarrow0.
\end{equation*}
Any strong \(H_a\) limit would be zero, contradicting
\(\mathcal E_a(u_k^E)\to\mathcal E_E^a(W_E)>0\).
\end{proof}

\subsection{A qualitative deficit estimate}

\begin{corollary}
\label{cor:deficit-intro}
Let \(a\in\mathcal B\), let \(\mathcal F_a\) be the normalized minimizer
set, and put
\(d_a(u)=\inf_{U\in\mathcal F_a}\|u-U\|_{H_a}\).
For every \(\delta>0\),
\begin{equation}
 \rho_a(\delta)=\inf\{\mathcal E_a(u)-S_a:
 \mathcal M(u)=1,\ d_a(u)\geq\delta\}>0,
 \label{eq:deficit-modulus}
\end{equation}
with \(\inf\varnothing=+\infty\).
\end{corollary}

\begin{proof}[Proof of \eqref{eq:deficit-modulus}]
Suppose that \(\rho_a(\delta)=0\) for some \(\delta>0\).
Choose \(u_k\in H_a\) with
\begin{equation*}
 \mathcal M(u_k)=1,\qquad
 \mathcal E_a(u_k)\longrightarrow S_a,\qquad d_a(u_k)\geq\delta.
\end{equation*}
Since \(a\in\mathcal B\), Theorem \ref{thm:threshold-intro} gives,
after a subsequence, \(b_k\in\mathbb R^n\) and
\(U\in\mathcal F_a\) such that
\begin{equation*}
 \|u_k-T_{b_k}U\|_{H_a}\longrightarrow0.
\end{equation*}
But \(T_{b_k}U\in\mathcal F_a\), so
\begin{equation*}
 0<\delta\leq d_a(u_k)
 \leq\|u_k-T_{b_k}U\|_{H_a}\longrightarrow0,
\end{equation*}
a contradiction.  This proves \eqref{eq:deficit-modulus}.
For \(\mathcal M(u)=1\) and \(d_a(u)>0\), the definition of
\(\rho_a\) gives
\begin{equation*}
 \mathcal E_a(u)-S_a\geq\rho_a(d_a(u)/2).
\end{equation*}
No rate for \(\rho_a(\delta)\) as \(\delta\downarrow0\) is obtained.
\end{proof}

\section{Coefficients and weights with the same tangent limits}
\label{sec:perturbations}

Let \(b,h\in C(\mathbb R^m)\), let \(V\in L^\infty(\mathbb R^m)\), and
suppose that, for some positive constants \(c_0,C_0,v_0\),
\begin{align}
 c_0|x|^2\leq b(x)&\leq C_0|x|^2,\qquad
 v_0\leq V(x)\leq C_0,\notag\\
 c_0(1+|x|^2)^{n/(Q-2)}\leq h(x)
 &\leq C_0(1+|x|^2)^{n/(Q-2)}.
 \label{eq:perturbation-comparability}
\end{align}
Assume, as \(x\to0\),
\begin{equation*}
 \frac{b(x)}{|x|^2}\longrightarrow1,\qquad h(x)\longrightarrow1,
\end{equation*}
and, as \(|x|\to\infty\),
\begin{equation}
 \frac{b(x)}{|x|^2}\longrightarrow1,\qquad
 \frac{h(x)}{(1+|x|^2)^{n/(Q-2)}}\longrightarrow1,\qquad
 V(x)\longrightarrow a_\infty>0.
 \label{eq:perturbation-tangents}
\end{equation}
All limits at infinity are uniform outside sufficiently large balls;
the limit for \(V\) is understood in essential supremum.
Define
\begin{equation*}
 \mathcal E_{b,V}(u)=\int
 (|\nabla_xu|^2+b(x)|\nabla_yu|^2+V(x)|u|^2),\qquad
 \mathcal M_h(u)=\int h(x)|u|^q,
\end{equation*}
and let \(S_{b,V,h}\) be the corresponding sharp quotient infimum.
The energy completion is the same underlying space as \(H_1\).

\begin{theorem}
\label{thm:perturbed-decomposition}
Under \eqref{eq:perturbation-comparability}--\eqref{eq:perturbation-tangents}
and the limits at zero above, Theorem \ref{thm:profile-intro} holds with
finite model \((\mathcal E_{b,V},\mathcal M_h)\), singular model
\((\mathcal E_G,\mathcal M_G)\), and escape model
\((\mathcal E_E^{a_\infty},\mathcal M_E)\).  The same three gauges give a
refined inequality and the stated power remainder bound.  The
constants may depend on \(m,n,c_0,C_0,v_0,a_\infty\); they are uniform
within classes with these quantities fixed.
Furthermore,
\begin{equation*}
 0<S_{b,V,h}\leq\min\{S_G,S_E(a_\infty)\}.
\end{equation*}
Strict inequality below both tangent constants implies attainment and
compactness of normalized minimizing sequences modulo \(y\)-translations.
At equality, every normalized minimizing sequence has one profile of
full mass of an admissible finite, singular, or escape type.
Every Palais--Smale sequence for
\(\mathcal E_{b,V}/2-\mathcal M_h/q\) has the finite critical-profile
decomposition and strong remainder of Theorem \ref{thm:ps-global-intro}.
\end{theorem}

\begin{proof}
Condition \eqref{eq:perturbation-comparability} gives constants
\(c,C>0\) such that
\begin{equation*}
 c\mathcal E_1(u)\leq\mathcal E_{b,V}(u)\leq C\mathcal E_1(u),
 \qquad c_0\mathcal M(u)\leq\mathcal M_h(u)\leq C_0\mathcal M(u).
\end{equation*}
Proposition \ref{prop:refined} therefore gives
\begin{equation*}
 \mathcal M_h(u)\leq C\mathcal E_{b,V}(u)^{q/2-\gamma}
 \max\{\Theta_F(u),\Theta_G(u),\Theta_E(u)\}^{\gamma}.
\end{equation*}
The same comparison yields
\begin{equation*}
 \max_\tau\Theta_\tau(u)\leq C\mathcal E_{b,V}(u),\qquad
 |\Theta_\tau(u)^{1/2}-\Theta_\tau(v)^{1/2}|
 \leq C\mathcal E_{b,V}(u-v)^{1/2}.
\end{equation*}
Lemma \ref{lem:uniform-local-rellich} gives local strong \(L^2\)
convergence of the pullbacks.  Norm comparison also transfers the inverse
estimate: a gauge of size \(\theta>0\) yields a profile with
\(\mathcal E_\tau(p)\geq c\theta\).

Fix a compact model set \(K\times K'\).
In a singular frame, the energy coefficients are
\begin{equation*}
 1,\qquad \frac{b(r_k\xi)}{r_k^2},\qquad r_k^2V(r_k\xi).
\end{equation*}
There are \(\varepsilon_k\to0\), depending on \(K\), such that
\begin{equation*}
 \left|\frac{b(r_k\xi)}{r_k^2}-|\xi|^2\right|
 \leq\varepsilon_k|\xi|^2,\qquad
 r_k^2\|V(r_k\cdot)\|_\infty\leq C_0r_k^2,\qquad
 \|h(r_k\cdot)-1\|_{L^\infty(K)}\longrightarrow0.
\end{equation*}
For a bounded singular pullback \(v_k\) and a compact model test
\(p_L\), the quadratic pairing error is bounded by
\begin{equation*}
 \varepsilon_k\mathcal E_G(v_k)^{1/2}\mathcal E_G(p_L)^{1/2}
 +C_0r_k^2\|v_k\|_{L^2(K\times K')}\|p_L\|_2
 \longrightarrow0.
\end{equation*}
Thus its limiting energy and mass are \(\mathcal E_G\) and
\(\mathcal M_G\).

In an escape frame, put \(R_k=|X_k|\).  The energy coefficients and
mass density are
\begin{equation*}
 1,\qquad \frac{b(X_k+\xi)}{R_k^2},\qquad V(X_k+\xi),
 \qquad\rho_k(\xi)=R_k^{-nq/2+n}h(X_k+\xi).
\end{equation*}
The limits at infinity give
\begin{align*}
 \left\|\frac{b(X_k+\cdot)}{R_k^2}-1\right\|_{L^\infty(K)}
 +\|\rho_k-1\|_{L^\infty(K)}&\longrightarrow0,\\
 \|V(X_k+\cdot)-a_\infty\|_{L^\infty(K)}&\longrightarrow0.
\end{align*}
The norm involving \(V\) is understood as an essential supremum.
The error in a compact quadratic pairing is consequently bounded by
\begin{equation*}
 \varepsilon_k\|v_k\|_{H^1(K\times K')}
                  \|p_L\|_{H^1(K\times K')}\longrightarrow0.
\end{equation*}
This gives \(\mathcal E_E^{a_\infty}\) and \(\mathcal M_E\).
A finite frame is an isometry for \(\mathcal E_{b,V}\) and
\(\mathcal M_h\), since the coefficients depend only on \(x\).

In each frame, the mass density converges uniformly on compact sets
to its positive model density.  The pulled-back sequence is bounded
in local \(L^q\) and converges almost everywhere along a subsequence.
For a fixed compact core \(p_L\), Brezis--Lieb splitting gives
\begin{equation*}
 \mathcal M_h(u_k)-\mathcal M_h(u_k-P_{k,L}p)
 =\mathcal M_\tau(p)-\mathcal M_\tau(p-p_L)+o_k(1).
\end{equation*}
To replace the density by its limit, use
\begin{equation*}
 \left|\int_{K\times K'}(\rho_k-\rho_\tau)
       (|v_k|^q-|v_k-p_L|^q)\right|
 \leq C\|\rho_k-\rho_\tau\|_{L^\infty(K)}
       \bigl(\|v_k\|_{L^q(K\times K')}^q+\|p_L\|_q^q\bigr)
 \longrightarrow0.
\end{equation*}
The quadratic pairing limits give energy subtraction and lower
semicontinuity of the model energies.

The relative coordinate maps are unchanged.  As in Section 4, support
separation handles separated positions.  Absolute continuity of the
fixed profile's energy and mass integrals handles finite--singular pairs
and separated singular scales.  The coefficient limits above identify
the singular and escape pairings.
For a general profile \(p\), choose compact \(p_L\) with
\(p_L\to p\) in its model space.  Norm comparison gives
\begin{equation*}
 \limsup_k\mathcal E_{b,V}(R_kp-P_{k,L}p)
 \leq C\mathcal E_\tau(p-p_L)\longrightarrow0.
\end{equation*}
The global embedding makes this error tend to zero in
\(\mathcal M_h^{1/q}\), so decoupling extends to full profiles.
Apply the diagonal choices in Lemmas
\ref{lem:synchronized-realizations} and \ref{lem:one-profile-subtraction}
to these norm, pairing, mass and weak-test limits.
Theorem \ref{thm:three-family-criterion} gives the decomposition and
the power estimate for the remainder.

For nonzero compact smooth \(\phi,\psi\), the singular and escape
recovery quotients give
\begin{equation*}
 S_{b,V,h}\leq\frac{\mathcal E_G(\phi)}{\|\phi\|_q^2},\qquad
 S_{b,V,h}\leq\frac{\mathcal E_E^{a_\infty}(\psi)}{\|\psi\|_q^2}.
\end{equation*}
Take the two infima.  Positivity follows from norm comparison and
the global embedding.
For a normalized minimizing sequence, let \(m_j\) be its profile
masses.  The splittings and the vanishing remainder give
\begin{equation*}
 \sum_jm_j=1,\qquad
 S_{b,V,h}\geq\sum_jS_{\tau_j}m_j^{2/q}
 \geq S_{b,V,h}\sum_jm_j^{2/q}\geq S_{b,V,h}.
\end{equation*}
As in Theorem \ref{thm:threshold-intro}, strict concavity leaves one
profile of full mass and a strong energy remainder.
Strict inequality below the tangent constants excludes the singular
and escape types, so the finite profile attains \(S_{b,V,h}\).

Put \(\mathcal I=\mathcal E_{b,V}/2-\mathcal M_h/q\).
For a Palais--Smale sequence,
\begin{equation*}
 \left(\frac12-\frac1q\right)\mathcal E_{b,V}(u_k)
 =\mathcal I(u_k)-\frac1q\langle\mathcal I'(u_k),u_k\rangle
 \leq C+o(1)\mathcal E_{b,V}(u_k)^{1/2}.
\end{equation*}
Thus it is bounded, and
\(\mathcal E_{b,V}(u_k)-\mathcal M_h(u_k)=o(1)\).
Testing on compact model sets, with the pairing and density limits
above, gives the model Euler equations.  The density argument in
Lemma \ref{lem:ps-profile-criticality} extends them to each model
space.  Every nonzero critical profile satisfies
\begin{equation*}
 \mathcal E_\tau(p)=\mathcal M_\tau(p)
 \geq S_\tau^{q/(q-2)}>0.
\end{equation*}
As in Theorem \ref{thm:ps-global-intro}, only finitely many profiles
occur.  Their terminal mass vanishes.  Subtracting the energy and mass
identities gives
\begin{equation*}
 \mathcal E_{b,V}(w_k)=\mathcal M_h(w_k)+o(1)\longrightarrow0.
\end{equation*}
\end{proof}

\begin{corollary}
\label{cor:binding-perturbations}
Fix \(a\in\mathcal B\).  Let \((b_j,V_j,h_j)\) satisfy the hypotheses of
Theorem \ref{thm:perturbed-decomposition}, with \(a_\infty=a\), and suppose
\begin{equation}
 \sup_{x\ne0}\left|\frac{b_j(x)}{|x|^2}-1\right|
 +\|V_j-a\|_\infty
 +\sup_x\left|\frac{h_j(x)}{(1+|x|^2)^{n/(Q-2)}}-1\right|
 \longrightarrow0.
 \label{eq:small-perturbations}
\end{equation}
Then \(S_{b_j,V_j,h_j}\to S_a\), and these perturbed constants are
strictly below both tangent constants for large \(j\).
Their normalized minimizers are compact in \(H_a\) modulo
\(y\)-translations, with every subsequential limit in \(\mathcal F_a\).
This applies, in particular, to small changes supported in a fixed
annulus \(0<r<|x|<R\).  It gives an attained family near each explicit
parameter in Theorem \ref{thm:explicit-binding}.
\end{corollary}

\begin{proof}
By \eqref{eq:small-perturbations}, there are \(\delta_j\to0\) such that,
for all \(u\in H_a\),
\begin{align*}
 (1-\delta_j)\mathcal E_a(u)
 &\leq\mathcal E_{b_j,V_j}(u)\leq(1+\delta_j)\mathcal E_a(u),\\
 (1-\delta_j)\mathcal M(u)
 &\leq\mathcal M_{h_j}(u)\leq(1+\delta_j)\mathcal M(u).
\end{align*}
For \(\delta_j<1\), take infima to obtain
\begin{equation*}
 \frac{1-\delta_j}{(1+\delta_j)^{2/q}}S_a
 \leq S_{b_j,V_j,h_j}
 \leq\frac{1+\delta_j}{(1-\delta_j)^{2/q}}S_a.
\end{equation*}
Hence \(S_{b_j,V_j,h_j}\to S_a\).
The tangent constants are fixed, and
\(S_a<\min\{S_G,S_E(a)\}\).  Strict binding holds for large \(j\).

Let \(U_j\) be a minimizer with \(\mathcal M_{h_j}(U_j)=1\).
The same estimates give
\begin{equation*}
 \frac1{1+\delta_j}\leq\mathcal M(U_j)\leq\frac1{1-\delta_j},\qquad
 \frac{S_{b_j,V_j,h_j}}{1+\delta_j}
 \leq\mathcal E_a(U_j)
 \leq\frac{S_{b_j,V_j,h_j}}{1-\delta_j}.
\end{equation*}
Set \(v_j=U_j/\mathcal M(U_j)^{1/q}\).  Then
\begin{equation*}
 \mathcal M(v_j)=1,\qquad\mathcal E_a(v_j)\longrightarrow S_a,
 \qquad
 \|v_j-U_j\|_{H_a}
 =|\mathcal M(U_j)^{-1/q}-1|\mathcal E_a(U_j)^{1/2}
 \longrightarrow0.
\end{equation*}
Apply Theorem \ref{thm:threshold-intro} to \((v_j)\).
Its translated subsequential limits belong to \(\mathcal F_a\).
The last norm limit gives the same conclusion for \((U_j)\).
\end{proof}

\begin{remark}
The interval and endpoint assertions in Theorem
\ref{thm:connected-binding-intro} concern the prototype coefficients
\(b(x)=|x|^2\), \(V(x)=a\), and
\(h(x)=(1+|x|^2)^{n/(Q-2)}\).  Their proof uses the exact parameter
rescaling \eqref{eq:parameter-rescaling}.  Theorem
\ref{thm:perturbed-decomposition} gives the profile, minimizing-sequence
and Palais--Smale conclusions for the general coefficient class.
\end{remark}

\appendix
\section{Cutoffs and density arguments}
\label{sec:cutoffs}

Let $\rho(x,y)=(|x|^4+4|y|^2)^{1/4}$.  Choose
$\chi\in C_c^\infty([0,\infty))$ with $\chi=1$ on $[0,1]$ and $\chi=0$ on
$[2,\infty)$.  Put
\begin{equation*}
 \chi_R(x,y)=\chi(\rho(x,y)/R).
\end{equation*}
Then
\begin{equation}
 |\nabla_G\chi_R|\leq CR^{-1}
 \mathbf1_{\{R<\rho<2R\}}.
 \label{eq:cutoff-gradient-bound}
\end{equation}

\begin{lemma}
\label{lem:homogeneous-cutoff}
For every $u\in D_G^{1,2}$,
\begin{equation*}
 \chi_Ru\longrightarrow u\quad\text{in }D_G^{1,2}
 \quad(R\to\infty).
\end{equation*}
\end{lemma}

\begin{proof}
For smooth \(v\), the product rule, H\"older's inequality and
\eqref{eq:cutoff-gradient-bound} give
\begin{equation*}
 \mathcal E_G(\chi_Rv)
 \leq C\mathcal E_G(v)
 +CR^{-2}|\{R<\rho<2R\}|^{1-2/q}\|v\|_q^2
 \leq C\mathcal E_G(v).
\end{equation*}
Here \(|\{R<\rho<2R\}|\asymp R^Q\),
\(Q(1-2/q)=2\), and the homogeneous Sobolev inequality was used.
Choose \(v_\ell\in C_c^\infty\) with \(v_\ell\to u\) in
\(D_G^{1,2}\).  Apply the bound to \(v_\ell-v_h\).  Then
\((\chi_Rv_\ell)\) is Cauchy in \(D_G^{1,2}\) and converges to
\(\chi_Ru\) in \(L^q\).  Its limit belongs to \(D_G^{1,2}\), and
the product rule passes to the limit:
\begin{equation*}
 \nabla_G(\chi_Ru-u)
 =(\chi_R-1)\nabla_Gu+u\nabla_G\chi_R.
\end{equation*}
The first term tends to zero in \(L^2\) by dominated convergence.
For the second term,
\begin{align*}
 \int|u|^2|\nabla_G\chi_R|^2
 &\leq CR^{-2}
 \left(\int_{R<\rho<2R}|u|^q\right)^{2/q}
 |\{R<\rho<2R\}|^{1-2/q}\\
 &\leq C\left(\int_{R<\rho<2R}|u|^q\right)^{2/q}
 \longrightarrow0.
\end{align*}
\end{proof}

\begin{lemma}
\label{lem:weighted-bl}
Let $f_k\to f$ almost everywhere and suppose
\begin{equation}
 \sup_k\int\omega^{q/Q}|f_k|^q<\infty.
 \label{eq:weighted-bl-bound}
\end{equation}
Then
\begin{equation*}
 \int\omega^{q/Q}
 \bigl(|f_k|^q-|f_k-f|^q-|f|^q\bigr)\longrightarrow0.
\end{equation*}
\end{lemma}

\begin{proof}
Put \(g_k=\omega^{1/Q}f_k\) and \(g=\omega^{1/Q}f\).
Then \(g_k\to g\) almost everywhere, and
\eqref{eq:weighted-bl-bound} gives
\begin{equation*}
 \sup_k\|g_k\|_q^q
 =\sup_k\int\omega^{q/Q}|f_k|^q<\infty,
 \qquad \|g\|_q^q\leq\liminf_k\|g_k\|_q^q<\infty.
\end{equation*}
The Brezis--Lieb lemma \cite{BrezisLieb1983} yields
\begin{equation*}
 \int\bigl(|g_k|^q-|g_k-g|^q-|g|^q\bigr)\longrightarrow0.
\end{equation*}
Since \(|g_k-g|^q=\omega^{q/Q}|f_k-f|^q\), this is the asserted
weighted identity.
\end{proof}

\end{document}